\documentclass[11pt]{article}
\usepackage{amsmath,amssymb,amscd,mathrsfs}

\newtheorem{propo}{Proposition}[section]
\newtheorem{defi}[propo]{Definition}

\newtheorem{lemma}[propo]{Lemma}
\newtheorem{corol}[propo]{Corollary}

\newtheorem{theo}[propo]{Theorem}

\newcommand{\ld}{,\ldots ,}

\newcommand{\ra}{ \rightarrow }

\newcommand{\lan}{ \langle }
\newcommand{\ran}{ \rangle }

\newcommand{\diag}{\mathop{\rm diag}\nolimits}

\newcommand{\Id}{\mathop{\rm Id}\nolimits}
\newcommand{\Irr}{\mathop{\rm Irr}\nolimits}

\newcommand{\al}{\alpha}

\newcommand{\ep}{\varepsilon}

\newcommand{\lam}{\lambda }

\newcommand{\up}{^{-1}}

\newcommand{\si}{\sigma }
\newcommand{\om}{\omega }

\newcommand{\CG}{{\mathbf G}}

\def\d12{{_{12}}}

\def\au{{automorphism }}

\def\ei{{eigenvalue }}
\def\eis{{eigenvalues }}

\def\f{{following }}

\def\ho{{homomorphism }}
\def\hos{{homomorphisms }}

\def\ii{{if and only if }}

\def\ir{{irreducible }}

\def\irr{{irreducible representation }}

\def\itf{{It follows that }}

\def\mult{{multiplicity }}

\def\rep{{representation }}
\def\reps{{representations }}

 \newcommand{\GC}{{\mathbf G}}
\newcommand{\TC}{{\mathbf T}}
\newcommand{\med}{\medskip}
\newcommand{\bl}{\begin{lemma}\label}
\newcommand{\el}{\end{lemma}Proof. }

\date{}

\begin{document}

\title{On the Steinberg character of an orthogonal group \\ over a finite field}

\author{A. E. Zalesski}
\maketitle

\medskip
\centerline{Dedicated to Roderick Gow on the occasion of his 65th birthday}


\medskip
{\bf Abstract} We determine the \ir constituents of the Steinberg character of an  
orthogonal  group over a finite field restricted to the orthogonal  group of one less dimension.

\section{Introduction}

 The Steinberg character plays outstanding role in general theory of characters of Chevalley groups. It has many interesting features studied in numerious papers. In particular, experts are interested in branching rules for restrictions of the    Steinberg character to natural subgroups. In this paper we study  the branching rule for the special case described above.

Let $F_q$ denote the   field of $q$ elements, and let $V'$ be a vector space 
over $F_q$, endowed with a non-degenerate quadratic
form. Let $H=SO(V')$ be the special orthogonal  group. For a non-degenerate subspace $V$ of $V'$ of codimension 1, set $G=SO(V)$. Let $St_H$ denote the Steinberg character of $H$.  Theorems \ref{thm1} and \ref{thm2} are stated for   $q$  odd, but the  results remain valid for $q$ even if one takes for $G$ the group $\Omega(V)$ formed by the elements of spinor norm 1.

\begin{theo}\label{thm1}   $(1)$ Suppose that $\dim V$ is odd.  Then the restriction of $St_H$ to G is a \mult free character. 


$(2)$ Suppose that $\dim V$ is even. 
Then the multiplicities of the \ir constituents  of $St_H|_G$ does not exceed $2$. 
\end{theo}

We also determine the 
\ir constituens of the restriction $St_H|_G$ of $St_H$ to $G$.  
The nature of the \ir constituents of $St_H|_G$ is described as follows (the meaning of the term "regular character"   is as in the Deligne-Lusztig theory):  

\begin{theo}\label{thm2}   Let $\rho$ be an \ir constituent of $St_H|_G$. Then $\rho$ is either a regular character, or a constituent of the induced character $\si^G$, where $\si$ is an  \ir character of a maximal parabolic subgroup P of G.

More precisely, if $L$ is a  
Levi subgroup of $P$ then $\si$ is trivial on the unipotent radical of $P$ and $\si|_L$ is a regular character of L. 
Furthermore, $L=GL_m(q)\times SO(V_1)$ for some  $m>0$ and
non-degenerate subspace $V_1$ of $V$ (or $V_1=0$).    Then  $\si|_L=St_{GL_m(q)}\otimes \rho'$ for some regular character $\rho'$ of $SO(V_1)$. 
 \end{theo}

 If $G=SO_{2n+1}(q)$ then every regular character of $G$ is a constituent of $St_H|_G$ (Proposition \ref{pp5}).  Note that 
the group $G$ is not quasi-simple if $q$ is odd, and the above results are not valid for the derived subgroup of $G.$

The key ingredient of our consideration is the analysis of the Curtis dual $\om_G$ of
  $St_{H}|_G$ (see \cite[\S 8]{DM}, \cite[\S 71]{CR2} for the duality 
of generalized characters of finite reductive groups). In our situation $\om_G$ is a generalized character, 
whose restrictions to the maximal tori of $G$ determine the \ir constituents of
$St_{H}|_{G}$.  
 
The method used here for proving Theorems \ref{thm1} and \ref{thm2} is analogous to that developed in \cite{HZ09}. There
for $G=Sp_{2n}(q)$ and $ U_n(q)$,  the authors determined the \ir constituents  
of the character $\omega\cdot St$, where $\omega$ is the Weil character of $G$.
Note that, if $q$ is even then $SO_{2n+1}(q)\cong Sp_{2n}(q)$, and this case has been already treated in  \cite{HZ09}.

A priori, the method in  \cite{HZ09} does not seem to work for $G$ orthogonal, as an orthogonal group has no proper analog of the Weil character. The main idea of this paper is  the suggestion of the Curtis dual of $St_H|_G$ for the role  played in \cite{HZ09} by the Weil character. 

This deserves to be explained in more detail. First, for the purpose of \cite{HZ09} the values of the Weil character at the non-semisimple elements are irrelevant, and hence the Weil character can be replaced by any other generalizized character which coincides with the Weil one at the semisimple elements. This could be, for instance,  the Brauer lift of the Weil character. Secondly, it was observed in \cite{HZ2}  that the Brauer lift of $\om$ coincides with the Curtis dual of $\om\cdot St$. Thirdly, in the appendix by Brunat \cite{OB} in  \cite{HZ09}, it was shown that the Steinberg character of $U_{n+1}(q)$ restricted to $U_n(q)$ coincides with $\om\cdot St$.

These observations may lead to a hint that the Curtis dual of $St_H|_G$ could be a promissing replacement of the Weil character in \cite{HZ09}. However, I wish to admit here the significance of an unpublished manuscript by Gow \cite{Gw}  who 
suggested the difference of certain permutation characters of an orthogonal group for an analog of the Weil character. 
(The character in question was introduced in \cite[p. 413]{Gt} for a  different purpose.)
I have found out that the Gow character coincides at the semisimple elements with the Curtis dual of 
$St_H|_G$. 

Technically, the reasoning in \cite{HZ09} uses quite heavily certain very specific properties of the Weil character. A priori,
the Curtis dual $\om_G$ of $St_H|_G$ may not enjoy such properties. However, rather surprizingly, it does have a lot of properties analogous to those of the Weil character. The \f "multiplication theorem" is one of the most useful properties of  $\om_G$. 

\begin{theo}\label{mm8}   Let $V$  be an orthogonal space, and  $V=V_1\oplus V_2$, where  $V_1,V_2$ are non-degenerate. Let $G=SO(V)$, $G_i=SO(V_i)$ for $i=1,2$, and  let $\om_G,\om_{G_i}$ be the generalized characters as defined above. 

\med
$(i)$ Suppose that at least one of $\dim V_1,\dim V_2$ is even.
Then 
$\om_G(g)=  \om_{G_1}(g)\cdot
\om_{G_2}(g)$ for every  $g\in G_1G_2$.

\med
$(ii)$ Suppose that both $\dim V_1,\dim V_2$ are odd. Then  
$\om_G(g)= q\cdot \om_{G_1}(g)\cdot \om_{G_2}(g)$ for every 
$g\in G_1G_2$. \end{theo}

\medskip
The branching rule for the restriction of the cuspidal regular  \ir \reps $\phi$ of $H=SO(V')$ to $G=SO(V)$, where $q$ and $\dim V$ are  odd, 
is obtained in Reeder \cite{R}. The case where $\phi=St_H$ is somehow at the opposite extreme, as $St_H$ is cuspidal only if $H$
is abelian. It is stated in \cite[p. 573]{R} that one can deduce from unpublished work of Bernstein and Rallis that
the restriction $\phi|_G$  is \mult free for every \irr $\phi$ of $H$  provided $\dim V$ is odd.  


\medskip
{\bf Notation}. $F_q$ is the field of $q$ elements, where $q$ is a power of a prime $p$, The algebraic closure of $F_q$ is denoted by  $\overline{F}_q$.  All vector spaces under consideration are over $F_q$, unless otherwise is said explicitly.  Let $V$ be a  vector space over $F_q$. We set  $\overline{V}=V\otimes \overline{F}_q$. 
 By 
${\rm End}\,V$ we denote the ring of all $F_q$-\hos $V\ra V$. If $g\in {\rm End}\,V$ then $V^g$ is the subspace 
of vectors fixed by $g$ (unless otherwise is stated).  The general and special linear groups are denoted by $GL(V)$ and $SL(V)$, respectively. If $n=\dim V$ then we also use the notation $GL_n(q)$ for $GL(V)$  and $SL_n(q)$ for $SL(V)$. 
If $A,B,C$ are square matrices, we denote by $\diag(A,B,C)$ the block-diagonal matrix with consequtive blocks $A,B,C$
(similarly, for any number of blocks).   

We say that $V$ is an orthogonal space if $V$ is endowed with a quadratic form $Q$, say,  which is non-degenerate unless 
$q$ is even and $\dim V$ is odd; in the latter case the form $Q$ is assumed to be non-defective \cite{Die}. 
If $q$ is odd then the quadratic form can be replaced by the associated bilinear form. 
It is well known that $\overline{V}$ can be endowed with a  quadratic form $\overline{Q}$ such that 
$Q$ is the restriction of $\overline{Q}$ to $V$. 

Let  $V,V'$ be orthogonal spaces over $F_q$ defined by quadratic forms   $Q,Q'$, respectively. 
An embedding $e:V\ra V'$ is called {\it natural} if $Q$ is the restriction of $Q'$ to $e(V)$. 

A subspace $U$ of $V$ is called totally singular if $Q$ vanishes on $U$, and non-degenerate if the restriction of $Q$ to  $U$
yields a non-denenerate  associated bilinear form. 

Let $V$ be an orthogonal space. We say that $V$ is of Witt defect $d$  if $\dim V=2n$ is even and  maximal totally singular subspaces of  $V$ are of dimension $n-d$. Note that $n-d$ is called the Witt index of  $V$.

The isometry group of an orthogonal space  $V$ is denoted by $O(V)$, and we set  $SO(V)=O(V)\cap SL(V)$.
The subgroup of $SO(V)$ of elements of spinor norm 1 is denoted by $\Omega(V)$ \cite{Die, Gr}. Note that  $|SO(V):\Omega(V)|=2$. If $d=\dim V$ is odd, the group $SO(V)$ is also  denoted by $SO_d(q)$. If $d$ is even, there are two non-isomorphic non-degenerate orthogonal spaces. If $V$ is of Witt defect 0, we often  use $SO^+_d(q)$ for $SO(V)$, and  if the Witt defect equals 1 then we write $SO_d^-(q)$ for  $SO(V). $

Groups $G=SO(V)$ for $n>2 $ are groups with a split BN-pair, for which
one defines an irreducible character
called the Steinberg character, see \cite{CR2}; this is denoted here
by $St_{G}$. For uniformity, if $\dim V\leq 2$, we define $St_G=1_G$.
 The   Steinberg-plus character $St_G^+$ is defined in Section 4.

  If $X$ is a finite group, we denote by $\Irr X$ the set of \ir characters, by $\rho_X^{reg}$ the character of the regular \rep of $X$ and by $1_X$  the trivial character. For a subgroup  $Y$ of $X$ and a character (or a representation) 
 $\mu$ 
of $X$ the symbol $\mu|_Y$ denotes the restriction of $\mu$ to $Y$. If
$\lam $ is a character of $Y$, we write $\lam^G$ for the 
induced character. For the notion of Harish-Chandra induction a reader may consult \cite[\S  70]{CR2}. Note that we use 
$\otimes$ to express a \ir character of a direct product
of two groups in terms of characters of the multiple.

\medskip

Let $\GC$ be a reductive connected algebraic group over an
algebraically closed  field of characteristic $p>0$.
An algebraic group endomorphism $Fr:\GC\ra \GC$ is called  {\it Frobenius} if its fixed
point group $G := {\GC}^{Fr}$ is finite. A group $G$ is called a
{\it finite reductive group} or a group of Lie type in characteristic
$p$ if there exists a reductive connected algebraic group $\CG$
over a field of characteristic $p$ and a Frobenius endomorphism
$Fr:\GC\ra \GC$ such that $G=\CG^{Fr}$. A subgroup $T$ of $G$ is
called a maximal torus if there exists an $Fr$-stable maximal torus
${\mathbf T}$ of $\CG$ such that $T={\mathbf T}^{Fr}$. Note that saying that two maximal tori $T,T'$ of $G$ are $G$-conjugate means that the respective $Fr$-stable maximal tori in ${\mathbf G}$ are $G$-conjugate. For a maximal torus  $T={\mathbf T}^{Fr}$ we   set $W(T)=N_G({\mathbf T})/T$. 

For simple algebraic groups $\GC$ Frobenius  endomorphisms and the groups $G=\GC^{Fr}$ have been classified, see Carter \cite[\S 1.19]{Ca}.  Let ${\mathbf G}$ be the simply connected simple algebraic group of type $D_n$ or $B_n$. We identify ${\mathbf G}$  with the  $\overline{F}_q$-form of ${\mathbf G}$.   Then there is an algebraic group \ho  $\eta:
{\mathbf G}\ra SO(\overline{V})$ with finite kernel. If $q$ is odd, $\eta$ is surjective,  $SO(\overline{V})= \Omega(\overline{V})$ and the group $SO(\overline{V})$ is connected \cite[p. 258]{B}. If $q$ is even then $\eta({\mathbf G})=\Omega (\overline{V})$, so $\Omega (\overline{V})$ is connected.
 Therefore, for $q$ even  $SO(\overline{V})$ is not connected  as    $\Omega (\overline{V})$ has index 2 in $SO(\overline{V})$.
Thus,  if $q$ is odd then $SO(V) $ is a finite reductive group,  if $q$ is even then so is $\Omega (V)$. 

For a connected algebraic group $\GC$ with Frobenius  endomorphism $Fr$ one defines the relative rank $rel.rk$ of 
$\GC$ and the function $\ep_\GC=(-1)^{rel.rk (\GC)}$ called the sign of $G$, see Carter \cite[p. 197-199]{Ca}. If $\GC$ is not connected,
we define the sign to be that of $\GC^0$,  the connected component of $\GC$.  

For the notions of  dual groups $\CG^*$ of $\GC$ and $G^*$
of $G$, see \cite[Ch. 4]{C} and \cite[13.10]{DM}. (Note that
$G^*=(\CG^*)^{Fr^*}$, where $Fr^*$ is a suitable Frobenius
endomorphism of $G^*$. To simplify the notation, we shall use $Fr$
for $Fr^*$, which should not lead to any confusion.)

Notation for Deligne-Lusztig characters is introduced in Section 6. 

\section{Some properties of orthogonal groups}

\begin{lemma}\label{y9}  
Let $G=SO(V)$,
$q$ odd, or $\Omega(V)$, 
$q$ even. Suppose that $\dim V$ is even. Let $s\in  G$ be a semisimple element.
 Then $C_{O(V)}(s)\subset G$ \ii neither $1$ nor $-1$ is an \ei of $s$. 
\end{lemma}
Proof.  Let 
$\overline{V}$ the natural module for $O(2n,\overline{F}_q)$. We can asume that the form defining $G$ is the restriction to $V$ of that defining $O(2n,\overline{F}_q)$.
As $SO(V)=SO(\overline{V})\cap O(V)$ and $ \Omega (V)= \Omega (\overline{V})\cap O(V)$, respectively, if $q$ is odd or even, it follows that it suffices to prove that  $C_{O(\overline{V})}(s)\subseteq  SO(\overline{V})$  and $C_{O(\overline{V})}(s)\subseteq  \Omega (\overline{V})$, respectively. 

The "if" part. Let $\overline{V}=V_1\oplus\cdots \oplus  V_t$ be the decomposition of $\overline{V}$ as a direct sum of the homogeneous components of $s$ on $\overline{V}$. Then the restriction of $s$ to $V_i$ ($i=1\ld t$) is scalar $\al_i\cdot \Id$, say.
Suppose that  $\pm 1$ is not an \ei of $s$. Then each $V_i$ is totally  singular, $t$ is even and $V_1\ld V_t$ can be reordered so that $V_{2i-1}$ and $V_{2i}$ were dual. (That is, the Gram matrix of the bilinear form on $V_{2i-1}\oplus V_{2i}$ under a certain basis is $\begin{pmatrix}0&\Id\\ \Id&0\end{pmatrix}$.) Let $x\in C_{O(\overline{V})}(s)$. Then $xV_i=V_i$ for every $i=1\ld t$. Let $x_i$ denote the restriction of $x$ to $V_i$. Then $x_{2i}={}^tx_{2i-1}^{-1}$, where ${}^tx_{2i-1}$ means the transpose of $x_{2i-1}$.
\itf $\det x=1$, and hence, if $q$ is odd, then $x\in SO(\overline{V})$ as claimed. If $q$ is even then $|O (V):\Omega (V)|=2$, whereas $|GL(V_i):SL(V_i)|$ is odd. So the claim follows.

The "only  if" part. Suppose the contrary. Let $Y\neq 0$ be the $1$- or $-1$-eigenspace of  $s$ on $V$. Then $Y$ is non-degenerate.
Obviously, $C_{O(V)}(s)$ contains a subgroup isomorphic to $O(Y)$, and hence $C_{O(V)}(s)$ contains a matrix of determinant $-1$ if $q$ is odd, or of spinor norm $-1$ if $q$ is even. So $C_{O(V)}(s)$ is not contained in $G$.

\bl{os4} $(1)$ Let  $\GC=Sp_{2n}(\overline{F}_q) $. Then  $C_{{\mathbf G}}(s)$ is connected for any semisimple element $s\in \GC$.

 $(2)$  Let 
$G=SO_{2n}(\overline{F}_q)$ and $s\in G$ a semisimple element. Then $C_{{\mathbf G}}(s)$ is connected \ii either $1$ or $-1$ is not an \ei of $s$. In particular, if $q$ is even then $C_{{\mathbf G}}(s)$ is always connected.       
In addition, if $C_{{\mathbf G}}(s)$ is connected then $C_{{\mathbf L}}(s)$ is connected for every Levi subgroup ${\mathbf L}$
of $\GC$ containing $s$.
\el
(1)  As $Sp_{2n}(\overline{F}_q)$ is simply connected, the statement is a special case of 
\cite[Part E, Ch.2, 3.9]{SS}.

(2)  Let $\overline{V}$   be the natural module for ${\mathbf G}$.  If $\lam$ is an \ei of $s$, let $V(\lam)$ denote the $\lam$-eigenspace of $V$. If $\lam\neq \pm 1$ then $V(\lam)$ is totally   singular, and hence  $V({\lam\up})\neq 0$. Set  $V_{\lam}=V(\lam) +V(\lam\up)$. Then $V_\lam $ is non-degenerate. Moreover, it is well known that 
the common stabilizer of $V(\lam) ,V(\lam\up)$ is isomorphic to $GL(V(\lam))$.
Therefore, if $1$, $-1$ are not \eis of $s$ then $C_{O(V)}(s)$ is 
the direct product of $GL(V(\lam))$, and hence connected. 

Furthermore, suppose that  $1$ or $-1$ (but not each of them) is  an \ei of $s$, and let $Y$ be the $1$- or $-1$-eigenspace of $s$.
Then $C_{O(V)}(s)$ is the direct product of $GL(V(\lam))$ with $\lam\neq \pm 1$ and $O(Y)$, so   $C_{SO(V)}(s)$ is the direct product of $GL(V(\lam))$ with $\lam\neq \pm 1$ and $SO(Y)$. Again, this group is connected.

Finally, suppose that both the $1$- and $-1$-eigenspaces of $s$ are non-zero (so $q$ is odd), and    let $Y,Z$ be the $1$- or $-1$-eigenspaces of $s$, respectively. Then $C_{SO(V)}(s)$ contains a subgroup isomorphic to the direct product of $GL(V(\lam))$ with $\lam\neq \pm 1$, $SO(Y)$ and $SO(Z)$. This group is connected. However, $C_{SO(V)}(s)$ additionally contains all elements $g_1g_2$,
where $g_1$ (respectively,  $g_2$) acts trivially on $Y^\perp$ (respectively,  $Z^\perp$), and $\det g_1=\det g_2=-1$. \itf  $C_{S O(V)}(s)$ contains a connected subgroup of index 2, and hence is not connected.   

For the additional claim, observe that ${\mathbf L}$ coincides with the stabilizer in ${\mathbf G}$ of the direct sum of 
subspace $(U_1+U_1')\oplus \cdots \oplus (U_t\oplus U_t')+U' $ of $\overline{V}$, where $U'$ is non-degenerate or $\{0\}$, and $U_i,U_i'$
are totally singular and dual to each other ($i=1\ld t$). Therefore, ${\mathbf L}$ is isomorphic to the direct product of $SO(U')$ with
$GL(U_1)\times \cdots \times GL(U_t)$. As $s\in {\mathbf L}$, all $U',U_i,U_i'$ are stabilized by $s$. Let $s_i$ be the restriction of 
$s$ to $U_i$, and $s'$  the restriction of 
$s$ to $U'$. Then $C_{{\mathbf L}}(s)$ is isomorphic to the direct product of $C_{SO(U')}(s')$ with $C_{GL(U_1)}(s_1)\times \cdots \times C_{GL(U_t)}(s_t)$. Each group  $C_{GL(U_i)}(s_i)$ is well known to be connected. As $C_{{\mathbf G}}(s)$ is connected, either $1$ or $-1$ is not an \ei of $s$ on $\overline{V}$, and this remains true for $s'$. So  $C_{SO(U')}(s')$ is connected, whence the result.

\bl{os5}  Let $G=SO(V)$, $q$ odd, and let $t\in G$ be a semisimple element. Let $V^t$ be the $1$-eigenspace of t on V, and W the unique t-stable complement of $V^t$. Let $t_1$ be the restriction of t to $W$. Then $SO(V^t)\times C_{SO(W)}(t_1)$ is a subgroup of $ C_G(t)$ of index at most $2$, and the index equals $2$ \ii $V^t\neq 0$ and $-1$ is an \ei of $t$. \el
Let $Y$ be the $-1$-eigenspace of $t$ and $M$ the unique $t$-stable complement of $V^t+Y$ in $V$. Let $t_1$ be the projection of $t$ to $M$. Obviously,
$C_G(s)\subseteq O(V^t)\times O(Y)\times C_{O(M)}(t_1)$. By Lemma \ref{y9}, $C_{O(M)}(t_1)=C_{SO(M)}(t_1)$.
This implies the lemma if $V^t$ or $Y$ is the zero space. Otherwise, $C_G(t)=\lan O(V^t)\times O(Y)\times C_{O(M)}(t_1), g\ran$, where $g\in O(V^t)\times O(Y)\times O(M)$ is such that the projection of $g$ to $M$ is the identity, whereas the projections to $V^t$ and to $Y$ are not in $SO(V^t)$, $SO(Y)$, respectively. 
So the lemma follows.

\med

Let $\al\in\{1,-1\}$. In many cases it is convenient to use uniform notation for $SO_{2n}^+(q)$ and $SO_{2n}^-(q)$ by writing $SO^\al_{2n}(q)$   and interpreting $\al$ as the plus sign if $\al=1$ and the minus sign if $\al=-1$. 
 The \f is well known:

\bl{h85}  Let $G=SO^\al_{2n}(q)=SO(V)$ if q is odd, otherwise let $G=\Omega^\al _{2n}(q)$. Then $G$ contains a cyclic subgroup $T$ of order $q^n-\al$. Moreover, 
if $\al=-1$ then $T$ is irreducible, otherwise T stabilizes a maximal  totally singular subspace of V and acts on it irreducibly.  In addition, $T$ is a maximal torus of $G$.
\el
Suppose first that  $\al=1$. Then $V=U+U'$, where $U,U'$ are totally singular subspaces of $V$ of dimension $n$. Moreover, there are bases in $U,U'$ such that the matrix $\diag(g,{}^tg\up)$ belongs to $O(V)$ for every $g\in GL(U)$ \cite{Die}. (Here ${}^tg$ means the transpose of $g$.)  Note that $GL(U)\cong GL_n(q)$ contains an \ir element $h$ of order $q^n-1$. We set $T= \lan \diag(h,{}^th\up)\ran $.  
Furthermore, $\lan h\ran$ is self-centralizing in $GL(U)$. 

Let $\al=-1$. By Huppert \cite[Satz 3]{Hu}, $O(V)$ contains a cyclic \ir subgroup $T$ 
of order $q^n+1$.  In fact, $t\in G$. This is trivial if   $q$ is even as $|O(V):\Omega(V)|=2$ and $q^n+1$ is odd. Let $q$ be odd.  View $t$ as an element of the group $O(\overline{V})$. 
Then there is a basis $B$, say, of  $\overline{V}$, under which
the matrix of  $t$ is diagonal. Moreover,  $t$ does not have \eis $\pm 1$. (Indeed, 
$t$ acts
irreducibly on $V$, and hence no \ei of $t$ on $V$ belongs to $ F_q$.) 
This implies that every vector of $B$ is singular,
and hence $B$ can be assumed to be   a hyperbolic basis. (That is,
one can replace every element of $B$ by a suitable multiple so that the Gram matrix of the associated bilinear form under the new basis becomes block-diagonal
with blocks $\begin{pmatrix}0&1\cr 1&0\end{pmatrix}$.) It easily follows from this that  $\det t=1$.

  We show  that $T$ is self-centralizing. Set $Y=G$ if $q$ is odd, otherwise denote by $Y$ the symplectic group of the associated bilinear form on $V$.  It is well known that  there exists an 
 involutory anti-\au $\si$  of 
${\rm End}\, V$ such that       $Y=\{x\in {\rm End}\, V: x\si(x)=\Id\}$.  
Let $K$ be the centalizer of $T$ in ${\rm End}\, V$.  By Schur's lemma, $K$ is a field of order $q^{2n}$.
 Set  $T_1=
\{x\in K: x\si(x)=\Id\}$. Then $T_1\subset Y$ and  $C_{Y}(T_1)=T_1$. As $\si$ is an \au of $K$ of order 2, an easy Galois argument implies that $|T_1|=q^n+1$. As $T\subseteq T_1$, we have $T=T_1$. As $C_Y(T)=T$,  it follows that  $C_G(T)=T$. 

Finally,  let $D\cong GL(U)$ be the group $\{\diag(g,{}^tg\up): g\in GL(U)\}$ and${\mathbf D}\cong GL(\overline{U})$ . As  $T$ is a cyclic group, $T$ is contained in a maximal torus $T'$ of $G$ (respectively, $D$) if $\al=-1$ (respectively,  $\al=1$), see \cite[Ch. II,  1.1.1]{SS}. Let ${\mathbf T}'$ be a maximal torus $\GC$ (respectively, ${\mathbf D}$) such that $T'={\mathbf T}^{prime Fr}$.   As $C_{G}(T)=T$ (respectively, $C_D(T)=T$ unless $n=2,q\leq 3$), we have $T=T'$. If  $n=2,q\leq 3$ then it is easy to observe that $T'=T$.  
As the rank of the algebraic group ${\mathbf D} $ equals the rank of ${\mathbf G}=SO(\overline{F}_q)$, a maximal torus of $D$ remains maximal in $ G$.

\bl{fv1}  Let $\dim V$ be odd and let $T\subset  G  \subset O (V)$ be a maximal torus. Then 
$T$ stabilizes a non-degenerate subspace $V'$, say, of V, and hence $T\subset O(V')$, when $O(V')$
is viewed as a subgroup of $O(V)$. In addition, $T\subset SO(V')$ if q is odd, otherwise $T\subset \Omega(V')$.   
\el
Observe first that there is a vector $v\in V$ fixed by $T$. Indeed, it is
well known that such a vector exists in $V\otimes \overline{F}_q$.
Let $b_1\ld b_{n}$ be a basis in $V$. Consider linear equations
$t\sum_ix_ib_i=\sum_ix_ib_i$ for every $t\in T$ with respect to
indeterminates  $x_1\ld x_{n}$. By the above, these equations  have a common 
solution with $x_i\in \overline{F}_q$. Therefore, there exists a
solution with $x_i\in F_q$, as claimed. Set $V^T=\{y\in V:ty=y$ for all $t\in
T\}$. If $q$ is odd then $V^T$ is non-degenerate, and hence there is an
 anisotropic  vector $v\in V^T$, so we take $V'=v^\perp$. Suppose that $q$ is even.
By Maschke's theorem, $T$ is completely reducible on $V$, and hence there a $T$-stable
complement $V'$ of $\lan v\ran$ in $V$. One observes that $V'$ is non-degenerate,
so $T\subset O(V')$. The additional statement is well known.

 \begin{lemma}\label{no6}    Let $V$ be an orthogonal space, and let  $U_1,U_2$ be totally singular subspaces of $V$ of equal dimension. 
Then $gU_1=U_2$ for some $g\in G$, unless when $\dim U_1=\dim U_2=\dim V/2$ (and hence $V$ is of Witt defect $0$).
In the exceptional case there are  two $G$-orbits of the subspaces in question. \el
 By Witt's theorem, $hU_1=U_2$ for some $g\in O(V)$. Suppose that  $\dim U_1<\dim V/2$. 
 Let $W$ be a complement of $U_1$ in $U_1^\perp$. Then $W\neq 0$ is non-degenerate and $U_1\subset W^\perp$. Furthermore,
there is an element $x\in O(V)$ such that $xv=v$ for all $v\in W^\perp$ and 
(a) $\det x=-1$ if $q$ is odd and (b) the projection of $x$ into $O(W)$ is not in $\Omega(W)$ if $q$ is even.    Then, if $h\notin G$, we have  $hxU_1=U_2$ and
$hx\in G$ in (a). So the first claim follows. Suppose that $\dim U_1=\dim V/2$. Then 
the result follows  by  \cite[Lemma 2.5.8(ii)]{KL}.

\begin{lemma}\label{ag1}   Let ${\mathbf G}$ be a reductive algebraic group, ${\mathbf Z}$ a finite central subgroup of ${\mathbf G}$, ${\mathbf G}_1:={\mathbf G}/{\mathbf Z}$ and let $\eta:{\mathbf G}\ra {\mathbf G}_1$ be the natural homomorphism. Let $Fr$ be a Frobenius endomorphism of ${\mathbf G}$, $G:={\mathbf G}^{Fr}$ and $G_1:={\mathbf G}_1^{Fr}$. Let ${\mathbf T}$ be a maximal $Fr$-stable torus of ${\mathbf G}$, ${\mathbf T}_1=\eta({\mathbf T})$, $T={\mathbf T}^{Fr}$,  $T_1={\mathbf T}_1^{Fr}$ and $Z={\mathbf Z}^{Fr}$. Then $|G_1:\eta(G)|=|Z|$ and $G_1=\eta(G)\cdot T_1$. 
\end{lemma}
Proof. Obviously, $\eta(G)\subset G_1$,  $\eta(T)\subset T_1$. Furthermore,
$\eta(G)$ is a normal subgroup of $G_1$,  $|G|=|G_1|$ and $|T|=|T_1|$, see  \cite[4.2.3]{Ge}. Obviously, $Z=G\cap {\mathbf Z}=G\cap \ker\eta$. So $|G|=|G_1|$ implies $|G_1|=|\eta(G)|\cdot |Z|$.  

Observe that $\eta (T)=\eta(G)\cap T_1$. (Indeed, let $x\in \eta(G)\cap T_1$. As ${\mathbf T}_1=\eta({\mathbf T})$, there is $t\in {\mathbf T}$ such that $\eta(t)=x$. Let $g\in G $ with $\eta(g)=x$. Then $\eta(g)=\eta(t)$, or $\eta(t\up g)=1$. Therefore, $z:=t\up g\in {\mathbf Z}$. As ${\mathbf T}$
is a maximal torus of ${\mathbf G}$, it follows that ${\mathbf Z}\subset
{\mathbf T}$, and hence $z\in {\mathbf T}$. Let $t'=tz$. Then $g=tz=t'\in {\mathbf T}$. This means that $g\in {\mathbf T}$, and hence $g\in {\mathbf T}\cap G=T$. So $x\in \eta(T)$, as claimed.) 

Set $m=|T_1/(\eta(G)\cap T_1)|=|T_1/\eta(T)|$. As ${\mathbf Z}\subset {\mathbf T}$, we have $Z\subseteq T$, and hence $Z=T\cap \ker \eta$. As above, $|T|=|T_1|$ implies $|Z|=|T_1/\eta(T)|=m$. Therefore,
$\eta(G)\cdot T_1$ is of order  $|\eta(G)|\cdot m=|G_1|$, whence the result. 

\medskip
In the classical group theory an important role is played by
the notion of spinor norm \cite{Die, Gr}.   This is a \ho $\nu:O(V)\ra K$, where $K$ is an abelian group of exponent 2 and of order at most 4. It is well known that
 $\Omega(V)$ is a
subgroup of $SO(V)$ of index 2 unless   $q$ is even and $\dim V$ is odd, or $\dim V=1$,  see \cite{Die} or \cite{KL}. 
Furthermore, if $W$ is a non-degenerate subspace of $V$ then 
$\Omega(W)=O(W)\cap \Omega(V)$, when $O(W)$ is viewed as a subgroup of $O(V)$.

If $q$ is even and $\dim V$ is odd then $O(V)\cong Sp(U)$, where $U$ is a symplectic space of dimension $\dim V-1$,
see  \cite[Theorem 14.2]{Gr}; therefore in this case $O(V)$ is simple unless $(n,q)= (3,2) $ or $(5,2)$.

Note that $G=SO(V)$ has a unique subgroup of index 2,
unless $\dim V=1$, or $q$ is even, $\dim V$ is odd,  or $G=SO^+_4(2)$, see \cite[2.5.7]{KL}. With this exceptions,
 $\Omega(V)$ can be defined as the subgroup of index 2 in $SO(V)$.

\begin{lemma}\label{h01}  Assume that  $q$ is odd and $\dim V>1$. Let $T$ be a maximal torus in $SO(V)$.
Then $T$ contains an element of spinor norm $-1$.\end{lemma}
Proof. If $n=2$, the statement is trivial as  $SO(V)=T$.
Otherwise,  this follows from Lemma \ref{ag1}.  (Indeed, in notation of Lemma \ref{ag1} we have $|Z|=2$ and hence  $|G_1:\eta(G)|=2$. Therefore, $\eta(G)=\Omega(V)$  as $\Omega(V)$ is the only subgroup of inder 2 in $SO(V)$).

\begin{lemma}\label{hu0}  Let $\dim V=2n$, $q$ odd. Suppose that the Witt defect of V equals $1$.
Then $SO(V)\setminus \Omega(V)$
contains an element $t$ of order $q^n+1$ such that $C_G(t)=C_G(t^2)=\lan t\ran$.
\end{lemma}
Proof. If $n=1$ then $SO(V)$ is abelian, and the statement is trivial. Let $n>1$.
By Lemma \ref{h85}, $SO(V)$ contains a cyclic \ir subgroup $T$ 
of order $q^n+1$, which is a maximal torus in $SO(V)$.  Let $T=\lan t \ran$.

By Lemma \ref{h01}, $T$ is not contained in $ \Omega(V)$, and hence $t\notin  \Omega(V)$. 
We has shown in Lemma \ref{h85} that $T$ is self-centralizing, so $C_G(t)=T$. We show that $T=C_G(t^2)$.   
Let $K$ be as in the proof of  Lemma \ref{h85}. If $t$ is \ir in $V$ then the claim follows by Schur's lemma. 
Suppose that $t$ is reducible. Then $t$ belongs to a proper subfield $K_1$, say,  of $K$. Let   $K^\times, K_1^\times$ be the multiplicative groups of $K,K_1$, respectively. Then $t^2\in K_1$ implies $|K^\times/ K_1^\times|=2$, 
which is false. So the lemma follows.

\medskip
If $\GC$ is an algebraic group then $\GC^0$ denotes the connected component of $\GC$. 
 
\begin{lemma}\label{h00}  Assume that $q$ is odd. 
Then there exists $t\in  G\setminus \Omega(V)$ 
such that $|C_G(t)|_p=|C_G(t^2)|_p$ and $V^t=V^{t^2}$. In addition, $C_{\GC}(t)^0=C_{\GC}(t^2)^0$.
\end{lemma}
Proof. Set $G=SO(V)$ and $\dim V=n$. If $n=2$ then $G$ is abelian, so the statement is trivial. Let $n>2$.
In view of Lemma \ref{hu0}, we can assume that the Witt defect of $V$ equals 0. 
Let $W$ be a non-degenerate subspace of $V$
of Witt defect 1  and  $2k=\dim W$. We choose
$W$ so that  $2k=n-2$ if $n$ is even, otherwise $2k=n-1$. By
Lemma \ref{hu0},  there is an element $t\in SO(W)\setminus \Omega(W)$  of order $q^k+1$ acting on $W$ irreducibly.
Let us view $t$ as an element of $G$. Then $t\notin \Omega(V)$. 

If $G=SO^+_4(3)$ or $  SO_3(3)$ then  $|t|=4$, and the $-1$-eigenspace of $t^2$ on $V$ is of dimension 2.
So $C_G(t^2)$ is a subgroup of  $O(W)\times O(W^\perp)$, which has no element of order $3$. The equality $V^t=V^{t^2}$ is here obvious. In addition, the connected component of $C_{\GC}(t^2)$ is $SO(\overline{W})\times SO(\overline{W}^\perp)$, so the additional claim of the lemma  is true in this case.

With exception of the above two groups,  $t^2$ is \ir on $W$. Indeed, let $K$ be the enveloping algebra of $\lan t\ran$ in ${\rm End}\,W\cong {\rm Mat}(2k,q)$. By Schur's lemma,
$K$ is a field of order $q^{2k}$. If $t^2$ is reducible then $t^2$ belongs to a proper subfield $L$ of
$K$, in fact $K/L$ is a quadratic extension. Let $\gamma$ be the Galois \au of $K/L$.
Then the group $\{x\in K: x\gamma(x)=1\}$ is of order $q^k+1$ and  $t$ is a generator of this group.
Therefore, $t^2\in L$ \ii $q^k+1=2$ or 4. This implies $k=1,q=3$, and hence $G\in\{SO^+_4(3),  SO_3(3)\}$. 

Thus, $t^2$ is \ir on $W$. \itf $C_G(t)=C_G(t^2)$, as claimed. As all \eis of $t^2$ on $\overline{W}$ are distinct, it easily follows that  $C_{\GC}(t)^0=C_{\GC}(t^2)^0$.

\med Remark. One observes that if $G\notin\{ SO^+_4(3),  SO_3(3)\}$ then $|C_G(t)|_p=|C_G(t^2)|_p=1$. 

\begin{lemma}\label{h03}  Let $V$ be a non-degenerate subspace of an orthogonal space $U$,  let $G=SO(V)\subset D=SO(U)$ be a natural embedding, and let $t\in G$  be as in Lemma $\ref{h00}$.
Then $|C_D(t)|_p=|C_D(t^2)|_p$. In addition, if ${\mathbf D}:=SO(\overline{V})$ then $C_{{\mathbf D}}(t)^0=
C_{{\mathbf D}}(t^2)^0$.
\el
Let $ W=V^\perp$ so $U=V\oplus W$. Then  $W$ is a trivial $F_qG$-module. Therefore, $U^t=V^t+W=V^{t^2}+W=U^{t^2}$. Set $V'=(V^t)^\perp,U'=(U^t)^\perp$. Then $V'=U'$.  Let $t'$ be the projection of $t$ to $V'$. As $V^t=V^{t^2}$, it follows that  $t$ does not have \ei $-1$. 
By  Lemma \ref{os5}, $C_G(t)=SO(V^t)\times C_{SO(V')}(t')$. Furthermore, $|C_{G}(t)|_p=|C_{G}(t^2)|_p$,
so $|C_{O(V')}(t)|_p=|C_{O(V')}(t^2)|_p$, and hence $|C_{O(U')}(t)|_p=|C_{O(U')}(t^2)|_p$ as $U'=V'$.
We have $C_{  D}(t)=SO(U^t) \times C_{O(U')}(t')$ and $C_{D }(t^2)=SO(U^t) \times C_{O(U')}(t^{\prime 2})$. 
So the first statement of the lemma follows. Similarly, $C_{{\mathbf D}}(t)=SO(\overline{V}^t)\times C_{SO(\overline{V}')}(t')$, so $C_{{\mathbf D}}(t)^0=SO(\overline{V}^t)\times C_{SO(\overline{V}')}(t')^0$. By Lemma \ref{h00},
$C_{SO(\overline{V}')}(t')^0=C_{SO(\overline{V}')}(t^{\prime 2})^0$. As $C_{{\mathbf D}}(t^2)^0=SO(\overline{V}^{t^2})\times C_{SO(\overline{V}')}(t^{\prime 2})^0$, the second statement follows. 

\medskip

The \f result by Gow and Szechman \cite[Theorem 4.1]{GS} considerably simplifies the computational
aspect of our reasoning below. This is called the character comparison theorem in \cite{GS} and generalizes
an earlier result by Kn\"orr \cite[Proposition 1.1]{Kn}.

\begin{theo}\label{ct1}  {\rm (The comparison theorem)}
 Let X be a finite group, and $p$ a prime divisor of $|X|$. Let $\phi,\psi$ be generalized characters of X.
Suppose that $\phi(g)=\pm \psi(g)=\pm p^{m(g)}$ for every
$p'$-element $g$ of G, where $m(g)\geq 0$  is an integer. Suppose
also that $\phi(1)=\psi(1)$. Then there exists a linear character
$\lam$ of $X$ such that $\lam^2=1$ and $\phi(g)= \lam(g) \psi(g)$
for all $p'$-elements $g\in X$.
\end{theo}

\section{Orthogonal decompositions and maximal tori}
Let $V$ be the natural module for $G$.  We need to compute the restriction
$\om_{G}|_T$ for every maximal torus $T$ of $G$. In order to express the result in a
convenient and uniform way, we introduce a so called $T$-decomposition of $V$, 
see \cite{HZ09}. If $q$ is even, we assume $\dim V$ to be even. 
 Let

\begin{equation}\label{eq2} V=V_0\oplus V_1\oplus \cdots \oplus
V_{k}\oplus V_{k+1}\oplus \cdots \oplus V_{k+l},\end{equation} 
where

$(a)$ $k,l\geq 0$ and $\dim V_i$ is even for $i>0;$

$(b)$ $V_1\ld V_{k+l}$ are non-degenerate subspaces of $V$
orthogonal to each other, $V_0= (V_1\oplus \cdots \oplus
V_{k+l})^\perp$,  $V_0=0$ if $n$ is even and $\dim V_0= 1$ otherwise;  
 
 $(c)$  $V_1\ld V_k$ are of  Witt defect $0$;

 $(d)$ $V_{k+1}\ld V_{k+l}$ are of Witt defect $1;$

$(e)$ $\dim V_1\leq \cdots \leq \dim V_k$ and $\dim V_{k+1}\leq \cdots \leq \dim V_{k+l}$. 

\med
We call this  an {\it  orthogonal decomposition of} $V$.
Note that   the part  $V_1\oplus \cdots \oplus V_{k}$ or
$V_{k+1}\oplus \cdots \oplus V_{k+l}$ may be absent. We express
this by writing $k=0$ or $l=0$. Additionally, we set $$V'=
V_1\oplus \cdots \oplus V_{k},~~{\rm and}~~ V''= V_{k+1}\oplus
\cdots \oplus V_{k+l}.$$
Obviously, the Witt defect of $V'$ is 0.

The list of the dimensions $\dim V_1\ld \dim V_k$ and $\dim V_{k+1}\ld \dim V_{k+l}$ is an essential invariant of 
an orthogonal decomposition. There is a convenient way
(for our purposes)   to encode the list of these dimensions in terms of two functions $i\ra d_i$, $j\ra e_j$, where
$d_i,e_j\geq 0$ are integers. We often record these functions in the form $[1^{d_1}, 2^{d_2},\ldots ]$ and 
 $[1^{e_1},2^{e_2}, \ldots ]$. These are interpreted as follows. The entry $i^{d_i}$ tells us that in the list $V_1\ld V_k$ in (\ref{eq2})
there are $d_i$ subspaces of dimension $2i$.   Similarly, the entry $j^{e_j}$ tells us that in the list $V_{k+1}\ld V_{k+l}$ 
there are $e_j$ terms of dimension $2j$. Observe that some $d_i,e_j$ may be zeros. Obviously, $k=d_1+d_2+\cdots $ and $l=e_1+e_2+\cdots $.
Note that \begin{equation}\label{eq2a}\dim V=\sum_i 2id_i+\sum_i 2je_j~{\rm  if}~ n~ {\rm is ~even,~ and} ~\dim V=1+\sum_i 2id_i+\sum_i 2je_j~ {\rm if} ~n~{\rm  is~ odd}.\end{equation} 
Computing the Witt defect of $V$ in terms of $V_i$ \cite[2.5.11]{KL}, one observes that $l$ is even (respectively, odd)
if $G=SO_{2n}^+(q)$ (respectively, $G=SO_{2n}^-(q)$).

\bl{ca6}   For every orthogonal decomposition $(\ref{eq2})$ there is  a maximal torus T of G of order $\Pi_{i,j}(q^i-1)^{d_i}\cdot \Pi_j (q^j+1)^{e_j}$. In fact,
 
\begin{equation}\label{eqT} T=T_1\times \cdots\times T_k\times T_{k+1}\times\cdots \times T_{k+l},\end{equation}
where $ TV_i=V_i$ and $T_i$ is the restriction of T to $V_i$ for  $i=1\ld k+l$. 
\el
To every  decomposition (\ref{eq2}) one can correspond an  abelian
subgroup $T$ of $SO(V)$, which is the direct product $
T_1\times \cdots \times   T_k \times T_{k+1}\times \cdots\times
T_{k+l}$  as follows. Let  $T_i\subset SO(V_i)$ be a
maximal torus of $SO(V_i)$ of order  $q^{\dim V_i/2}-1$ for $1\leq i\leq k$  and of order  $q^{\dim V_i/2}+1$ for $k+1
\leq i\leq k+ l$ according with Lemma \ref{h85}. Let $G={\GC}^{Fr}$. Set ${\mathbf Y}=SO(\overline{V}_1)\times \cdots \times  SO(\overline{V}_{k+l})$ and $Y=SO(V_1)\times \cdots \times  SO(V_{k+l})$. Then ${\mathbf Y}$ and $Y$ are subgroups of ${\mathbf G}=SO(\overline{V})$ and of $G$, respectively (in the obvious sense).  Then $Y={\mathbf Y}\cap G$. It follows that $Fr({\mathbf Y})={\mathbf Y}$ and  $Y={\mathbf Y}^{Fr}$. Saying that $T_i$ is a maximal torus of $G_i$ means that $T_i={\mathbf T}^{Fr}_i$  for some  $Fr$-stable maximal torus ${\mathbf T}_i$ of $SO(\overline{V}_i)$. Therefore, ${\mathbf T}:={\mathbf T}_1 \times \cdots \times {\mathbf T}_{k+l}$ is a maxinal $Fr$-stable torus of ${\mathbf Y}$ and $T={\mathbf T}^{Fr}$. By dimension reason,  ${\mathbf T}$ is a  maximal $Fr$-stable torus of ${\mathbf G}$.
(One can also mimic the reasoning in \cite[\S 2.4]{HZ09}.)

 \med
Observe that every maximal torus $T$ of $G$ determines an orthogonal decomposition with properties described in Lemma \ref{ca6}.
This means (in more precise terms) that if  $T=\TC^{Fr}$ for an $Fr$-stable  maximal torus $\TC$ of the algebraic group $\GC$, then 
 $\TC$  determines an orthogonal decomposition in question, which we call the $T$-{\it decomposition of} $V$ in this paper.  For symplectic groups, which appear below as the dual groups of $SO_{2n+1}(q)$, this notion was introduced in \cite{HZ09}.

\begin{lemma}\label{my6}  If $\dim V$ is odd
 then there is a bijection between
the $G$-conjugacy classes of maximal $Fr$-stable tori in $\GC$ and
the 
pairs of functions $[1^{d_1}, 2^{d_2},\ldots ]$ and 
 $[1^{e_1},2^{e_2}, \ldots ]$ satisfying $(\ref{eq2a})$. 
\end{lemma}

Proof. This is well known  (the case with $q$ even is not excluded), see for instance \cite{HZ09}. (Note that
in this case there is a  bijection between the $G$-conjugacy classes
of maximal $Fr$-stable tori in $\GC$ and the conjugacy classes of
the Weyl group $W$ of $\CG$.) 

\med
Thus, if $\dim V$ is odd then   the maximal tori of $G$ are parametrized by the above functions. Note that this parametrization differs a bit from
the parametrization by bipartitions which are pairs $(\lam_1,\ldots \mu_1,\ldots )$, where $2\lam_i=\dim V_i$ for $i=1\ld k$,
and $2\mu_j=\dim V_{k+j}$ for $j=1\ld l$.

The situation is more complex if $n$ is even. However, if one wishes
to consider the maximal tori of $G$ up to conjugacy in the full
orthogonal groups, then  we have a similar statement:

\begin{lemma}\label{my7}  Suppose that $\dim V=2n$ is even.
 Then there is a bijection between
$O(V)$-conjugacy classes of maximal $Fr$-stable tori in $\GC$ and
the pairs of functions $[1^{d_1}, 2^{d_2},\ldots ]$ and 
 $[1^{e_1},2^{e_2}, \ldots ]$ satisfying $(\ref{eq2a})$. 

\end{lemma}

Proof. Note that $O(V)$ is a subgroup of $O(\overline{V})$ and $SO(V)=SO(\overline{V})\cap O(V) $.
Consider $G$ as a subgroup of $H=SO(W)$, where $\dim W=2n+1$. Then a maximal torus $T$ of $G$ remains maximal in $H$. Set $X=\{\diag(\det g\up , g): g\in O(V)\}$. Obviously, $X\cong O(V)$ and the restriction of $X$ to $V$ is exactly $O(V)$.
Then two maximal tori $T,T'$ of $G$ are conjugate in $H$ \ii they are conjugate in $X$. Indeed, tori ${\mathbf T}$, ${\mathbf T}'$ of ${\mathbf G}$ stabilizes a complement $\overline{U}$ of $\overline{V}$ in $\overline{W}$,  and $\overline{U}$ is the fixed point space of each  ${\mathbf T}$, ${\mathbf T}'$ on  $\overline{W}$.  Therefore, if  $h{\mathbf T}h\up ={\mathbf T}'$ for $h\in H$ then $h\overline{U}=\overline{U}$. Simirlarly, $h\overline{V}=\overline{V}$. Therefore, 
$hV=V$ and the restriction of $h$ to $V$ belongs to $O(V)$. As $h\in SO(W)$, it follows that $h\in X.$ Now 
the lemma follows from Lemma \ref{my6}. 

\medskip
We need some observations on the Weyl groups of the BN-pairs of type $B_n$ and $D_n$. 
We denote them here by $W(B_n)$ and  $W(D_n)$, respectively. 
Recall that  $W(B_n)$ is a semidirect product of a normal subgroup $A$ and of the symmetric group $S_n$; in addition,  $A$ is exponent 2 and rank $n$. 
In the natural realization as a group of $(n\times n)$-matrices   over the rationals, $W(B_n)$  is exactly the group of monomial matrices with non-zero entries $\pm 1$. So there is a basis $b_1\ld b_n$ of the underlying space such that the elements of $W(B_n)$   permute  the set $\pm b_1\ld \pm b_n$. If $w\in W(B_n)$, we denote by $\tilde w$
the projection of $w$ into $S_n$. Let  $A_0$ be the  subgroup of $A$ formed by matrices of determinant 1. Then   $|A:A_0|=2$, and   $W(D_n)$ 
is a semidirect product of  $A_0$ and  $S_n$.

 The conjugacy classes of $W(B_n)$ and  $W(D_n)$ are determined in \cite{Ca}.
Let $w\in W(B_n)$. The conjugacy classes of $B=W(B_n)$ are determined by  pairs of functions $i\ra d_i$, $j\ra e_j$,
and a representative $w$ of the respective conjugacy class is the permutation of the set $\{\pm b_1\ld \pm b_n\}$ obtained as follows. Let $k=\sum _{i} id_i$ and $l=\sum _{j} je_j$ (so $k+l=n$), and let $\pi\in S_k\times S_l\subset S_n$ be the permutation in  which the cycle of length $i$ occurs with \mult $d_i$, and the cycle of length $j$ occurs with \mult $e_j$. If $m\leq k$ then  
$w^t(b_m)\in \{b_{1}\ld b_k\}$ for any $t$, that is, $w$ acts by permuting the basis elements $b_1\ld b_k$.  If $m>k$ and $\tilde  w^r\cdot m$ is a cycle of length $t$ then $w$ can be chosen so that  $w^rb_m\in \{b_1\ld b_n\}$ for $r<t$ and $w^tb_m=-b_m$. In addition, $w\in W(D_n)$ \ii
$l$ is even.

\bl{ee1}  Set $D=W(D_n)$ and $B=W(B_n)$.  For $w\in B$ denote by $w^D$ (respectively, $w^B$) the $D$- (respectively, $B$-) conjugacy class of $w$. Let $i\ra d_i$, $j\ra e_j$ be the functions
which determine  $w^B$. 

$(1)$ $C_{D}(w)=C_{ B}(w)$ \ii $w\in D$,
$l=0$  
and all  $i$ with $d_i>0$  are even.  In the exceptional case $|C_{B}(w):C_{D}(w)|=2.$

$(2)$  Either
$w^D=w^B$, or $w\in D$, $l=0$, each cycle of  $\tilde w$ has even size 
 and  $w^B$ is a union of two $D$-conjugacy classes of equal size.  

\el Clearly, (1) and (2) are equivalent. For (2), if $w\in D$, the result is contained in Carter \cite[Proposition 25]{Ca}. 
Let $w\notin  W(D_n)$. As $|B:D|=2$, it follows that $C_{D}(w)\neq C_{B)}(w)$ \ii $ C_{ B}(w)$ contains an element 
$x\notin C_{D}(w)$. As $w\notin D$, we can take $x=w$.

\med
Our next aim is to obtain a
 formula for $W(T)$. Recall that $W(T)$  is defined to be 
$(N_{{\mathbf G}}({\mathbf T})/{\mathbf T})^{Fr}=N_{G}({\mathbf T})/T$.  
It is well known that there is a bijection between the
conjugacy classes of $Fr$-stable maximal  tori of ${\mathbf G}$ and the so called $Fr$-conjugacy
classes in $W({\mathbf G})$, the Weyl group of ${\mathbf G}$ \cite[3.23]{DM}. 
If $G=B_n(q),D_n^+(q)$ then $W({\mathbf G})=W(B_n)$ and  $W(D_n)$, respectively, and the $Fr$-conjugacy classes of $W({\mathbf G})$  are ordinary conjugacy classes.

\med The \f definition plays an essential role in what follows. 

\begin{defi}\label{de4}  Let T be a maximal torus in G, and let $i\ra d_i$, $j\ra e_j$ be the corresponding 
functions. We say that T is neutral if $l=0$, and exceptional if $l=0$ and all $i$ with $d_i>0$ are even.  
\end{defi}

\begin{lemma}\label{a33} $(1)$ Let $H=SO_{2n+1}(q)$ and  $T$ a maximal torus of H corresponding to the functions 
$i\ra d_i$, $j\ra e_j$. Then $|W(T)|=\Pi_{i} (2i)^{d_i}d_i!\cdot 
\Pi_{j} (2j)^{e_j}e_j!$.
\medskip

$(2)$ Let $G=SO^+_{2n}(q)$   and
 $T$ be a maximal torus of G corresponding to a function 
$i\ra d_i$, $j\ra e_j$. 
Then $|W(T)|=\Pi_{i} (2i)^{d_i}d_i!$
 if T is exceptional, 
otherwise $|W(T)|=\Pi_{i} (2i)^{d_i}d_i!\cdot 
\Pi_{j} (2j)^{e_j}e_j!/2$. 

\medskip

$(3)$ Let $G=SO_{2n}^-(q)$   and
 $T$ be a maximal torus of G corresponding to a function 
$i\ra d_i$, $j\ra e_j$. 
Then $|W(T)|=\Pi_{i} (2i)^{d_i}d_i!\cdot 
\Pi_{j} (2j)^{e_j}e_j!/2$.
 
\medskip
$(4)$ Let $G=SO(V)$, where $\dim V=2n$, and   ${\mathbf G}=SO_{2n}(\overline{F}_q)$. Let $Fr$ be a Frobenius endomorphism of  ${\mathbf G}$ such that  $G={\mathbf G}^{Fr}$.  Let
${\mathbf T}$ be a maximal $Fr$-stable torus of ${\mathbf G}$ and $T={\mathbf T}^{Fr}$. 
Then  $N_{G}({\mathbf T})\neq  
 N_{O(V)}({\mathbf T}) $, 
 except for the case where $G=SO^+(2n,q)$ and $T$ is exceptional. 
\end{lemma}
Proof. (1) See \cite[Lemma 2.3]{HZ09}. 

(2) Let $T=T_w$, where $w\in W(\GC)\cong W(D_n)$. Set $B=W(B_n)$, $D=W(D_n)$ as in Lemma \ref{ee1}. Recall that $W(T)\cong C_{W(\GC)}(w)=C_{D}(w)$ \cite[3.3.6]{C}. We view $G$ as a natural  subgroup of $H$,  and also consider $T$ as a maximal torus of $H$. Then $T=T_w$, where $w$ is viewed as an element of $W({\mathbf H})\cong W(B_n)=B$.  Then $W_{H}(T)\cong C_{{\mathbf H}}(w)$, where we use the subscript $H$ in $W_{H}(T) $ to indicate that $T$ is viewed as a torus in $H$. 
As above, for $w\in D$ denote by $w^D$ (respectively, $w^B$) the conjugacy class of $w$ in $D$ (respectively,  $B$). 
Recall that  $|B:D|=2$. So  $w^B\subset D$, and we have $|C_{D}(w)|=\frac{|D|}{|w^D|}$ as well as   $|C_{B}(w)|=\frac{|B|}{|w^B|}$. By Lemma \ref{ee1}, 
either $w^D=w^B$ or $w^B$ is a union of two  $D$-conjugacy classes of equal size.
In the former case $|C_{D}(w)|=\frac{|D|}{|w^D|}=\frac{|B|}{2\cdot |w^B|}=\frac{1}{2}|C_{B}(w)|$, in the latter case,  we have $|C_{D}(w)|=|C_{B}(w)|. $
So (2) follows from (1) and the fact  that $w^D\neq w^B$ \ii all $i$ with $d_i>0$ are even. 

(3) Now $G=SO_{2n}^-(q)$. 
The $G$-conjugacy classes of maximal tori are in bijection with the $Fr$-conjugacy classes in $W({\mathbf G})\cong W(D_n)=D$. These are defined in terms of the action of $Fr$ on $W(\GC)=D$. This is known to be realized as the conjugation by some element 
 $a\in B\setminus D$, see \cite[\S 11]{St}  where the action in question is explicitly described. Furthermore, an  $Fr$-conjugacy class of $D$ is of shape $\{x\up waxa\up :x\in D\}$, which  is the orbit of $w$ under the action $w\ra x\up w axa\up $, $x\in D$. By \cite[3.3.6]{C}, $W(T_w)$ is isomorphic to the point stabilizer of $w$ in $D$ under this action.  The mapping $w\ra wa$ transforms the action in question to the $D$-conjugacy action on $B\setminus D$. Therefore, $W(T_w)\cong C_D(wa)$. 
  By Lemma \ref{ee1},
for every $w\in Da$ the $D$-conjugacy class coincides with $B$-conjugacy class. 
This implies the equality $2\cdot |C_{D}(w)|= | C_{B}(w)|$ as $|B:D|=2$.
So the result follows. 

(4) Set $\overline{G}=O(V)$. The natural embedding $G\ra H=SO_{2n+1}(q)$ extends to the embedding  $\overline{G}=O^\pm _{2n}(q)\ra H$
in an obvious way. Let 
 ${\mathbf H}=SO_{2n+1}(\overline{F}_q)$. 
Then the above embedding extends to an embedding $e: O_{2n}(\overline{F}_q)\ra {\mathbf H}$. 
Observe that $e(N_{O_{2n}(\overline{F}_q)}({\mathbf T}))=N_{\mathbf H}({\mathbf T})$. Indeed, $e({\mathbf T} )$
fixes a unique a  singular line $U$ on  $\overline{V}$ (the natural ${\mathbf H}$-module), and hence $N_{\mathbf H}({\mathbf T})$ stabilizes  $U$. 
As the stabilizer of $U$ in ${\mathbf H}$ coincides with  $e(O_{2n}(\overline{F}_q))$, the claim follows.
Recall that $W(T)=(N_{\mathbf G}({\mathbf T})/{\mathbf T})^{Fr}$ and $W_H(T)=(N_{\mathbf H}({\mathbf T})/{\mathbf T})^{Fr}$.  Clearly, $(N_{\mathbf H}({\mathbf T})/{\mathbf T})^{Fr}\cong N_{ G}({\mathbf T})/T$.  
If $G=SO^+_{2n}(q)$ then 
 $(4)$ follows in this case by comparison of  the items (1) and (2) above. If $G=SO^-_{2n}(q)$ then (3) implies that $N_{ G}({\mathbf T}) \not\subset N_{\overline{ G}}({\mathbf T}) $, which yields the result.

\bl{tt5}  $(1)$ Let $G=SO(V)\cong  SO_{2n}^-(q) $. Then the functions $i\ra d_i,j\ra e_j$, where $l=\sum _j e_j$ is odd,  determine a maximal torus in $G$ up to conjugacy. In particular, G has no neutral maximal torus. 

$(2)$ Let $G=SO(V)\cong  SO_{2n}^+(q) $. Then the functions $i\ra d_i,j\ra e_j$, where $l=\sum _j e_j$ is even, determine a maximal torus in $G$ up to conjugacy, unless this is 
exceptional. In the latter  case there are two non-conjugate maximal tori corresponding to the function $i\ra d_i$. 
\el 
Let $T,T'$ be maximal tori in $G$ corresponding to 
the same pair of functions $i\ra d_i,j\ra e_j$. Let ${\mathbf T},{\mathbf T}'$ be maximal $Fr$-stable tori of ${\mathbf G}$
such that $T={\mathbf T}^{Fr}, T'={\mathbf T}^{\prime Fr}$. By Lemmas \ref{my6} and \ref{my7}, $g{\mathbf T}g\up= {\mathbf T}'$ for some 
$g\in O(V)$. If ${\mathbf T}$ is not exceptional, then $N_{O(V)}({\mathbf T})$ is not contained in $SO(V)$ (Lemma \ref{a33}(4)), so there is $x\in O(V)\setminus SO(V)$ such that $x{\mathbf T}x\up={\mathbf T}$. If $g\notin SO(V)$ then $gx\in SO(V)$, and $gx{\mathbf T}(gx)\up={\mathbf T}'$. So in this case the lemma follows, in particular, (1) holds (see Definition \ref{de4}). 
Suppose that $T$ is exceptional. Then $V$ is of Witt defect 0, and
 $SO(V)\cong  SO_{2n}^+(q) $. 
Let ${\mathcal S}$ denote the $O(V)$-orbit of ${\mathbf T}$.
If ${\mathcal S}$ coincides with the $G$-orbit of ${\mathbf T}$  then $N_{O(V)}({\mathbf T})\neq N_{G}({\mathbf T})$. This is not the case by    Lemma \ref{my7}(4). 

\med Remark. Lemma \ref{tt5} is equivalent to saying that an $O(V)$-conjugacy class of maximal tori in $SO(V)$ consists of a single $SO(V)$-conjugacy class if the tori of ${\mathcal S}$ are non-exceptional, otherwise it  consists of two $SO(V)$-conjugacy classes.

\section{The Steinberg-plus character and its Curtis dual}


For  a finite reductive group ${\mathbf H}$ one can define an \ir character $St_H$ called the Steinberg character. We use for $St_H$ the definition in \cite[9.3]{DM}, which simultaneously gives the values of $St_H$.
Specifically,  $St_H$ is defined to vanish at all non-semisimple elements, 
whereas for the semisimple elements $h\in H$ one defines 
\begin{equation}\label{eq1}
 St_H(h)=\ep_{{\mathbf H}}\ep_{C_{{\mathbf H}}(h)}\cdot |C_H(h)|_p,\end{equation}
where $C_H(h)$ is the centralizer of $h$ in $H$,
$|C_H(h)|_p$ is the $p$-part of the order $|C_H(h)|$ of
$C_H(h)$ and $p$ is the defining characteristic of ${\mathbf H}$.   For the meaning of $\ep_{{\mathbf H}},\ep_{C_{{\mathbf H}}(h)}$,  see Notation.

Let $G=SO_{2n-1}(q)$.  A priory, $St_H|_G$ may depend on the choice of $SO^+ _{2n}(q)$ or $SO^- _{2n}(q)$ for $H$. 
We show that, in fact,  $St_H|_G$ does  not depend on the choice of $H\in \{SO^\pm _{2n}(q)\}$. 
Furthermore, if $q$ is odd, there are two non-equivalent embeddings $G\ra H^+$ (as well as $G\ra H^-$) obtained as follows. 
Let $H=SO(V)$  and let $Q$ be the quadratic form defining $H$. Let $v,v'\in V$ be non-singular vectors such that $Q(v)/Q(v')$ is a non-square in $F_q$. Then $G\cong {\rm Stab} _{H}(v)\cong {\rm Stab} _{H}(v')$, but the groups $ {\rm Stab} _{H}(v)$, $ {\rm Stab} _{H}(v')$ are not conjugate in $H$. Therefore, $St_H|_G$ may depend on the choice ${\rm Stab} _{H}(v)$ or ${\rm Stab} _{H}(v')$ for $G$. In fact, this is not the case. 
  This means that the character $St_H|_G$ depends on $G$ only,  and we shall call it the {\it Steinberg-plus character} and denote by $St^+_G$.  We keep this term for the restriction  $St_H|_G$ when  $G=SO^\pm _{2n}(q)$ and $H=SO_{2n+1}(q)$ (in this case there is no ambiguity in the definition). 

\begin{propo}\label{wd1}  Set $G=SO_{2n-1}(q)$,  $H^+=SO^+_{2n}(q)$ and $H^-=SO^-_{2n}(q)$.
Then $St_{H^+}|_G=St_{H^-}|_G$. In addition, both the left and the right hand sides do not depend on the choice of a natural embedding $G\ra H^+$ and $G\ra H^-$. \end{propo}

Proof. It suffices to show that the equality holds at the semisimple elements of $G$ as both the characters vanish at all non-semisimple elements.

Recall that the Steinber character of any finite reductive group $X$ is of defect 0, and hence remains \ir under reduction modulo $p$. Furthermore, the \rep obtained is the restriction to $X$ of a \rep $\tau$, say, of the algebraic group ${\mathbf X}$ \cite[Theorem 43]{St}. We apply this to  $X=H^+$ and $X=H^-$, and observe that ${\mathbf X}$ is the same in both the cases. (One can assume that ${\mathbf X}=SO_{2n}(\overline{F}_q)$ and ${\mathbf G}=SO_{2n}( \overline{F}_q))$. In addition, the \f diagrams are commutative:

\begin{center}
$\,\,~G\ra H^+$ $~~~~~~~~~~~~~~~~~~~~~~~~G\ra H^-$

$\,\downarrow$ $~~~~\downarrow$ $\,\,~~~~~~~~~~~~~~~~~~~~~~~~~\downarrow$ $~~~~\downarrow$

$\,{\mathbf G}\ra {\mathbf X}$  $~~~~~~~~~~~~~~~~~~~~~~~~~~{\mathbf G}\ra {\mathbf X}$
\end{center}
Therefore, $G\ra H^+\ra {\mathbf X}$ and $G\ra H^-\ra {\mathbf X}$ yield the same embedding $G\ra  {\mathbf X}$. 
Furthermore, the groups ${\rm Stab} _{H}(v)$ and ${\rm Stab} _{H}(v')$ (where $v,v'\in V$ are chosen as in the paragraph prior the lemma) are conjugate in  $SO_{2n}(\overline{F}_q)$. Therefore, $\tau|_G$ is independent from the choice 
of the embedding $G\ra H$, and lemms follows. 

\med
We denote by
$\om_{G}$ the Curtis dual of the character $St^+_G=St_H|_G$
\cite[8.8]{DM} or \cite[Ch.8, \S 2]{C}. So $\om_{G}$ is a generalized
character of $G$. 

\bl{ax2}  Let $\Phi$  be a class function on $G$ vanishing at all non-semisimple elements, and 
let $\phi$ be the Curtis dual of $\Phi$. Let $g\in G$. Then
$\phi(g)=\Phi(g)/St_G(g)$ if $g$ is semisimple.   In general,  if $g=su$, where $s$ is semisimple, $u$ is unipotent and $su=us$ then $\phi(g)=\phi(s)$.
\el
This  follows from Brou\'e \cite[Theorem 1 and \S
2${}^{\circ}$]{Br}, see also \cite[Proof of Lemma 4.1]{HZ2}, where
this result is deduced from Carter \cite[7.6.4]{C}.

As $St_H|_G$ vanishes at all non-semisimple
elements of $G$, it follows that
$\om_{G}(g)=\frac{St_H(g)}{St_G(g)}$ for all semisimple elements
$g\in G$ (see \cite{Br} or \cite[Lemma 4.1]{HZ2}). Furthermore, if $g=su$, where $s$ is semisimple, $u$
is unipotent and $su=us$ then $\om_{G}(g)=\om_{G}(s)$.

\begin{lemma}\label{hk8}  Let 
V be the natural $F_qG$-module.
 Then $\om_{G}(g)=\pm \om_{SO (V^g)}(1)=\pm q^{[\dim V^g/2]}$.
In particular, $\om_G(1)=q^{[\dim V^g/2]}$.
\end{lemma}

Proof. Suppose first that $g=1$. If $\dim V=2n$ then $St_G(1)=|G|_p=q^{n(n-1)}$ whereas $St_H(1)=|H|_p=q^{n^2}$. If $\dim V=2n+1$ then $St_G(1)=|G|_p=q^{n^2}$ whereas $St_H(1)=|H|_p=q^{n(n+1)}$. In both the cases $\om_G(1)=q^n$. as requireded. In general, let $V\subset W$, where $W$ is an orthogonal space of dimension $1+\dim V$ and  $H=SO(W)$.
We can write $V=V^g\oplus (1-g)V$ and $W=W^g\oplus (1-g)W$. Note that $(1-g)V=(1-g)W$. Set $g'=g|_{(1-g)V}$. As $g$ fixes no non-zero vector in $(1-g)V$
and $(1-g)W$, we have $|C_G(g)|_p=|SO(V^g)|_p\cdot |C_{SO((1-g)V)}(g')|_p$ 
(at least for $p>2$), and $|C_H(g)|_p=|SO(W^g)|_p\cdot |C_{SO((1-g)V)}(g')|_p$. 
 Therefore, $\om_{G}(g)=\frac{St_H(g)}{St_G(g)}=\pm \frac{|SO(W^g)|_p}{|SO(V^g)|_p}=
\pm \om_{SO(V^g)}(1)$, and the lemma follows. 

\med
Examples. Here we consider some degenerate cases where
$n=1,2,3$.

\med
(1) $n=1$. Then $G=\{1\}$ and $H=SO_2^\pm (q)$. Therefore,
$\om_{G}=St_H(1)\cdot 1_G$. As $St_H(1)=|H|_p$, we have $St_H(1)=1$, so
$\om_{G}= 1_G$.

\med
(2) $n=2$. The group $G=SO^{\pm }_2(q)$ is abelian of order $q\mp 1$. In particular, $G$
has no non-trivial $p$-element, so $|C_G(g)|_p=1$ for $1\neq g\in G$.
It follows that $St_G=1_G$.
Therefore, $\om_{G}=St_G^+$. Recall that $H=SO_3(q)$ and $St_H(1)=q$. One observes that
$C_H(g)$ is a $p'$-group
for every $1\neq g\in G$, so $St_H(g)=\pm 1$. Then  
the $F_q$-rank of
$SO_2^+ (q)$ equals 1, while the $F_q$-rank of
$SO_2^- (q)$ equals 0. As $G$ is abelian,  $C^0_{{\mathbf H}}(g)\cong {\mathbf G}$
 for every $1\neq g\in G$. So for $1\neq g\in G$
 we have

$$\om_{G}(g)=St_H(g)=\begin{cases}1&if ~G= SO_2^+(q)\cr
-1&if ~G= SO_2^-(q).\end{cases}$$
One easily deduces from this, that $\om_{G}=\rho_G^{reg}+1_G$
if $G= SO_2^+(q)$, and $\rho_G^{reg}-1_G$
if $G= SO_2^-(q)$.

\med
(3) $n=3$. Here  $G= SO_3(q)\cong PGL_2(q)$, and $H= SO_4^+(q)$.

Group $G$ has two maximal tori $T_1,T_2$ of order $q-1$ and $q+1$,
respectively. Let $V$ be the natural module for $H$. One observes that if
$1\neq t\in G$ is a semisimple then $(t-1)V$ and $ V^t$ are
2-dimensional non-degenerate subspaces of $V$ orthogonal to each
other, and $V=(t-1)V\oplus V^t$. It follows from this that
$|C_G(t)|_p=1$. Furthermore, it is easy to observe that
$St_H(t)=1$, respectively,  $-1$ if $t\in T_1$, respectively, $t\in T_2$. 
 In addition, $C_G(t)=1$, respectively,  $-1$ if
$t\in T_1$, respectively, $t\in T_2$. 
In fact $\om_{G}(g)=1$, respectively,
$-1$ if $t\in T_1$, respectively. $t\in T_2$. (Note that $\om_{G}(u)=q$
for every unipotent element $u\in G$, in particular, $\om_{G}(1)=q$.)

\subsection{Multiplication theorem}

\begin{lemma}\label{hk9}  Let $V$  be an orthogonal space,
and  $V=V_1\oplus V_2$, where  $V_1,V_2$ are non-degenerate
subspaces of $V$ orthogonal to each other. Set $G_1=SO(V_1)$, $G_2=SO(V_2)$. Let $g\in G$
be a semisimple element such that  $gV_i=V_i$, 
and set $g_i=g|_{V_i}$ for $i=1,2$. 

\med
$(i)$ Suppose that not both $\dim V_1,\dim V_2$ are odd. Then
$\om_G(g)=\pm \om_{G_1}(g_1)\cdot \om_{G_2}(g_2)$.

\med
$(ii)$ Suppose that both $\dim V_1,\dim V_2$ are odd. Then
$\om_G(g)=\pm q\cdot \om_{G_1}(g_1)\cdot
\om_{G_2}(g_2)$.\end{lemma}

Proof. 
By Lemma \ref{hk8}, $\om_G(g)=\pm q^{[(\dim V^g)/2]}$. Note that
$V^g=V_1^g\oplus V_2^g$. So $\om_G(g)=\pm q^{[(\dim V^g_1+ \dim
V^g_2)/2]}$. It is well known that $\dim V-\dim V^g$ is even whenever $g\in SO(V)$. So the result follows from 
the \f observation.
 If $a,b\geq 0$ are integers then $[a/2]+[b/2]=\frac{a+b}{2}-1$
if $ab$ is odd, otherwise
$[a/2]+[b/2]=[(a+b)/2]$.  

\med
The \f lemma is a key ingredient of our proof of Theorem \ref{mm8}

\begin{lemma}\label{mt6} 
Let $V$  be an orthogonal space, and  $V=V_1\oplus V_2$, where  $V_1,V_2$ are non-degenerate.

\med
$(i)$ Suppose that at least one of $\dim V_1,\dim V_2$ is even.
Then there is a linear character $\lam$ of $G_1G_2$ such that
$\lam^2=1$ and $\om_G(g)=\lam(g) \cdot \om_{G_1}(g)\cdot
\om_{G_2}(g)$ for every semisimple element $g\in G_1G_2$.

\med
$(ii)$ Suppose that both $\dim V_1,\dim V_2$ are odd. Then there
is a linear character $\lam$ of $G_1G_2$ such that $\lam^2=1$ and
$\om_G(g)=\lam(g) \cdot q\cdot \om_{G_1}(g)\cdot \om_{G_2}(g)$ for every semisimple element $g\in G_1G_2$. 

In addition, $\lam(u)=1$ for all unipotent elements $u\in G_1G_2$.\el
 Applying Theorem \ref{ct1} to $X=G_1G_2$, $\phi=\om_G|_{G_1G_2}$ and to
 $\psi=\om_{G_1}\cdot \om_{G_2}$, we have the equality required at the semisimple elements $g.$
 Lemma 
\ref{hk9}
 guarantees that the hypothesis of
Theorem \ref{ct1} holds. 

Furthermore, let $\lam'$ be the  linear character of $G_1G_2$ such that 
$\lam(s)=\lam'(s)$ at the semisimple elements and $\lam'(y)=1$ for all unipotent elements of $G_1G_2$. 
Let   $x\in G_1G_2$ and $x=gu$, where $g$ is semisimple, $u$ is unipotent and $gu=ug$. We know that $\om_G(x)=\om_G(g)$ and $\om_{G_i}(x)=\om_{G_i}(g)$ for $i=1,2.$ So in (1) we have $\om_G(x)=\om_G(g)=
\lam(g) \cdot \om_{G_1}(g)\cdot \om_{G_2}(g)=\lam(g) \cdot \om_{G_1}(x)\cdot \om_{G_2}(x)=\lam'(x) \cdot \om_{G_1}(x)\cdot \om_{G_2}(x)$, and similarly in (2).

\med Remarks. (1) If $\dim V_2=1$ then $\om_{G_2}=1_{G_2}$, so
$\om_G(g)=\lam(g) \cdot \om_{G_1}(g)$. If $\dim V_2=2$ then $G_2$
is a torus, but $\om_{G_2}\neq 1_{G_2}$. (2) If one replaces
$G,G_1,G_2$ by the
subgroup $G',G_1',G_2'$ of index 2 in $G,G_1,G_2$, respectively,   
 then $\lam|_{G_1'G_2'}=1_{G_1'G_2'}$.

\bl{vg5}  The character $\lam$ in Lemma $\ref{mt6}$ is
trivial. \el
If $q$ is even then semisimple elements of $G$ are of odd order. As $\lam^2=1$ and $\lam$ is trivial at every unipotent element, we conclude that $\lam$ is trivial.

Let  $q$ be odd. Suppose the contrary. Then $\lam_i\neq 1$ for $i\in \{1,2\}$. Fix any $i$ with $\lam_i\neq 1$.

Suppose first that $\dim V_i$ is even for $i=1,2$.
Let $G=G_1\times G_2$, and let $\lam$ be a character of $G_1G_2$ of
order 2. Then $\lam=\lam_1\lam_2$, where $\lam_i=\lam|_{G_i}$.

Note that $G_i$   has a unique subgroup of index 2 \cite[2.5.7]{KL}. As
$\Omega(V_i)$ has index 2 in $G_i$, it follows that if
$\lam_i(G_i)\neq 1$ then $\lam_i(x)=1$ $(x\in G_i)$ \ii $x\in \Omega(V_i)$.
By Lemma \ref{h00}, there is $t\in SO(V_i)\setminus \Omega(V_i)$ such that 
$|C_{G_i}(t)|_p=|C_{G_i}(t^2)|_p$ and $V^t_i=V_i^{t^2}$. View $t$ as an element of $G$. Then $\om_G(t)=\lam_i(t)\om_{G_i}(t)=-\om_{G_i}(t)$ as $\lam_i^2=1$. By Lemma \ref{h03},  
 $|C_{G}(t)|_p=|C_{G}(t^2)|_p$ and $|C_{H}(t)|_p=|C_{H}(t^2)|_p$; in addition $C_{{\mathbf G}}(t)^0=
C_{{\mathbf G}}(t^2)^0$ and $C_{{\mathbf H}}(t)^0=
C_{{\mathbf H}}(t^2)^0$. Therefore,   $St_H(t)=St_H(t^2)$ and $St_G(t)=St_G(t^2)$. This yields  $\om_{G}(t)=St_H(t)/St_G(t)=\om_{G}(t^2)$. Then $-\om_{G_i}(t)=\om_G(t)=\om_{G}(t^2)=\lam_i(t^2)\om_{G_i}(t^2)=\om_{G_i}(t^2)=\om_G(t)$, which is a contradiction.

\begin{corol}\label{vv3}  Let  V be an orthogonal space of odd dimension n
over $F_q$,  and let $V'$ be a non-degenerate subspace of dimension $n-1$ in V. 
 Set $G=SO(V)$ and $X=SO(V')$
if q is odd, otherwise set  $G=\Omega(V)$ and $X=\Omega(V')$. View $X$ as a subgroup of $G$. 
Then $\om_G|_X=\om_X$.

\end{corol}
Proof. If $q$ is odd  then 
the corollary straightforwardly follows from Lemmas \ref{mt6} and \ref{vg5}.
Suppose that $q$ is even. In this case we have $V=V_1+V'$, where $\dim V_1=1$ 
however, $V_1$ is not non-degenerate. Nonetheless, the formula $\om_G(g)=\pm q^{[(1+\dim V^g_2)/2]}=\pm q^{\dim V^g_2/2} $ remains true as $\dim V^g_2$ is even. This implies $\om_G(g)=\pm \om_X(g)$
for $g\in X$. The reasoning in the proofs of Lemmas \ref{mt6} and \ref{vg5}
remain valid, whence the result. 

\med
{\bf Proof of Theorem} \ref{mm8}. The result follows from Lemmas \ref{mt6} and \ref{vg5}.

\med
The standard embedding $GL_n(q)\ra SO^+_{2n}(q)$ is well known. This
comes from fixing a maximal totally singular subspace $W$ of $V$.
Every basis $b_1\ld b_n$ of $W$ can be extended to a basis $b_1\ld
b_{2n}$ of $V$ such that the Gram matrix of this basis is
$\begin{pmatrix}0&\Id_n\\ \Id_n&0\\ \end{pmatrix}$. Then the
matrix $\diag(g,{}^tg\up)$ belongs to $SO^+_{2n}(q)$ for every $g\in
GL_n(q)$. (Here ${}^tg$ denotes the transpose of $g$.)

If $n$ is even (respectively, odd) then there is an embedding $U(n,q)$ into $O^+_{2n}(q)$ (respectively,  $O^-_{2n}(q)$).   
 These can be obtained as follows.  Let $W$ be the natural $F_{q^2}G$-module for $H=U_n(q)$. 
 There is a surjective embedding of
$h:W\ra V $, where $V$ is an orthogonal $F_q$-space  of dimension
$2n$,  such that every non-degenerate one-dimensional subspace of $W$ goes
to a two-dimensional   anisotropic  subspace. Moreover, $h$ regards
the orthogonality relation, that is, if subspaces $X,Y$ of $W$ are
orthogonal then so are $h(X),h(Y)$. Therefore, $V$ is an
orthogonal  sum of two-dimensional anisotropic  subspaces. It
 follows  from \cite[2.5.11]{KL} that the Witt defect of $V$ equals
$0$ if $n$ is even, otherwise equals $1$. Therefore, $h$ yields an
embeddings $e:U_n(q)\ra O^+_{2n}(q)$ if $n$ is even, and $U_n(q)\ra
O^-_{2n}(q)$ if $n$ is odd. In addition,  $V^{e(g)}=h(W^g)$  and $2\dim W^g=\dim
h(W^g)=\dim V^{e(g)}$ for $g\in U_n(q)$. 
In fact,  $e(U_n(q))$ is contained in $G$. This is trivial if $q$ even, moreover,  in 
this case $e(U_n(q))$ is contained in $\Omega^\pm_{2n}(q)$ as $|U_n(q):SU_n(q)|$ is odd.
Let $q$ be odd. Note that $U_n(q)$ contains a central element $z$, say, of order greater than 2.
Then $e(z)$ does not have \ei $\pm 1$, and hence the claim follows from Lemma \ref{y9}.

Furthermore,  if $t,n$ are odd
then there is an embedding $U(n,q^t)$ into $U_{tn}(q)$ \cite{Hu}. In
particular,  $U_1(q^t)$ embeds into $U_t(q)$ and hence into  
$SO^-_{2t}(q)$.  

\med
 We use Gerardin's
definition of the Weil \rep of a unitary group. This is the one
afforded by the character, which is defined for   $g\in U_m(q)$ by
$\chi(g)=(-1)^m(-q)^{\dim W^g}$. 

\bl{uu1}  $(1)$ Let $G=SO_{2n}^+(q)$ and $X\subset G$ a standard subgroup isomorphic to
$GL_n(q)$.   Then  there is a
linear character $\lam$ of $X$ such that  $\om_G(g)=\lam(g) \pi_X(g)$ for all semisimple elements $g\in X$, where $\pi_X$
is the character of the permutation \rep of $X$ associated with the action of
X on the vectors of $F_q^n$.

\med
$(2)$ Let $G=SO_{2n}^+(q)$   or  $SO_{2n}^-(q)$. Let
$Y\subset G$ be a subgroup isomorphic to
$U_n(q)$. 
Then there is a
linear character $\lam$ of $Y$ such that
$\om_G(g)=\lam(g)\cdot \phi(g)$ for all semisimple elements $g\in Y$, where $\phi$ is
the character of  the  Weil \rep of $Y$.
\el (1) By the formula for $g\in G$ in by Lemma \ref{hk8} we have  $\om_G(g)=\pm q^{\dim
V^g/2}$.
 As the Witt index of $V$ is equal to $n$,
there is a basis $b_1\ld b_{2n}$ in $V$ with Gram matrix
$\begin{pmatrix}0&\Id_n\\ \Id_n&0\end{pmatrix}$. Then the elements
$g:=\diag(x,{}^Tx\up)$ with $x\in GL_n(q)$ preserve the Gram
matrix so $g\in SO^+_{2n}(q)$. Let $U_1=\lan b_1\ld b_n\ran$,
$U_2=\lan b_{n+1}\ld b_{2n}\ran$. Then $V^g=U^g_1+U^g_2$.
Therefore, $\om_G(g)=\pm q^{\dim U_1^g}$. Note that
$\pi(g)=q^{\dim U_1^g}$. So (1) follows from the comparison
theorem.

\med
(2) Let $g\in U_n(q)$ be a semisimple element.  As  the Weil character value at  $g$ is  
$(-1)^n(-q)^{\dim W^g}$ and $\om_G(e(g))=\pm q^{\dim V^{e(g)}/2}$ (Lemma \ref{hk8}), 
 the result again follows  by the comparison
theorem. 

\bl{ir1}
 The character $\lam$ in Lemma $\ref{uu1}$ is
trivial.
\el If $n>2$ then the derived subgroup $G'$ of $G=SO^\pm
_n(q)$ has index 2, so $\lam^2=1_G$. We show
that $\lam=1_G$. Suppose the contrary.

\med
(a) $X=GL_n(q)$.  Let $t\in G$ be of order $q^n-1$. Then it is
known that $t\notin G'$, and hence $\lam (t)=-1$ and  $\lam
(t^2)=1$. 
As
both $t$ and $t^2$ are regular semisimple in $G$ and in $H$, we have
$C_H(t)=C_H(t^2)$ and $C_G(t)=C_G(t^2)$. Therefore, we have $$\om
_{G}(t)=St_H(t)/St_G(t)=
\frac{\ep_{\mathbf H}\ep_{C^0_{\mathbf H}(t)}|C_H(t)|_p}{\ep_{\mathbf G}\ep_{C^0_{\mathbf G}(t)}|C_G(t)|_p}
=\om _{G}(t^2).$$ Now $\pi (t)\lam(t)=\om _{G}(t)=\om
_{G}(t^2)=\pi(t^2)=\pi(t),$ a contradiction.

\med
(b) $X=U_1(q^n)$. Let $t\in G$ be of order $q^n+1$. Then it is
known that $t\notin G'$, and hence $\lam (t)=-1$ and  $\lam
(t^2)=1$. 
Let $\phi$ denote the Weil \rep of $U(1,q).$ Then
$\phi(t)=-(-q^n)^0=\phi(t^2)=-1$. As both $t$ and $t^2$ are
regular semisimple in $G,H,$ we have $C_H(t)=C_H(t^2)$ and
$C_G(t)=C_G(t^2)$. Therefore, we have $$\om
_{G}(t)=St_H(t)/St_G(t)=
\frac{\ep_{\mathbf H}\ep_{C^0_{\mathbf H}(t)}|C_H(t)|_p}{\ep_{\mathbf G}\ep_{C^0_{\mathbf G}(t)}|C_G(t)|_p}
=\om _{G}(t^2).$$ Now $\phi (t)\lam(t)=\om _{G}(t)=\om
_{G}(t^2)=\phi(t^2)=\phi(t)$, a contradiction.

\section{Maximal tori and Curtis dual}

\subsection{Character formula}



In this section we give a formula that describes in a convenient
way the restriction $\om_{G}|_{T}$. Recall that maximal tori in $G$ are determined up to $O(V)$-conjugation by the orthogonal decompositions $(\ref{eq2})$, which in turn yields an "orthogonal" decomposition $(\ref{eqT})$ for the corresponding torus. The formula depends on the $O(V)$-conjugacy class of a torus in question, rather than on the $SO(V)$-conjugacy class. This is not surprizing as $SO_{2n+1}(q)$ contains $O^\pm_{2n}(q)$, and hence $St^+_G$ is invariant under $O(V)$ for $G=SO(V)=SO^\pm_{2n}(q)$.

\med
 First, we deal with two
special cases:

\bl{mt2}  $(1)$ Let $G=SO_{2n}^+(q)$ and $T$ a  maximal torus  of
order $q^n-1$ (that is, T corresponds to an orthogonal  decomposition with $k=1,l=0$).
Then $\om_G|_T=\rho_T^{reg}+1_T$.

\med
$(2)$ Let $G=SO_{2n}^-(q)$ and $T$ a   maximal torus  of order $q^n+1$ (that is, T corresponds
to an orthogonal   decomposition with $k=0,l=1$). Then $\om_G|_T=\rho_T^{reg}-1_T$.
\el
 This follows from Lemma \ref{uu1}.

\medskip
Next, as in \cite{HZ09},  Theorem \ref{mm8}   allows us to express the
character of $\om_G|_T$ in terms of the characters
$\om_{G_i}|_{T_i}$. In turn,    $\om_{G_i}(T_i)$ can be  expressed
in terms of the regular character $\rho_{T_i}^{reg}$ and the
trivial character $1_{T_i}$. (For a group $H$ we denote by
$\rho^{reg}_H$ the character of the regular \rep of $H$.)

\begin{theo}\label{g1}  
Let $G=SO(V)$ if q is odd, and $G=\Omega(V)$ if q is even.  Let $T=T_1\otimes \cdots \otimes  T_k\otimes  T_{k+1}\otimes  \cdots \otimes T_{k+l}$
be a maximal torus of G as in $(\ref{eqT})$. Then 
 $$\om_G|_T = 
(\rho^{reg}_{T_1}+1_{T_1})\otimes\cdots
\otimes(\rho^{reg}_{T_k}+1_{T_k})\otimes (\rho^{reg}_{T_{k+1}}
-1_{T_{k+1}})\otimes\cdots
\otimes(\rho^{reg}_{T_{k+l}}-1_{T_{k+l}}) .$$
\end{theo}

Proof. If $\dim V$ is even, the result follows from 
Theorem \ref{mm8} and Lemma  \ref{mt2}. Suppose that  $\dim V$ is odd. Then,  by Lemma \ref{fv1}, $T$ stabilizes a non-degerate subspace  $V'$, say, of $V$.  So $T\subset X$, where $X\cong SO(V)$ if $q$ is odd 
and $X\cong \Omega(V')$ if $q$ is even. 
So the result follows from that for the case where $\dim V$ is even.

\begin{corol}\label{c22} 
The restriction of $\om_{G}$ 
to every maximal torus $T$ of $G$ is a proper character of T.
\end{corol}

\med
Examples. (1) $G=SO_3(q)$. We have two maximal tori $M_1,M_2$ of order $q-1,q+1$, respectively.
We have $\om_G|_{M_1}=
\rho_{M_1}^{reg}+1_{M_1}$ and $\om_G|_{T_2}=
\rho_{M_2}^{reg}-1_{M_2}$. 
 
\med
(2) $G=SO^-_4(q)$, $q$ odd.
Up to conjugacy, there are two maximal  tori $M_1,M_2$ of order
$q^2-1,q^2+1$, respectively.  They correspond to the decompositions with
$k=l=1$ and $k=0,l=1$, respectively. (That is, $\dim V_1=\dim V_2=1$ in (\ref{eq2}) in  the former case, 
and  $\dim V_l=4$ with $l=1$ in the latter case. Note that $M_1$ is not a cyclic group.)  So $M_1=T_1\times T_2$, where $|T_1|=q-1$, $|T_2|=q+1$.  
We have $\om_G|_{M_1}=(\rho_{T_1}^{reg}+1_{T_1})\otimes (\rho_{T_2}^{reg}-1_{T_2})$ and
$\om_G|_{M_2}=\rho_{M_2}^{reg}-1_{M_2}$. (Note that $G\cong PSL_2(q^2)$ if $q$ is even,
 otherwise $G\cong  PSL_2(q^2)\times \lan -\Id\ran$.  Indeed, $\Omega^-_4(q^2)\cong  PSL_2(q^2)$ \cite[2.9.1]{KL},
and $-\Id\in SO^-_4(q^2)\setminus  \Omega^-_4(q^2)$. As $|SO^-_4(q^2):  \Omega^-_4(q^2)|=2$, the claim follows.)

\med
(3) $G=SO^+_4(q)$. We have four
maximal tori $M_1,M_2,M_3,M_4$ of orders
$(q-1)^2,(q+1)^2,q^2-1,q^2-1$, respectively. They correspond to the
decompositions with $(k,l)=(2,0)$, $(0,2)$, and remaining two tori correspond   to $(k,l)=(1,0)$,
that is, $\dim V_k=4$ for $k=1$. We have
$\om_G|_{M_1}=(\rho_{T_1}^{reg}+1_{T_1})(\rho_{T_2}^{reg}+1_{T_2})$,
$\om_G|_{M_2}=(\rho_{T_1}^{reg}-1_{T_1})(\rho_{T_2}^{reg}-1_{T_2})$,
and $\om_G|_{M_3}=\om_G|_{M_4}=\rho_{T_1}^{reg}-1_{T_1}$ in notation of Theorem \ref{g1}. (Note that $G\cong SL_2(q)\times SL_2(q)$ if $q$ is even,
 otherwise $G\cong \cong (SL_2(q)\circ SL_2(q))\cdot 2$.)
 
 \med
(4) Let $G=SO^\al_{2n}(q)$, $q$ odd. Then $\om_G(-\Id)=\al$. Indeed, if $\al=1$ then $-\Id$ belongs to a maximal torus with $k=n$, $l=0$,  and 
the claim follows from Theorem \ref{g1} as $l=0$. Alternatively, one can use Lemma \ref{tt5}(1). If $\al=-1$ then 
$-\Id$ belongs to a maximal torus with $k=n-1$, $l=1$,  and again apply  
  Theorem \ref{g1}.

\begin{lemma}\label{g2}  Let G and $\nu$ be as in Theorem $\ref{g1}$.
Let $\theta = \theta_1 \otimes \cdots \otimes \theta_{k+l}$ be an
irreducible character of $T = T_1 \times \cdots \times T_{k+l}$.

\medskip
$(1)$ If $\theta_{k+j} =1_{T_{k+j}}$ for some $j>0$ then
$\theta$ does not occur as an irreducible constituent of
$\om_{G}|_T$ (that is, $(\om_G|_T,\theta )=0$).

\medskip
$(2)$ Suppose that $\theta_{k+j} \neq 1_{T_{k+j}}$ for every
$j=k+1, \ldots , l$. Let $k(\theta)$ be the number of $0 \leq i \leq
k$ such that $\theta_i=1_{T_i}$. Then $(\om_G|_T,\theta )=
2^{k(\theta)}$.

\medskip
$(3)$ Suppose that $\theta_{i} \neq 1_{T_{i}}$ for every $1 \leq
i \leq k + l$. Then $(\om_G|_T,\theta )= 1$. 
\end{lemma}

Proof. This follows from Theorem \ref{g1}. 

\med
Lemma \ref{g2} can be restated in terms of the dual group $G^*.$  Note that
the dual of $SO^+_{2n}(q)$ (respectively, $SO^-_{2n}(q)$, respectively,
$SO_{2n+1}(q)$) is isomorphic to $SO^+_{2n}(q)$ (respectively,
$SO^-_{2n}(q)$, respectively, $Sp_{2n}(q)$), see \cite[p. 120]{C}. Let
$V^*$ denote the natural module for $G^*$, and we assume that 
 $V^*$ endows a bilinear or quadratic form defining $G^*$.

A torus $T^*$  of $G^*$  dual to $T$ has the same
structure as $T$, and in fact there is an orthogonal decomposition
$V^*=V_1^*\oplus \cdots \oplus V_k^*\oplus V_{k+1}^*\oplus\cdots
\oplus V^*_{k+l}$ with the same properties for $T^*$ as those for
$T$ on $V$. The only difference is that we do not have $V_0^*$
anymore. In particular, $\dim V_i^*=\dim V_i$ and $|T^*_i|=|T_i|$
for $i=1\ld k+l$. In addition, $T^*=T_1^*\times \cdots \times
T_{k+l}^*$, 
and  $T^*_i$ can be viewed as a
maximal torus of $SO(V_i^*)$. This also tells us that if
$\theta=\theta_1\otimes \cdots \otimes \theta_{k+l}$ is a linear
character of $T$ then the corresponding element $s\in T^*$ can be
expressed as $(s_1\ld s_{k+l})$, where $s_i\in T_i^*$ for $i=1\ld
k+l$. In addition, one observes that $\theta_i=1_{T_{i}}$ \ii
$s_i=1$.

Now, we can restate Lemma \ref{g2} as follows:


\begin{lemma}\label{d2g}  Let $s\in T^*$  be a semisimple
element and a maximal torus in $ G^*$, which correspond to a pair
$(T,\theta)$ under the duality $G\ra G^*$. Let $s =\diag (s_1,\ldots , s_{k+l})$,
where $s_i$ means the restriction of $s$ to $V_i^*$. 

\medskip
$(1)$ If $s_{k+j} =1$ for some $j>0$ then $(\om_G|_T,\theta )=0$ (that is, $\theta$ does not occur
as an irreducible constituent of $\om_{G}|_T)$.

\medskip
$(2)$ Suppose that $s_{k+j} \neq 1 $ for every $j=k+1, \ldots ,
l$. Let $m(\theta)$ be the number of $0 \leq i \leq k$ such that
$s_i=1$. Then $(\om_G|_T,\theta )= 
2^{m(\theta)}$.

\medskip
$(3)$ Suppose that $s_{i} \neq 1 $ for every $1 \leq i \leq k +
l$. Then $(\om_G|_T,\theta )=1$. 
\end{lemma}

\begin{corol}\label{nc4}  $(\om_G|_T,\theta )=0$ \ii the $1$-eigenspace of s on $V^*$ is non-zero and not
 contained in $V^*_1\oplus\cdots \oplus V^*_k$.
\end{corol}

\begin{corol}\label{hh1}  Let $s\in G^*$, and let $T^*$ be a maximal torus of $ G^*$ containing $s$. 
Let $\theta\in \Irr T$ corresponds to $s\in T^*$. Suppose that    $s$ does not have \ei $1$ on the natural module
for $G^*$. Then $(\om_G|_T, \theta)=1$. 
\end{corol}

\begin{corol}\label{hh2}  Let $s\in G^*$, and let $T^*$ be a maximal torus of $ G^*$ containing $s$. 
Let $\theta\in \Irr T$ corresponds to $s\in T^*$. Let $V^*$ be the natural module for $G^*$, let  $V^*_s$ be the $1$-eigenspace of $s$ on $V^*_s$ and let  $V^*=V^*_1\oplus \cdots \oplus V^*_k\oplus V^*_{k+1}\oplus \cdots \oplus V_k\oplus V^*_{k+l}$ a $T^*$-decomposition of $V^*$. 
Suppose that  $(\om_G|_T, \theta)>1  $. Then  $V^*_s\neq 0$ and $V_s^*\subseteq V_1^*\oplus \cdots \oplus V^*_k$. 
In particular,   $V_s^*$ is of Witt defect $0$. 
\end{corol}


\subsection{The characteristic $2$ case}

Let ${\mathbf G}$ be a simple algebraic group of rank $r$ 
and let $\lam_1\ld \lam_r$ denote the fundamental weights of 
${\mathbf G}$. The \ir \reps of ${\mathbf G}$ are parametrized by the dominant weights $\lam$ of ${\mathbf G}$,
which are  linear combinations $a_1\lam_1+\cdots+a_r \lam_r$ with non-negative coefficients $a_1\ld a_r$.
We write $\phi_\lam$ for the \irr of ${\mathbf G}$ corresponding to $\lam$.
A dominant weight is called $q$-restricted if all $a_1\ld a_r$ do not exceed $q-1$. Let $Fr$ be a Frobenius endomorphism of ${\mathbf G}$ and $G={\mathbf G}^{Fr}$. If  $G$ is not of type ${}^2B_2(q)$, ${}^2G_2(q)$, ${}^2F_4(q)$ then the restriction $\phi_\lam$ to $G$ is \ir whenever $\lam$ is $q$-restricted. Moreover, the \ir \reps of $G$ are parametrized by the $q$-restricted dominant weights. So we denoted by $\beta_\lam$ the Brauer character of $\phi_\lam|_G$.

\begin{propo}\label{p99} Let q be even and $G=Sp_{2n}(q)\cong \Omega_{2n+1}(q)$. Then 
$St^+_G= \beta_{(q-1)\lam_n}\cdot St_G$.

\end{propo}

Proof. 
It suffices to observe that $\om_G(g)$ coincides with  $ \beta_{(q-1)\lam_n}(g)$ for every semisimple element $g\in G.$ This follows by comparison of Theorem \ref{g1} with the formula in \cite[Proposition 4.12]{HZ09}.

\med 

For ${\mathbf G}_1=SO_{2n}(\overline{F}_2) $ with $n\geq 4$,  let  $\mu_1\ld \mu_n$ denote the 
fundamental weights of  ${\mathbf G}_1$.  Denote by $\Delta_t$ the set of weights $\nu_1+2\nu_2+\cdots +2^{t-1}\nu_{t-1}$, where $\nu_i\in \{\mu_{n-1},\mu_{n}\}$ for $i=1\ld t$.   Let $\lam_1\ld\lam_n$ and  $\mu_1\ld \mu_n$ denote the 
fundamental weights of  ${\mathbf G}\cong SO_{2n+1}(\overline{F}_2)$. 

\bl{qq2}  Let ${\mathbf G}_1 \subset {\mathbf G}\cong SO_{2n+1}(\overline{F}_2) $ be the natural enbedding. 
Then the restriction to 
${\mathbf G}_1$ of the \irr $\phi_{(q-1)\lam_n}$ of ${\mathbf G}$ is the direct sum of the \reps $\phi_{\nu}$ when   $\nu$ runs over $\Delta_t$.\el
Let $\ep_1\ld \ep_n$ be the weights of  ${\mathbf G}$ introduced in \cite[Planchee III]{Bo}. The weights of $\phi_{\lam_n}$ are well known to be  $\pm\ep_1\ld \pm \ep_n$  in terms of $\ep_1\ld \ep_n$.  
   They form a single $W({\mathbf G})$-orbit. Viewing $W({\mathbf G}_1)$ as a subgroup of $W({\mathbf G})$, these weights form two $W({\mathbf G}_1)$-orbits. Note that a maximal torus of ${\mathbf G}_1$ remains a maximal torus in ${\mathbf G}$. \itf  
the restriction  $\phi_{\lam_n}|_{{\mathbf G}_1}$ consists of at most two composition factors. The weights of $\phi_{\mu_n}$ and $\phi_{\mu_{n-1}}$ are also known, and one easily deduces from this that the composition factors are exactly $\phi_{\mu_n}$ and $\phi_{\mu_{n-1}}$. Observe that the restriction  $\phi_{\lam_n}|_{{\mathbf G}_1}$ can be obtained from the restriction  $\phi_{\lam_n}|_{O_{2n}(\overline{F}_2)}$. As $\phi_{\mu_n}$ is not invariant under $O_{2n}(\overline{F}_2)$, it follows that
$\phi_{\lam_n}|_{O_{2n}(\overline{F}_2)}$ is  irreducible, and, by Clifford's theorem,  $\phi_{\lam_n}|_{{\mathbf G}_1}$ is the direct sum of two \ir constituents. This implies the lemma for $q=2$.  

 In general, let $q=2^t>2$. 
 Then $q-1=1+2+\cdots +2^{t-1}$. It follows that
$\phi_{(q-1)\lam_n}=\phi_{\lam_n}\otimes \phi_{2\lam_n}\otimes \cdots \otimes \phi_{2^{t-1}\lam_n}$.
Therefore, the restriction of $\phi_{(q-1)\lam_n}$ to ${\mathbf G}_1$ is the tensor product of the restrictions of $\phi_{2^{i}\lam_n}$ to ${\mathbf G}_1$ for $i=0\ld t-1$. As $\phi_{2^{i}\lam_n}|_{{\mathbf G}_1 }=\phi_{2^{i}\mu_{n-1}}\oplus \phi_{2^{i}\mu_n}$, the lemma follows. (Note that the \reps  $\phi_{\nu}$ for $\nu\in \Delta_t$   are \ir by  \cite[Theorem 43]{St}.)

\med
Let $n\geq 4$ and $G_1=SO_{2n}^+(q)$ or $SO_{2n}^-(q)$, so $G_1=\GC^{Fr}$ for a suitable Frobenius
 endomorphism $Fr$ of $\GC=SO_{2n}(\overline{F}_2)$. Denote by $st_{G_1}^+$ and $st_{G_1}$ the projective
 modules with the characters $St^+_{G_1}$ and $St_{G_1},$ respectively.

\begin{propo}\label{yy8}  $st_{G_1}^+=\sum_{\nu\in\Delta_t}\phi_\nu\otimes st_{G_1}$.
\end{propo}

Proof. It is well known that $\phi_\nu\otimes st_{G_1}$
 is a projective module (possibly decomposable).  By Corollary \ref{vv3}, the Curtis dual  $\om_{G_1}$ of $St^+_{G_1}$ coincides with $\om_G|_{G_1}$.  
By Lemma \ref{qq2} and the above remarks, $\beta_{(q-1)\lam_n}|_{G_1}$ coincides with the sum of 
$\beta_{\nu}$ for all $\nu\in\Delta_t$. Therefore, $St^+_{G_1}$ is the direct sum of   
$\beta_{\nu}\cdot St_{G_1}$ ($\nu\in\Delta_t$), where we keep $\beta_{\nu}$ for the Brauer lift
of the Brauer character $\beta_{\nu}$. \itf  $st_{G_1}^+$ and $\sum_{\nu\in\Delta_t}\phi_\nu\otimes st_{G_1}$
have the same Brauer characters. It is well known that projective modules with 
 the same Brauer characters are isomorphic. So the result follows. 


\section{ $s$-components of the Steinberg-plus character}

In this section we recall some general facts of character theory
of finite groups of Lie type and results from the paper \cite{HZ2}
which contains some approach to the analysis of
characters of $G$ vanishing on the non-semisimple  elements.

Our main references for the character theory of groups of Lie type
are  \cite{C} and \cite{DM}. One of the principal notion of
the theory is that of the dual groups $\CG^*$ of $\GC$ and $G^*$ 
of $G$, see \cite[Ch. 4]{C} and \cite[13.10]{DM}. (Recall  that
to simplify  notation, 
we  keep $Fr$
for $Fr^*$ (the Frobenius endomorphism of $\GC^*$ defining $G^*$).) 
The duality
yields a bijection between the maximal tori $T$ of $G$ and maximal tori
$T^*$ of $G^*$ such that $T^*$ is naturally identified with $\Irr
T$, the set of irreducible characters of $T$. The group $G$ acts
on the set of maximal $Fr$-stable tori of $\CG$ by conjugation, and
this induces the action of $G$ on the set of pairs $({\mathbf T},
\theta)$, where $T={\mathbf T}^{Fr}$ and $\theta\in\Irr T\cong T^* $.
By \cite[13.13]{DM}, $G$-orbits of the pairs $({\mathbf
T},\theta)$ are in bijection with the $G^*$-orbits of the pairs
$(s,{\mathbf T}^*)$, where $s\in {\mathbf T}^{*Fr}$.

To every pair $({\mathbf T},\theta)$ the theory corresponds a
generalized character $R_{{\mathbf T},\theta}$ of $G$, called a
Deligne-Lusztig character of $G$. If $({\mathbf T}',\theta')$ is
another pair then $R_{{\mathbf T},\theta}=R_{{\mathbf
T}',\theta'}$ if $({\mathbf T},\theta)$ and $({\mathbf
T}',\theta')$ are $G$-conjugate, otherwise $(R_{{\mathbf
T},\theta},R_{{\mathbf T}',\theta'})=0$, where $(\cdot, \cdot)$
means the usual inner product of functions on $G$
\cite[11.15]{DM}.  The duality allows us to parametrize
$R_{{\mathbf T},\theta}$ by the $G^*$-orbits of the pairs $(s\in {\mathbf T}^{*})$.  We
 denote by ${\mathbf R }_s$ the set of pairs $({\mathbf
T},\theta)$,
where $s\in G^* $ is a fixed semisimple element,  whereas
$T^*$ vary within $C_{G^*}(s)$ such that $C_{G^*}(s)$-conjugate tori are counted once. 
In other words,  ${\mathbf
R}_s $ consists of the pairs $({\mathbf T},\theta)$ such that the dual $T^*$ of $T$
contains $s$ and $\theta$ corresponds to $s$ under the isomorphism $\Irr T\ra T^*$.
In addition, ${\mathbf R }_s$ contains at most one representative of the $G$-orbit of every pairs $(T,\theta)$. 
The \ir constituents of the
Deligne-Lusztig characters that belong to ${\mathbf R }_s$ form 
the Lustzig rational series usually denoted by ${\mathcal E}_s$.
This yields a partition $\Irr G=\cup_{s}{\mathcal E}_s$, when $s$
runs over representatives of the semisimple conjugacy classes in
$G^*$ \cite[14.41]{DM}.

\medskip
Recall that $W(T):=N_{G}(\TC)/T  $. The group $N_G(\TC)$  acts on $\TC$ by conjugation and stabilizes $T=G\cap\TC$.
This yields an action of $W(T)$ on $\Irr T$. If  $\theta$ is an irreducible
character of $T$ we set 
$W(T)_\theta=C_{W(T)}(\theta )$, that is, $W(T)_\theta$ is the stabilizer of $\theta$ in $W(T)$. 
In addition, recall that $\ep_{\GC}:=(-1)^{r}$, where $r$ is
the relative rank of a connected algebraic group $\GC$, which
therefore is meaningful for $\TC$ as well (consult \cite[pp. 64,66]{DM}).

\begin{propo}\label{or1}  Let $\Phi$  be a class function on $G$ vanishing at all non-semisimple elements. 
Let $\phi$ be the Curtis dual of $\Phi$. Then
$\Phi=\sum _{({\mathbf T},\theta)} 
\frac{(\phi|_T,\theta)}{|W(T)_{\theta}|}\ep_{{\mathbf
G}}\ep_{{\mathbf T}}R_{{\mathbf T},\theta}$, 
where the sum is over representatives of the $G$-orbits
of $(T,\theta)$. \end{propo}
Proof.   This  is \cite[Lemma 2.1]{HZ2}.

\med 

Applying Proposition \ref{or1} to the Steinberg-plus
character $St_G^+$ of $G=SO(V)$ and the Curtis dual $\om_G$ of $St_G^+$, we
have

\begin{equation}\label{eq12}
St_G^+=\sum_{(T,\theta)}\frac{(\om_G|_T,\theta)}{|W(T)_{\theta}|}\ep_{{\mathbf
G}}\ep_{{\mathbf T}}R_{{\mathbf T},\theta},\end{equation} \noindent where the
sum is over representatives of the $G$-orbits of $(T,\theta)$. For
every fixed semisimple element $s\in G^*$ we denote by $St^+_s$
the partial sum consisting from the terms with $({\mathbf
T},\theta)\in {\mathbf R}_s$. Thus, $St_G^+=\sum_{s}St^+_s$, with
the sum over representatives of the $G^*$-conjugacy classes of
semisimple elements of $G^*$. By the above comments, if $t$ is a
representative of another semisimple conjugate class of $G^*$ then
$St^+_s$ and $St^+_t$ has no common \ir constituents. 
It follows that $St^+$ is \mult free \ii so are $St^+_s$ for all
semisimple elements $s\in {\mathbf G}^*$. 
The argument depends on certain properties of
$s$. For convenience of further references, we record:

\begin{equation}\label{eq11}
St ^+_s=\sum_{(T,\theta)\in {\mathbf R}_s} \frac{(\om_G|_T,\theta)}{|W(T)_{\theta}|}\ep_{{\mathbf
G}}\ep_{{\mathbf T}} R_{{\mathbf T},\theta}.\end{equation} Obviously, if 
$St_s^+\neq 0$ then  $St_s^+$ is the sum of characters that belong to ${\mathcal E}_s$. Our goal is to determine these characters. However, we first
identify the cases where  $St^+_s=0$. 

\med
From now  on we specify $G$ to be $SO(V)$. 
Let $V^*$ be the natural module for $G^*$, $T^*$ is a maximal torus of $G^*$ and let $
V^*_1\oplus \cdots \oplus
V^*_{k+l}$ be the respective $T^*$-decomposition of $V^*$. If $s\in G^*$ is semisimple, we denote by $V^*_s$
the $1$-eigenspace $V_s^*$ of $s$ on $V^*$. If $s\in T^*$ then we denote by $T^*_s$ the projection of $T^*$ to 
$V^*$. So $T^*_s$ is a maximal torus in $SO(V^*_s)$.
Recall that, if $\dim V=2n+1$ is odd then 
$G^*=Sp_{2n}(q)$, otherwise $G^*\cong G.$ \itf $\dim V_s^*$ is even. (Indeed, if $G^*=Sp_{2n}(q)$ then 
$V_s^*$ is non-degenerate, and hence of even dimension. Let $G^*=SO(V^*)$ so $\dim V^*$ is even.   Viewing $s$ as an element of $\overline{V}^*$, one can choose a basis $b_1\ld b_{2n}$ of $\overline{V}^*$ such that $sb_i=\nu_ib_i$ for $i=1\ld 2n$ and some $\nu_i\in \overline{F}_q$. If $\nu_i\neq \pm 1$ then $\nu_i\up$ is an \ei of $s$, and one easily observes that the multiplicities of $\nu_i,\nu_i\up$ as \eis of $s$ are equal. As $\det s=1$, it follows that the \mult of $-1$ is even (or equals zero). Therefore, the \ei 1 must be of even multiplicity, as claimed.)

\begin{propo}\label{pp5}  
 $St^+_s=0$ \ii  $V^*_s$ has Witt defect $1.$ In particular, if $\dim V$ is odd  then 
  $St^+_s\neq 0$ for every semisimple element $s\in  G^*$. 
\end{propo}
Proof. As distinct characters $R_{{\mathbf T},\theta}$ with $({\mathbf T},\theta)\in {\mathbf R}_s$ 
are linearly independent, $St^+_s= 0$ if and only if every term in (\ref{eq11}) equals 0, that is,
$(\om_G |_T,\theta)=0$ for every $({\mathbf T},\theta)\in {\mathbf R}_s$. By Corollary \ref{nc4}, 
$(\om_G |_T,\theta)=0$ \ii   $V_s^*$  is contained in $V^*_1\oplus \cdots \oplus
V^*_{k}$, see (\ref{eq2}). This is equivalent to saying that  $T^*_s$ is a neutral torus in $SO(V^*_s)$  (Definition \ref{de4}). 
As $\dim V^*$ is even, $SO(V^*_s)$ has no neutral torus \ii    $V_s^* $ is of Witt defect $1$. 


\med Remark. The second statement of Proposition \ref{pp5} can be proven in a more conceptual way. This follows from the fact that  $St^+-\Gamma_G$ is a proper character of $G$, where $\Gamma_G$ denotes the Gelfand-Graev character of $G=SO(V)$. Indeed,   $\Gamma_G$ is \mult free and consists of $q^n$ irreducible constituents called the regular characters, see \cite[14.42,14.39]{DM}. In addition, as  $\Gamma_G$ vanishes at the non-identity semisimple elements of $G$ and $St_G^+$ vanishes at the non-semisimple ones, it follows that $(St_G^+, \Gamma_G)=St_G^+(1)/St_G(1)=q^n$. So every \ir constituents of  $\Gamma_G$ must occur in $St_G^+$.

\med
The \f result from \cite[Theorem 1.3]{HZ2} plays an essential role here:

\begin{theo}\label{th4}    Let $s\in G^*$   be semisimple and
let $\phi$ be a generalized character of     $G$.
 Suppose that $(\phi |_T,\theta)\in\{0,1\}$
for all $(\TC , \theta)\in {\mathbf R}_s$ corresponding to the conjugacy class of
$s$. Then $\phi \cdot St$ has at most one constituent from
${\mathcal E}_s$, and this is
 a unique  regular character.     \end{theo}

We shall apply this theorem to the situation where
$\phi=\om_G$ and $\phi\cdot St_G=St_G^+$.

\begin{propo}\label{pp4}  If s does not have \ei $1$ on $V^*$
then $St^+_s$ is a regular character in ${\mathcal E}_s$.
\end{propo}

Proof. By Lemma \ref{d2g}, if $(\TC , \theta)\in {\mathbf R}_s$ and $s$ does not have \ei 1 on $V^*$, then  $(\om_G|_T, \theta)=1$. Then the result follows from Theorem \ref{th4}.

\med
It is convenient to state here the \f technical fact:

\bl{tf5}  Let $V^*_s$ be the $1$-eigenspace of s on $V^*$. Suppose that  $V^*_s\neq 0$ and $St^+_s\neq 0$.
Then  $V^*_s$ has Witt defect $0$. 
\el  Let $T^*$ be a maximal torus containing $s$.  Obviously, $T^*V^*_s=V^*_s$. Therefore, a $T^*$-decomposition of $V_s^*$ can be extended to that of $V^*$. 
In other words, there is a $T^*$-decomposition $V_1^*\oplus\cdots \oplus V_{k+l}^*$ of $V^*$ such that   $V_s^*$
is the sum of some $V^*_i$ $(1\leq i\leq k+l)$. If $i>k$ for some $V_i^*\subseteq V^*_s$ then $St^+_s= 0$ by Proposition \ref{pp5}. This contradicts the assumption. Therefore, $V_s^*$ is contained in $V_1^*\oplus\cdots \oplus V_{k}^*$, and hence $V_s^*$ is  the direct sum of non-degenerate subspaces of Witt defect 0. So  $V_s^*$ is of Witt defect 0.

\section{The constituents of the Steinberg-plus character}

Recall that  $G=SO(V)$ if $q$ is odd, otherwise $G=\Omega(V)$, and $G^*$ the dual group.
As above, $V^*$ denote the natural module for $G^*$. Observe  that $G^*=Sp(V^*)$ if $\dim V$ is odd, 
 otherwise $G^*=SO(V^*)$. For a semisimple element $s\in G^*$ we keep
notation $V^*_s$ for  the 1-eigenspace  of $s$ on
$V^*$. 
 Then $V^*_s$  is non-degenerate of even dimension. In this section we deal with the case, where $V^*_s\neq 0$ and  $St^+_s\neq 0$. By Lemma \ref{pp4}, $V^*_s$ contains a totally singular (or totally isotropic) subspace $U$ of dimension $(\dim V^*_s)/2$.
Let $P_U$ be the stabiliser of $U$ in $G^*$, and $L^*$ a Levi subgrouop of $P_U$. The group $L^*$ and its dual group $L\subset G$ plays an important part  in the reasonong below. Observe that there is a bijection between parabolic subgroups of $G$ and $G^*$, and that the one corresponding to $P_U$ can also be defined as the stabilizer of a totally singular subspace$U'$, say, of $V$.
Similarly, $L$ is the stabilizer of the direct sum of two totally singular subspaces of $V$, each of dimension equal to $\dim U$. 
Let $V' $   be a complement of $U'$ in $U^{\prime \perp}$.  Then $L\cong GL(U') \times SO(V')$ if $q$ is odd, and $L\cong GL(U) \times \Omega (V')$ if $q$ is even. If $\dim V$ is odd then 
$L\cong GL(U') \times Sp(V')$. This observation is frequently  used below without reference. Furthermore,
we usually write $L=G_1\times G_2$, where  $G_1\cong GL(U')$ and $G_2\cong SO(V')$.

Note that  $V^*=V^*_s\oplus V^{*\perp }_s$.
We write  $s=\diag(s_1,s_2)$, where $s_1$, respectively, $s_2$, is the restriction of  $s$ to $V^*_s$, 
respectively,  $V^{*\perp }_s$.

We define  a subgroup   $X^*=X^*_1\times X^*_2$ as follows. If $G^*=Sp(V^*)$   then we set  $X_1^*=Sp(V_s^*)$, $X^*_2= Sp(V_s^{*\perp})$, otherwise, 
  $X_1^*=SO(V^*_s)$, $X_2^*=  SO(V_s^{*\perp})$ for $q$  odd, and $X^*$
is $X_1^*=\Omega(V^*_s)$, $X_2^*=\Omega(V_s^{*\perp})$ for $q$ even. 

Observe that we also have $\overline{V}^{*}=\overline{V}_s^{*}\oplus \overline{V}_s^{*\perp}$. Then the stabilizer of $\overline{V}_s^{*}$ in $O(\overline{V}^{*})$ is $O(\overline{V}_s^{*})\times O(\overline{V}_s^{*})$ if $O(\overline{V}^{*})$ is orthogonal. The connected component of this algebraic group is ${\mathbf X}^*=SO(\overline{V}_s^{*})\times SO(\overline{V}_s^{*})$ if $q$ is odd, and hence $X^*={\mathbf X}^{*Fr}$. Similarly, in the other cases.

As above, $X^*$ is the dual group of $X\cong SO(U')\times SO(V')$. It is clear that $G_2=X_2$. Similarly,  $L^*=G^*_1\times G^*_2$, where $ G^*_2 =X^*_2$.

\begin{lemma}\label{x9}  Let $s\in G^*$ be a semisimple element such that $V^*_s\neq0$. 
Then $C_{G^*}(s)\subseteq X^*$ \ii either $\dim V$ is odd or $-1$ is not an \ei of s. In particular, if $q$ is even then
$C_{G^*}(s)\subseteq X^*$ for any s.
\end{lemma}
Proof. Set $Y=C_{G^*}(s)$.  Then $YV^*_s=V^*_s$ and $YV^{*\perp}_s=V^{*\perp}_s$. If  $\dim V$ is odd then  
 $G^*=Sp(V^*)$ and $X^*= Sp(V^{*}_s)\times Sp(V^{*\perp}_s)$. So the lemma follows in this case.

Suppose that $\dim V$ is even. Let    $s=\diag(s_1,s_2)$ as above.      Then
 $Y\subset O(V^{*}_s)\times C_{ O(V^{*\perp}_s)}(s_2)$. As 1 is not an \ei of $s_2$,
by Lemma \ref{y9}, $C_{O(V_s^{*\perp})}(s_2)\subset 
SO(  V_s^{*\perp})$ \ii   $-1$ is not an \ei of $s_2$.
Let $y\in Y$, and $y=\diag(y_1,y_2)$
for $y_1\in O(V^{*}_s)$, $y_2\in  O(V^{*\perp}_s)$. 
If $q$ is even then $y_2\in X_2^*$, and hence  $y_1\in X_1^*$ as the spinor norm of $y$ is the product of the spinor norms of $y_1$ and $y_2$. So $y\in X^*$ is this case. 

Suppose first that $q$ is odd. Then $y\in SO(V^*)$ \ii  $\det y_1=\det y_2$.
Therefore, $y\notin X^*$ \ii $\det y_1=\det y_2=-1$. 
 As $1$ is not an \ei of $s_2$, the lemma follows in this case too.

\med

Next suppose that $s$ has \ei 1 on $V^*$ and $St^+_s\neq 0$.  Let $U$ be a maximal totally singlar (or totally isotropic)   subspace of $V^*_s$. Let $P_U$ be the stabilizer of $U$ in $G^*$.
If $u:=\dim U<\dim V/2$ then $P_U$ is determined by $u$ up to $G^*$-conjugacy.

\begin{lemma}\label{no9}  Let $s\in G^*$ be a semisimple element such that $V^*_s\neq0$. 
Suppose that $St^+_s\neq 0$.

$(1)$
 Let $U_1,U_2$ be maximal totally   singular  subspaces of $V^*_s$. Then  $gU_1=U_2$ for 
some $g\in C_{G^*}(s)$, unless $\dim V$ is even and $-1$ is not an \ei of s. 
In the exceptional case there are two $C_{G^*}(s)$-orbits of 
maximal totally singular subspaces of $V^*_s$. 

$(2)$ Let $P_i$ $(i=1,2)$ be the stabilizer of 
$U_i$ in $G^*$. Then there is a Levi subgroup of $P_i$ $(i=1,2)$ stabilizing $V_s^*.$

$(3)$ Let  $L^*_i$ be a Levi subgroup of $P_i$ stabilizing $V_s^*$. Then $gL^*_1=L^*_2$ for 
some $g\in C_{G^*}(s)$, 
unless $\dim V$ and $\dim U_1$ are even and $s$ does not have \ei $-1$. In the exceptional case 
$L^*_1$, $L^*_2$ are $C_{G^*}(s)$-conjugate \ii $gU_1=U_2$ for some $g\in C_{G^*}(s)$.

 \el 
 By Lemma $ \ref{tf5}$, $V^*_s$ is of Witt defect $0.$ 

(1) Let ${\mathcal U}$ be the set of all maximal totally singular  subspaces of $V^*_s$. If $C_{G^*}(s)$ is not contained in $X^*$ then the restriction of $C_{G^*}(s)$ to $V^*_s$ coincides with $O(V_s^*)$. As $O(V^*_s)$ is transitive on 
${\mathcal U}$, in this case $C_{G^*}(s)$ is transitive on 
${\mathcal U}$, and the claim follows. Suppose that $C_{G^*}(s)\subseteq X^*$. By Lemma \ref{x9}, this happens \ii
either $\dim V$ is odd or  $-1$ is not an \ei of $ s$. In the latter case the restriction of $C_{G^*}(s)$ to $V^*_s$ coincides with $SO(V_s^*)$. As $\dim V^*_s$ is even,    the claim now follows  from Lemma \ref{no6}. 

Suppose that $\dim V$ is even  and  $-1$ is an \ei of $s$. Then, by Lemma \ref{x9},  $C_{G^*}(s)\not\subseteq X^*$, and hence 
 the restriction of $C_{G^*}(s)$ to $V_s^*$ coincides with $O(V^*_s)$. So   
$C_{G^*}(s)$ is transitive on 
${\mathcal U}$.

(2) Let $V^*_s=U_i\oplus U_i'$, where $U'_i$ is an arbitrary totally  singular subsgace of  $V^*_s$.
Then the stabilizer $L_i^*$ of $U_i$ and $U'_i$ in $G^*$ is known to be a Levi subgroup of $P_i$. Obviously,
$L_i^*$ stabilizes $V_s^*$,
whence (3).   (Note that all Levi subgroups of $P_i$ stabilizing $V_s^*$
can be obtained in this way.)

(3)    Let $P_i'$ denote the stabilizer of $U_i$ in $X_1^*$, 
and let $L_i^{*\prime}$ be the projection of $L^*_i$ to $X^*_1$ (so $L_i^{*\prime}\cong G_1$). Then $P_i'$ is a parabolic subgroup of $X^*_1$ and $L_i^{*\prime}$ is a Levi subgroup of $P_i'$. 

Suppose first that $\dim V$ is odd. Then $C_{G^*}(s)$ is  contained in $X^*$ by Lemma \ref{x9}. Note that here $X^*=X^*_1\times X_2^*$, where $X_1^*=Sp(V^*_s)$ and $G_2^*=Sp(V_s^{*\perp})$. In addition, 
$L_i^*=G_1^*\times G_2^*$, where $G_1^*\cong GL(U)\subset X_1^*=Sp(V^*_s)$, $G_2^*=X^*_2=Sp(V_s^{*\perp})$. Therefore, $g'U_1=U_2$ for some $g'\in C_{G^*}(s)$. So we can assume that $P_1=P_2$. Then  $L_1^*,L^*_2$ are $C_{G^*}(s)$-conjugate \ii $L_1^{*\prime}$, $L_2^{*\prime}$ are $X^*_1$-conjugate. This is the case as Levi subgroups of any parabolic subgroup are conjugate. 

Suppose that $\dim V$ is even. Recall that $L_i^*=G_1^*\times G_2^*$, where $G_1^*\cong GL(U)\subset X_1^*=SO(V^*_s)$, $G_2^*=SO(V_s^{*\perp})$ (respectively, $X_1^*=\Omega (V^*_s)$, $G_2^*=\Omega(V_s^{*\perp})$). So
$G_2^*=X_2^*$ does not depend on the choice of $L^*_i$. Observe  that $C_{G^*}(s)$  coincides either with $X_1\cdot C_{X_2^*}(s_2)$ or, for $q$ odd,  with $\lan g_1g_2\ran ( X_1\cdot C_{X_2^*}(s_2))$, where $g_1\in O(V_s^*)$, $g_2\in C_{O(V_s^{*\perp})}(s_2)$ and $\det g_1=\det g_2=-1$.  So in the latter case the restriction of   $C_{G^*}(s)$ to $V^*_s$ coincides with $O(V_s^*)$. Clearly, this happens \ii $C_{G^*}(s)$ is not contained in $X^*$.

 As  $X_2^*$ is normal in $O(V_s^{*\perp})$, it follows that $L^*_1,L^*_2$ are $C_{G^*}(s)$-conjugate \ii 
$L_2^{*\prime}=hL_1^{*\prime}$ for some $h$ 
from the projection of $C_{G^*}(s)$ to $V^*_s$. As Levi subgroups of any parabolic subgroup are conjugate, it follows that 
$hL_1^{*\prime}=L_2^{*\prime}$ for some $h\in X_1^*$, unless $U_1,U_2$ are not in the same $X_1^*$-orbit.
The latter happens \ii $\dim U_1$ is even (Lemma \ref{no6}). Note that $g'U_1=U_2$ for some $g'\in O(V^*_s)$. It follows that, if  $C_{G^*}(s)\not\subseteq X^*$ then $g'U_1=U_2$ for some $g'\in C_{G^*}(s)$. So in this case (3) follows. 
Thus, we are left with the case where $C_{G^*}(s)\subseteq X^*$ and $\dim U_1$ is even. In this case (3) follows by Lemma \ref{x9}.

\med The \f lemma is one of the key points of our argument.

\begin{lemma}\label{no8}  Let $s\in G^*$ be a semisimple element, and let 
$T^*$ be a maximal  torus in $G^*$ containing $s.$ Suppose that $V^*_s\neq 0$. Then the \f are equivalent:

$(1)$ 
$(\om_G|_T,\theta)\neq 0;$ 

$(2)$ $T^*$ is $C_{G^*}(s)$-conjugate to a subgroup of $P_U$, where U is some  totally singular  subspace of $V^ *_s$  of dimension  $\dim V^*_s/2$. 

$(3)$ $T^*$ is $C_{G^*}(s)$-conjugate to a Levi subgroup $L^*$ of $P_U$ such that 
$L^*V^*_s=V^*_s$.
 
In addition,  suppose that $-1$ is not an \ei of $s$. Let $T^*$, $T^{\prime *}$ be two maximal tori of $L^*$. If $T^*$, $T^{\prime *}$ are $C_{G^*}(s)$-conjugate then they are $C_{L^*}(s)$-conjugate. 
\el
$(1)\ra (3)$ By Lemma \ref{no9} (see the proof of item (1)),  $T^*V^*_s=V^*_s$ and any $T^*$-decomposition of $V_s^*$ can be extended to  a  $T^*$-decomposition $V_1^*\oplus\cdots \oplus V_{k+l}^*$  of $V^*$. Moreover,  $V^*_s$ is contained in $V_1^*\oplus\cdots \oplus V_{k}^*$ (Corollary  \ref{hh2}).
Each $V^*_i$ with $i\leq k$ is the direct sum of two totally singular  $T^*$-stable subspaces of equal dimension. Let $U$ be the sum of any such subspaces, chosen by one from every $V_i^*\subset V^*_s$. Then $U$   is totally  singular,  $T^*$-stable and  of dimension $\dim V^*_s/2$. Moreover, it is easy to observe that there is a $T^*$-stable  totally singular (or totallyisotropic) subspace $U'\subset V_s^*$ such that $ V_s^*=U\oplus U'$.  Then the stabilizer of both $U,U'$ in $G$
is a Levi subgroup of $P_U$, as required. 

 $(2)\ra (1)$ follows from Lemma \ref{d2g}, and $(3)\ra (2)$ is trivial.

For the additional statement, let $L^*=G^*_1\times G^*_2$, where $G_1^*$ is the projection of $L^*$ to $V_s^*$ and $G_2^*$ is the projection of $L^*$ to $V_s^{*\perp}$. Then let  $T^*=T_1^*\times T_2^*$ with $T_1^*\subset G_1^*$, $ T_2^*\subset G_2^*$, and similarly, let  $T^{\prime *}=T_1^{\prime *}\times T_2^{\prime *}$.  As $-1$ is not an \ei of $s$, we have $C_{G^*}(s)\subseteq X^*=X_1^*\times X_2^*$ (Lemma \ref{x9}). Since $X_2^*=G_2^*$, we have 
$C_{L^*}(s)=G_1^*\times C_{G_2^*}(s_2)$ and $C_{G^*}(s)=X_1^*\times C_{G_2^*}(s_2)$. \itf 
 $T^*,T^{\prime *}$ are $C_{L^*}(s)$-conjugate \ii $T_1^{*},T_1^{\prime *}$ are $G_1^*$-conjugate. As  $T^*,T^{\prime *}$ are  $C_{G^*}(s)$-conjugate, it follows that  $T_1^{*},T_1^{\prime *}$ has the same partition function as maximal tori of $X^*_1$, and hence as those in $G^*_1\cong GL(U)$. Therefore, $T_1^{*},T_1^{\prime *}$  are conjugate in $G_1^*$, as claimed.

\med
Remark. Strictly speaking, we have to prove that if ${\mathbf T}^*,{\mathbf T}^{\prime *}$ are $Fr$-stable tori in the algebraic groups $C_{{\mathbf L}^*}(s)$ such that $T^*={\mathbf T}^{*Fr}$ and  $T^{\prime *}=({\mathbf T}^{\prime *})^{Fr}$, and if  ${\mathbf T}^*$ 
is $C_{G^*}(s)$-conjugate to ${\mathbf T}^{\prime *}$ then ${\mathbf T}^*$  
is $C_{L^*}(s)$-conjugate to ${\mathbf T}^{\prime *}$. Clearly, this only requires  routine changes of the above reasoning. 


\med
	Let ${\mathcal L}$ be the set of all Levi subgroups $L^*$ of $P_U$, when $U$ runs over all maximal totally  singular subspaces of $V_s^*$, satisfying the condition $L^*V_s^*=V_s^*$. By Lemma \ref{no9}, 
${\mathcal L}$  forms a single $C_{G^*}(s)$-orbit, unless $\dim V$ and $\dim U$ are even and $-1$ is not an \ei of $s$.
In the latter case ${\mathcal L}$  consists of two $C_{G^*}(s)$-orbits. With this notation we have the \f refinement of Lemma \ref{no8}: 

\begin{lemma}\label{aa2b} 
Suppose $G^*=SO^\pm(2n, q)$, $n$ even and let $s=\Id$. 
 Let $L^*_1,L^*_2\in{\mathcal L}$ be two non-conjugate subgroups. Let $T^*_i$  be an maximal torus of $L^*_i$ $(i=1,2)$ corresponding to the same  function $i\ra d_i$. 
Then $T^*_1$ is not $C_{G^*}(s)$-conjugate to $T^*_2$ \ii   $T^*_i$ are exceptional tori in $ G^*$.
\end{lemma}

Proof. 
By Lemma \ref{tt5}, non-exceptional neutral maximal tori of $G^*$ with the same function $i\ra d_i$ are conjugate in $G^*$
as stated,
whereas  exceptional neutral maximal tori of $G^*$ with the same function $i\ra d_i$ form two $G^*$-conjugacy classes. 
 Suppose the contrary, that $T^*_1$ and $ T^*_2$ are exceptional but $G^*$-conjugate.  As every neutral maximal torus of $G^*$
is conjugate to that in $L^*_1$ or $L^*_2$, and maximal tori in each $L^*_1,L^*_2$  with the same function $i\ra d_i$ are conjugate, it follows that all neutral maximal  tori in $G^*$   with the same function $i\ra d_i$ are conjugate in $G^*$, which is a contradiction.

\begin{lemma}\label{2es}  Suppose that  $(\om_G|_T,\theta)\neq0$ and, for $s\in T^*$ corresponding to $\theta$, let  $V^*_s\neq 0$.

$(1)$ Suppose that ${\mathcal L}$ consists of a single 
$C_{G^*}(s)$-orbit, and
let $L^*\in {\mathcal L}$. 
Then 
maximal tori of $C_{G^*}(s) $ are $C_{G^*}(s)$-conjugate to tori in  $L^*$. 
 
$(2)$  Suppose that ${\mathcal L}$  consists of two $C_{G^*}(s)$-orbits, and let
$L^*_1,L^*_2\in {\mathcal L}$  be their representatives.
Then every 
maximal torus  $T^*$  of $C_{G^*}(s)$ is $C_{G^*}(s)$-conjugate to a torus of $L_1^*$ or $L^*_2$.
 Moreover, $T^*$ is $C_{G^*}(s)$-conjugate to a torus of each $L^*_1$ and $L^*_2$ \ii $T_1^*$ is non-exceptional.

\end{lemma}

Proof. Let $T^*$ be maximal torus of $X^*\subset C_{G^*}(s)$. By Lemma \ref{no8},  $T^*$  is $C_{G^*}(s)$-conjugate to  a torus in 
a Levi subgroup $L^*\in {\mathcal L}$. So (1) is immediate. Consider (2). As  in the case (1), it follows that 
every maximal torus of $X^*$ such that $(\om_G|_T,\theta)\neq0$ is $X^*$-conjugate to a torus in $L^*_1$ or $ L^*_2$, whence the first claim of (2). 

Furthermore, $T^*\subset X^*=X_1^*X_2^*$. Let $T^*=T_1^*T_2^*$, where $T^*_1\subset X_1^*, T^*_2\subset X^*_2$. 
By Lemma \ref{no9}, in this case $\dim V$ and $\dim V_s^*/2$ are even and $-1$ is not an \ei of $s$.
So $C_{G^*}(s)\subseteq X^*$ by Lemma \ref{x9}. Therefore, maximal torus $T^*,T^{\prime *}$ of $X^*$ are $C_{G^*}(s)$-conjugate \ii they are
$X^*$-conjugate, and hence \ii  $T^*_1,T^{\prime *}_1$ are $X_1^*$-conjugate. If $T_1^*$ is non-exceptional in 
$X^*_1$, then $T^*_1,T^{\prime *}_1$ are conjugate in $X^*_1$ (Lemma \ref{aa2b}), and the statement is true in this case. 

Suppose that $T^*_1$ is exceptional. Denote by $L^{\prime *}_1,L_2^{\prime *}$ the projections of $L^{ *}_1, L^{ *}_2$ to $X^{*}_1$. Then $L^{ *}_1, L^{ *}_2$ are $C_{G^*}(s)$-conjugate \ii $L^{\prime *}_1,L_2^{\prime *}$ are $X^*_1$-conjugate. Therefore, the second statement in (2) follows from Lemma \ref{aa2b} applied to $X_1^*\cong SO(V_s^*)$.

\bl{yr1}  Let $G=SO(V)$, $L^*\in {\mathcal L}$ and  let $s\in T^*\subset L^*$. 
 Suppose that $\dim V$ is odd, or  $\dim V$ is even and $-1$ is not an \ei of $s$. 
 Then $W(T)_{\theta}=W_X(T)_{\theta}.$
 
\el  Recall that $G^*=Sp(V^*)$ if $\dim V$ is odd.  Note that $C_{G^*}(s)$ is contained in $X^*_1\cdot  X^*_2$, where $X^*_1=Sp(V_s^*)$ and $X^*_2=Sp(V_s^{*\perp})$.  (Indeed, $X^*_1\cdot  X^*_2$ coincides with the stabilizer of $V^*_s$ in $G^*$.) 
This implies the  
claim, as $W(T)_{\theta}$ is canonically isomorphic to $W_{C_{G^*}(s)}({T^*})=C_{W(T^*)}(s)$, 
which is exactly 
$N_{C_{{\mathbf G}^*}(s)}({\mathbf T}^*)/T^*$. Next, suppose $\dim V$ is even.  Observe that $W(T)_\theta\cong C_{W(T^*)}(s)=N_{C_{G^*}(s)}(\mathbf T)/T$ and $W_X(T)_\theta\cong C_{W_{X^*}(T^*)}(s)=N_{C_{X^*}(s)}(\mathbf T)/T$. As $C_{G^*}(s)\subset X^*$ by Lemma \ref{x9}, we have $C_{G^*}(s)=C_{X^*}(s)$, and the lemma follows.

\med

Recall that we write $L^*=G^*_1\times G^*_2$, where $G_1^*\cong GL(U)$ and $G^*_2\cong X_2^*$.)   If $T^*$ is a maximal torus of $L^*$, we write $T^*=T^*_1\times T^*_2$, where $T^*_r$ is a maximal torus in $G^*_r$ for $r=1,2$.  Similarly, 
If $T$ is a dual torus in $L$, we write $T=T_1\times T_2$, where $T_r$ is a maximal torus in $G_r$
 

We denote by ${\mathbf R}^L_s$ the set of representatives of the $L$-conjugacy classes of pairs $({\mathbf T},\theta)$ with ${\mathbf T}\subset {\mathbf L}$. They are in bijection with   $C_{L^*}(s)$-conjugacy classes of maximal tori $T^*\subset C_{L^*}(s)$.

\begin{lemma}\label{pr9}  Let $G=SO(V)$ with $\dim V$ odd, and let $s\in G^  *\cong Sp(V^*)$ be a semisimple element. Let $V_s^*$ be the $1$-eigenspace of s on  $V^*$. Suppose that $V_s^*\neq 0$,  and let U be a maximal  totally isotropic  subspace of 
$V_s^*$, $P_U$ the stabilizer of U in G,  and let $L^*$ be a Levi subgroup of   $P_U$. Then

\begin{equation}\label{eq17}St^+_s=
\sum_{({\mathbf T},\theta)\in {\mathbf R}_s}\frac{(\om_G|_T,\theta)}{|W(T)_\theta|}\ep_{{\mathbf
G}}\ep_{{\mathbf T}} R_{{\mathbf T},\theta}=
 \sum_{({\mathbf T},\theta)\in {\mathbf R}^{L}_s}\frac{1}{|W_{L}(T)_\theta|}\ep_{{\mathbf
G}}\ep_{{\mathbf T}} R_{{\mathbf T},\theta}.\end{equation}\end{lemma}
Proof.  Recall that  $(\om_G|_T,\theta)=0$, 
if $T$ is not conjugate to a torus in $L$ (Lemma \ref{no8}). So, by Lemma \ref{yr1}, the elements of ${\mathbf R}_s$ corresponding to the non-zero terms of (\ref{eq17}) can be chosen in ${\mathbf R}^L_s$.  As above, let $s=\diag(s_1,s_2)$, where $s_1$ is the projection of $s$ to $V_s^*$ and $s_2$  is the projection of $s$ to $V_s^{*\perp}$.

  Let $x\in X$, and let $x=x_1x_2\in X$ with $x_i\in X_i$, $i=1,2$.  By Lemmas \ref{mt6} and \ref{vg5},  $\om_G(x)=\om_{X_1}(x_1)\om_{X_2}(x_2)$, and hence $(\om_G|_T,\theta)=(\om_{X_1}|_{T_1}, \theta_1)\cdot (\om_{X_2}|_{T_2}, \theta _2)$, where
 $\theta_r\in \Irr T_r$ corresponds to $s_r  $, $r=1,2$. As $s_1=\Id$, we have $\theta_1=1_{T_1}$. 
As $s_2$ does not have eigenvalue 1 on $V_s^{*\perp}$, by Lemma \ref{d2g}(3), we have  $(\om_{X_2}|_{T_2},\theta _2)=1$. Furthermore, the torus $T_1$ is neutral in $X_1$. By Lemma \ref{d2g}(2),  $(\om_{X_1}|_{T_1},\theta _1)=2^{m(T_1)}$,
where $m(T_1)$ is the number of the parts in the partition that determines $T_1$ (because $\theta_1=1_{T_1}$). Therefore,  $(\om_G|_T,\theta)=(\om_{X_1}|_{T_1}, \theta_1)=2^{m(T_1)}$. 

Let $T_1$ correspond to a function $i\ra d_i$ (with zero function $j\ra e_j$). Then  $\sum _id_i=m(T_1)$ and $|W_{X_1}(T_1)|=\Pi_{i}(2i)^{d_i}d_i !$ by Lemma \ref{a33}  whereas $|W_{G_1}(T_1)|=\Pi_{i}i^{d_i}d_i !$.
So $|W_{X_1}(T_1)|=2^{m(T_1)}\cdot |W_{G_1}(T_1)|$.   

Finally, by Lemma \ref{yr1}, $W(T)_{\theta}=W_X(T)_{\theta}.$ We recall that $W_X(T)_\theta\cong C_{W(T^*)}(s)=W_{X^*_1}(T_1^*)\times C_{W_{X^*_2} (T_2^*)}(s_2)$ as $s_1=1$. By the above $|W_{G_1}(T_1)| =2^{m(T_1)}\cdot |W_{G_1}(T_1)|$.   As $L^*=G^*_1\times G^*_2$ and $G^*_2=X^*_2$,
 the lemma follows.

\medskip
Let  $G=SO(V)$, where $\dim V$ is even. 
 In Lemma \ref{dx5} below, we assume that $-1$ is an \ei of $s$, and hence, by Lemma \ref{no9}, 
${\mathcal L}$ forms a single $C_{G^*}(s)$-orbit. Statement (A) of Lemma \ref{dx5} refines further Lemma \ref{2es}(1) to obtain a statement similar to that in Lemma \ref{no8} for the case where $-1$ is not an \ei of $s$. 


\bl{dx5}  Let $G=SO(V)$, where $\dim V$ is even, $L\in {\mathcal L}$, 
and  $s\in T^*\subset L^*$. Let ${\mathcal S}$ be the set of maximal  $Fr$-stable  tori ${\mathbf T}^{'*}\subset  C_{L^*}(s)$  that are $ C_{G^*}(s)$-conjugate to ${\mathbf T}^*$. Let $Y$ be the $-1$-eigenspace 
of $s$ on $V^*$, and $T^*_3$ the restriction of $T^*$ to Y. Suppose $Y\neq0$.  
 
\med
$(A)$ The set ${\mathcal S}$   
consists of two $  C_{L^*}(s)$-orbits \ii $T^*_3$ is exceptional  but $T^*_1$ is not exceptional.

\medskip
$(B)$ The \f statements are equivalent:
\med

$(1)$ $N_{C_{  G^*}(s)}({\mathbf T}^*)$ is  contained in $X^*$. 

$(2)$ 
either $T_1$  or $T_3$ is exceptional;

$(3)$ The subgroup $W_X(T)_{\theta}= W(T_1)\times W_{X_2}(T_2)_{\theta_2}$ coincides with  
$W(T)_{\theta}$.  

\med
$(C)$ The \f statements are equivalent:
\med

$(4)$ $N_{C_{  G^*}(s)}({\mathbf T}^*)$ is  not contained in $X^*$. 

$(5)$  $T_1$ and $T_3$ are non-exceptional;

$(6)$ The subgroup $W_X(T)_{\theta}= W(T_1)\times W_{X_2}(T_2)_{\theta_2}$  has index $2$  in
 $W(T)_{\theta}$.  
\el 
(A) Set $M:=(V_s^*+Y)^\perp$.
Let $s_2=\diag(-\Id,s_4)$, where $s_4$ is the restriction of $s$ to $M$. Then
 $C_{X^*}(s)=X^*_1\times X_3^*\times X_4^*$
where $X_3^*=SO(Y)$ and $X_4^*=C_{SO(M)}(s_4)$, whereas $C_{G^*}(s)=\lan C_{X^*}(s), g\ran$, where $g=\diag(g_1,g_2,\Id)$, $\det g_1=\det g_2=-1$ and $g_1,g_2$ are the projections of $g$ to $U,Y$, respectively (see Lemma \ref{y9}).
Also $C_{L^*}(s)=G_1\times X_3\times X_4^* $. Furthermore, $T^*=T^*_1 \times T^*_3 \times T^*_4$, where $T^*_3,T^*_4$ are projections of $T^* $ to $Y,M$, respectively. Similarly, let $T^{\prime*}=T_1^{\prime*}\times T_3^{\prime*}\times T_4^{\prime*}$. Note that $T_1^{*},T_1^{\prime*}$ are conjugate in $O(V_s^{*})$, so they correspond to the same partition function as tori in $X_1^{*}$, and hence in $G_1^{*}$. 
 We know that the maximal tori corresponding to the same partition are conjugate in $G^*_1\cong  GL(U)$. As  $X_4^*$ is common for $C_{L^*}(s)$ and $C_{G^*}(s)$, it suffices to look at $T^*_3$. Again, $T^*_3$ and $T_3^{\prime *}$ are conjugate in $X^*_3$ unless $T^*_3$ is exceptional, see Lemma \ref{aa2b}.  

So we are left with the case where $T^*_3$
is exceptional. 
As $T^*_1$ and $T_1^{\prime *}$ are conjugate in $G^*_1\cong GL(U)$,   we can assume $T^*_1=T_1^{\prime *}$. 
 If $T^*_1$ is non-exceptional then $N_{O(V^*_s)}({\mathbf T}^*_1)$ contains an element $g_1$ with $\det g_1=-1$ (Lemma \ref{a33}). Therefore, $C_{G^*}(s)$ contains an element $g=\diag(g_1,g_2,\Id)$ with $\det g_2=-1$. \itf the restriction of $C_{G^*}(s)$ to $Y$ coincides with $O(Y)$, and hence all maximal tori of $X_3^*$ whose partition function is the same as that of $T^*_3$ are in the same $O(Y)$-orbit  (Lemma \ref{my7}). This implies that, given a maximal torus $D$, say, in this orbit, there is a torus in  ${\mathcal S}$ whose restriction to  $Y$ is $D$. However, these tori $D$ form two $SO(Y)$-orbits (Lemma \ref{my7}). \itf that ${\mathcal S}$ consists of  two $C_{L^*}(s)$-orbits.  

Finally, suppose that  $T^*_1$ is exceptional. Let $h\in C_{G^*}(s)$ be an element such that $h{\mathbf T}^{ *}h\up={\mathbf T}^{\prime *}$. Then $h =\diag(h_1,h_3,h_4)$, where $h_1\in O(V^*_s)$, $h_3\in O(Y)$
and $h_4\in O(M)$. As above, we can assume that ${\mathbf T}^*_1={\mathbf T}_1^{\prime *}$ and ${\mathbf T}^*_4={\mathbf T}_4^{\prime *}$.
 Then $N_{O(V_s^*)}({\mathbf T}_1^*)\subset  SO(V_s^*)$ and $N_{O(M)}({\mathbf T}_1^*)\subset SO(M)$. Therefore, $\det h_1=1$ and $\det h_4=1$. As $h\in G^*=SO(V^*)$, it follows that $\det h_3=1$ and hence $h_3\in SO(Y)=X_3^*$. As $X_3^*$ is a multiple of  $L^*$, we can assume that ${\mathbf T}^*_3={\mathbf T}_3^{\prime *}$. So (A) 
follows.

(B) Observe first that $(1)$ and $ (2)$ are equivalent. Let ${\mathbf T}^*$ be a maximal  $Fr$-stable torus of ${\mathbf L}^*$ such that $T^*={\mathbf T}^{*Fr}$. 
Note that $(N_{C_{  {\mathbf G}^*}(s)}({\mathbf T}^*))^{Fr}=N_{C_{  G^*}(s)}({\mathbf T}^*).$ 
  Observe that $N_{C_{  G^*}(s)}({\mathbf T}^*)\subset C_{  G^*}(s)\subset O(V_s^*)\times O(V_s^{*\perp})$. Furthermore,   $N_{C_{  G^*}(s)}({\mathbf T}^*)$ is  not contained in $X^*=X^*_1X^*_2=SO(V_s^*)\times SO(V_s^{*\perp})$ \ii there are  $g_1\in O(V_s^*)$ and $g_2\in C_{O(V_s^{*\perp})}(s_2)$ such that $\det g_i=-1$ and $g_i$ normalizes ${\mathbf T}_i$, $i=1,2$. 
Equivalently, \ii $N_{X_1^*}({\mathbf T}^*_1)\neq N_{O(V_s^*)}({\mathbf T}_1^*)$ and 
$N_{C_{X_2^*}(s_2)}({\mathbf T}^*_2)
\neq N_{C_{O(V_s^{*\perp})}(s_2)}({\mathbf T}_2^*)$. 
By Lemma \ref{a33}(4), $N_{X_1^*}({\mathbf T}^*_1)\neq N_{O(V_s^*)}({\mathbf T}_1^*)$  
 \ii $T_1$ is not exceptional. 
 Let $g_3$ be the projection of $g$ to $Y$. Then $g_2\in X^*_2$ \ii $g_3\in SO(Y)$.
By Lemma \ref{a33}(4), $N_{C_{X_2^*}(s_2)}({\mathbf T}^*_2)
\neq N_{C_{O(V_s^{*\perp})}(s_2)}({\mathbf T}_2^*)$ 
\ii  $T_3$ is non-exceptional. 
This implies the equivalence of $(1)$ and $ (2)$.

The equivalence of $(1)$ and $ (3)$ is obvious. (Recall that $W_X(T)_\theta\cong N_{C_{X^*}(s)}({\mathbf T}^*)/T^*$.) So (B) follows.  

(C) The statements (3), (4) and (5) are the negations of  (1), (2) and (3), respectively,
except for the additional statement on the index  $|W(T)_\theta:W_{X}(T)_\theta |$. This coincides with the index 
$|N_{C_{G^*}(s)}({\mathbf T}^*):N_{C_{X^*}(s)}({\mathbf T}^*)|$.  As $|C_{G^*}(s):C_{X^*}(s)|\leq 2$, it suffices to observe that $N_{C_{G^*}(s)}({\mathbf T}^*)\neq N_{C_{X^*}(s)}({\mathbf T}^*)$. This is stated in (4).

\bl{om1}  Let $T=T_1\times T_2$, let $i\ra d_i$ be the function defining $T_1$ and $m(T_1)=\sum d_i$. 
Set $t=2^{m(T_1)}$ if $T_1$ is exceptional, 
otherwise $t=2^{m(T_1)-1}$. 
 
If  both the tori $T^*_1,T_3^*$ are non-exceptional  then
$|W(T)_\theta|=2\cdot |W_X(T)_\theta|=2^{m(T_1)}\cdot |W_L(T)_\theta|$, otherwise $|W(T)_\theta|=|W_X(T)_\theta|=t\cdot |W_L(T)_\theta|.$
\el 
The first equality in  both the cases   follows from Lemmas \ref{yr1}  and \ref{dx5}. Furthermore, $W_X(T)_\theta=W_{X_1}(T_1)\times W_{X_2}(T_2)_{\theta_2}$  and $W_L(T)_\theta=W_{G_1}(T_1)\times W_{X_2}(T_2)_{\theta_2}$, as $G_2=X_2$ and $s_1=1$. In addition, $|W_{G_1}(T_1)|=\Pi_i i^{d_i}i!$, whereas $|W_{X_1}(T_1)|=\Pi_i (2i)^{d_i}i!/2$, 
unless all parts of  the partition defining $T_1$ are even, in which case 
$|W_{X_1}(T_1)|=\Pi_i (2i)^{d_i}i!$ (Lemma \ref{a33}).  Therefore, $W_X (T)_\theta=t\cdot W_L(T)_\theta$. 

\bl{ff1}  
Keep the notation of Lemma $\ref{dx5}$.
Denote by $\nu(s,T^*)$ the number of $C_{L^*}(s)$-orbits of tori $T^{\prime *}\subset C_{L^*}(s)$ that are $C_{G^*}(s)$-conjugate to $T^*$. 

$(1)$ Suppose that  $T^*_1$ is exceptional and either $Y=0$ or $T^*_3$ is  exceptional.  Then 
$\nu(s,T^*)=1$ and  $\frac{(\om_G|_T,\theta)}{|W(T)_\theta|}= \frac{1}{ |W_L(T)_\theta|}.$

$(2)$ Suppose that  $T^*_1$ is exceptional but  $T^*_3$ is not. Then 
 $\nu(s,T^*)=1$ and  $\frac{(\om_G|_T,\theta)}{|W(T)_\theta|}=\frac{1}{  |W_L(T)_\theta|};$

$(3)$ Suppose that  $Y\neq 0$ and $T^*_3$ is exceptional  but  $T^*_1$  is not exceptional.  Then 
$\nu(s,T^*)=2$ and  $\frac{(\om_G|_T,\theta)}{|W(T)_\theta|}=\frac{2}{ |W_L(T)_\theta|}.$

$(4)$ Suppose that both $T^*_1$ and  $T^*_3$ are  non-exceptional. Then 
$\nu(s,T^*)=1$ and  $\frac{(\om_G|_T,\theta)}{|W(T)_\theta|}=\frac{1}{  |W_L(T)_\theta|}.$

\smallskip
$(5)$ Suppose that $Y=0$ and $T_1$ is non-exceptional. Then 
$\nu(s,T^*)=1$ and
 $\frac{(\om_G|_T,\theta)}{|W(T)_\theta|}=\frac{2}{ |W_L(T)_\theta|}.$ 
\el In the notation of Lemma \ref{dx5} 
$\nu(s,T^*)$ is the number of $C_{L^*}(s)$-orbits in ${\mathcal S}$. So 
 the formulae for $\nu(s,T^*)$ follow from that lemma. The formulae for $(\om_G|_T,\theta)$  follow from Lemma \ref{om1} and the fact that $(\om_G|_T,\theta)=2^{m(T_1)}$ (Lemma \ref{d2g}). 


\medskip 
We shall use Lemma \ref{ff1} to transform formula (\ref{eq11}). 

\begin{propo}\label{pr8} Let $G=SO (V)$, where $\dim V$ is even,
 and let $s\in G^*$ be a semisimple element such that $V_s^*\neq0$. 

$(1)$ Suppose that $-1$ is not an \ei of s and $\dim U$ is odd.  Then

\begin{equation}\label{eq55}
St_s^+=\sum_{({\mathbf T},\theta)\in {\mathbf R}_s}\frac{(\om_G|_T,\theta)}{|W(T)_\theta|}\ep_{{\mathbf
G}}\ep_{{\mathbf T}} R_{{\mathbf T},\theta}
=2\cdot  \sum_{({\mathbf T},\theta)\in {\mathbf R}^{L}_s}\frac{1}{|W_{L}(T)_\theta|}\ep_{{\mathbf
L}}\ep_{{\mathbf T}} R_{{\mathbf T},\theta}. 
\end{equation}

\med $(2)$ Suppose that $-1$ is not an \ei of s and $\dim U$ is even.  Then

\begin{equation}\label{eq27}
St_s^+=\sum_{({\mathbf T},\theta)\in {\mathbf R}_s}
\frac{(\om_G|_T,\theta)}{|W(T)_\theta|}\ep_{{\mathbf
G}}\ep_{{\mathbf T}} R_{{\mathbf T},\theta}=$$
$$
=\sum_{({\mathbf T},\theta)\in {\mathbf R}^{L_1}_s}
\frac{1}{|W_{L_1}(T)_\theta|}\ep_{{\mathbf
L_1}}\ep_{{\mathbf T}} R_{{\mathbf T},\theta}~~ + \sum_{({\mathbf T},\theta)\in {\mathbf R}^{L_2}_s}
\frac{1}{|W_{L_2}(T)_\theta|}\ep_{{\mathbf
L_2}}\ep_{{\mathbf T}} R_{{\mathbf T},\theta}.
\end{equation}
(\noindent The above expression is understood so that if $T$ is not $C_{G^*}(s)$-conjugate to $L_i$ for $1\leq i\leq 2$
then there is no term in the corresponding sum.)

\med $(3)$ Suppose that $-1$ is an \ei of s. Then 

  \begin{equation}\label{eq54}
St_s^+=\sum_{({\mathbf T},\theta)\in {\mathbf R}_s}\frac{(\om_G|_T,\theta)}{|W(T)_\theta|}\ep_{{\mathbf
G}}\ep_{{\mathbf T}} R_{{\mathbf T},\theta}
= \sum_{({\mathbf T},\theta)\in {\mathbf R}^{L}_s}\frac{1}{|W_{L}(T)_\theta|}\ep_{{\mathbf
L}}\ep_{{\mathbf T}} R_{{\mathbf T},\theta}. 
\end{equation}


\end{propo}
Proof. The terms in (\ref{eq11}) with $(\om_G|_T,\theta)=0$ can be droped. If $(\om_G|_T,\theta)\neq 0$ then, by Lemma \ref{no8} (and the definition of ${\mathcal L}$ prior Lemma \ref{aa2b}), $T^*$ is $C_{G^*}(s)$-conjugate to a torus in some $L\in {\mathcal L}$. If (a) ${\mathcal L}$ consists of a single $C_{G^*}(s)$-orbit then we can fix $L\in {\mathcal L}$. Furthermore, if (b) $ {\mathcal L}$ consists of a single $C_{L^*}(s)$-orbit then the elements of $({\mathbf T},\theta)\in {\mathbf R}_s$ are in bijection with ${\mathbf R}_s^L$. Note that $\ep_{{\mathbf
G}}=\ep_{{\mathbf
L}}$ for every Levi subgroup $\ep_{{\mathbf
L}}$ of $\ep_{{\mathbf
G}}$.

Suppose that  (1) holds. 
 As $\dim U$ is odd,   (a) holds by Lemma \ref{no9}. As $-1$ is not an \ei of $s$,  (b) holds too.   
So the result follows  by Lemma \ref{ff1}(5). 

$(2)$ 
In this case
we have two Levi subgroups $L^*_1,L^*_2\in {\mathcal L}$  which are not
 conjugate in $C_{G^*}(s)$, see Lemma \ref{2es}.  So a maximal  torus $T^*$ of $C_{G^*}(s)$ is $C_{G^*}(s)$-conjugate either (i) to a torus of $L^*_1$ but not $L^*_2$, or  to a torus of $L^*_2$ but not $L^*_1$, or (ii)   to a torus of each $L_1^*, L_2^*$. The  option (i) happens if and only if  $T^*_1$ is exceptional (Lemma \ref{aa2b}). 

 If $T^*_1$ is not exceptional then we are in the case (5)  of Lemma \ref{ff1}. As $|W_{L_2}(T)_\theta|=|W_{L_2}(T)_\theta|$, when $T$ is viewed as a torus in $L_1^*$ or $L_2^*$,  we have:  $$\frac{(\om_G|_T,\theta)}{|W(T)_\theta|}\ep_{{\mathbf
G}}\ep_{{\mathbf T}} R_{{\mathbf T},\theta}=
\frac{1}{|W_{L_1}(T)_\theta|}\ep_{{\mathbf
L}_1}\ep_{{\mathbf T}} R_{{\mathbf T},\theta} +  
\frac{1}{|W_{L_2}(T)_\theta|}\ep_{{\mathbf
L}_2}\ep_{{\mathbf T}} R_{{\mathbf T},\theta} .$$   

If $T^*_1$ is exceptional then $T^*$ is $C_{G^*}(s)$-conjugate  to a torus either in $L^*_1$ or in $L^*_2$, but not in both of them. 
As we are in case (1) of Lemma \ref{ff1},  the term with this $T$ occurs only in one of the sums in (\ref{eq27}). So again 
(by the convention in the parentheses of (2)) we have: 
$$\frac{(\om_G|_T,\theta)}{|W(T)_\theta|}\ep_{{\mathbf
G}}\ep_{{\mathbf T}} R_{{\mathbf T},\theta} =\frac{1}{|W_{L_1}(T)_\theta|}\ep_{{\mathbf
L}_1}\ep_{{\mathbf T}} R_{{\mathbf T},\theta} +  
\frac{1}{|W_{L_2}(T)_\theta|}\ep_{{\mathbf
L}_2}\ep_{{\mathbf T}} R_{{\mathbf T},\theta} .$$ This implies the second equality of (\ref{eq27}) as 
  $\ep_{{\mathbf L}_1}\ep_{{\mathbf L}_2}=\ep_{{\mathbf G}}$.

\medskip
(3) As $-1$ is an \ei of $s$, all  subgroups of ${\mathcal L}$ are $C_{G^*}(s)$-conjugate by Lemma \ref{no9}. 
So $L\in {\mathcal L}$ can be fixed. 
If $\dim U$ is odd, then    either (3) or (4) of Lemma \ref{ff1} holds. If (3) of Lemma \ref{ff1} holds then ${\mathcal L}$ consists of two $C_{L^*}(s)$-orbits, and hence the element $({\mathbf T}, \theta)\in{\mathbf R}_s$ corresponds to two elements $({\mathbf T}, \theta), ({\mathbf T}', \theta')$ in ${\mathbf R}^L_s$. Therefore,

  \begin{equation}\label{eq77}\frac{(\om_G|_T,\theta)}{|W(T)_\theta|}\ep_{{\mathbf
G}}\ep_{{\mathbf T}} R_{{\mathbf T},\theta} =\frac{2}{|W_{L}(T)_\theta|}\ep_{{\mathbf
G}}\ep_{{\mathbf T}} R_{{\mathbf T},\theta}=\frac{1}{|W_{L}(T)_\theta|}\ep_{{\mathbf
L}}\ep_{{\mathbf T}} R_{{\mathbf T},\theta} +  
\frac{1}{|W_{L}(T')_{\theta'}|}\ep_{{\mathbf
L}}\ep_{{\mathbf T}'} R_{{\mathbf T}',\theta'} .\end{equation}
If (4) of Lemma \ref{ff1} holds then ${\mathcal L}$ consists of a single $C_{L^*}(s)$-orbit. Therefore, we have

  \begin{equation}\label{eq88}
\frac{(\om_G|_T,\theta)}{|W(T)_\theta|}\ep_{{\mathbf
G}}\ep_{{\mathbf T}} R_{{\mathbf T},\theta} =\frac{1}{|W_{L}(T)_\theta|}\ep_{{\mathbf
G}}\ep_{{\mathbf T}} R_{{\mathbf T},\theta}.\end{equation}
So (\ref{eq54}) holds.

\med
Finally suppose that $\dim U$ is even. 
 We show that (\ref{eq54}) holds. 
 
Let ${\mathcal S}$ be as in Lemma \ref{dx5}. By Lemma \ref{dx5}, ${\mathcal S}$ consists of two $C_{L^*}(s)$-orbits 
\ii $T_1^*$ is non-exceptional and $T^*_3$ is exceptional. In this case (3) of Lemma \ref{ff1} holds.  Then 
 there is a single class of $C_{G^*}(s)$-conjugacy class of tori  $T\subset C_{L^*}(s)$, 
which splits in two $C_{L^*}(s)$-conjugacy classes. As  $\frac{(\om_G|_T,\theta)}{|W(T)_\theta|}=\frac{2}{ |W_L(T)_\theta|}, $ we  get (\ref{eq77}).

So we are 
left with the situation where ${\mathcal S}$ consists of a single $C_{L^*}(s)$-orbit. Then we have one of the cases 
(1), (2) or (4) of Lemma \ref{ff1}.

 If (1) of Lemma \ref{ff1} holds then there are two $C_{G^*}(s)$-conjugacy classes of tori  $T\subset C_{L^*}(s)$, 
and each of them form a single class of $C_{L^*}(s)$-conjugate tori.  So   (\ref{eq88}) holds.

Suppose that  (2) of Lemma \ref{ff1} holds.  Then there is a single class of $C_{G^*}(s)$-conjugacy class of 
tori  $T\subset C_{L^*}(s)$, and all these tori are  $C_{L^*}(s)$-conjugate.  So again  (\ref{eq88}) holds.

Suppose that  (4)  of Lemma \ref{ff1} holds.  Then there is a single class of $C_{G^*}(s)$-conjugacy class 
of tori  $T\subset C_{L^*}(s)$, and all these tori are  $C_{L^*}(s)$-conjugate.  So   (\ref{eq88}) holds.

This implies the statement (3) of the  proposition.  

\med 
Next we show that the right hand side of the equalities in Lemmas \ref{pr9} and  \ref{pr8} can be expressed in terms 
of regular characters of $L$. Recall that if $L$ is a Levi subgroup of a parabolic subgroup $P$ of $G$ and 
$\lam$  is a character of $L$ then $\lam^{\#G}$ denotes  the Harish-Chandra induced character. This is exactly  the induced character $\lam_1^{G}$, where $\lam_1$ is the inflation of $\lam$ to $P$ via the projection $P\ra L$. Note that 
$\lam^{\#G}$  does not depend on the choice of   $P$ \cite[70.10]{CR2}. Note that if $L$ is abelian then every character of $L$ is regular (by convention). 

Similarly to the usage of the notation ${\mathbf R}_s$, for $s\in L^*$ we denote by ${\mathbf R}^L_s$ the set of representatives of the $L$-conjugacy classes of pairs $({\mathbf T},\theta)$ with ${\mathbf T}\subset {\mathbf L}$. They are in bijection with   $C_{L^*}(s)$-conjugacy classes of maximal tori $T^*\subset C_{L^*}(s)$. The argument is based on the fact that $ {\mathbf T}\subset  {\mathbf L}$ then $R_{{\mathbf T},\theta}=(R^L_{{\mathbf T},\theta})^{\# G}$, see  \cite[7.4.4]{C}.

\begin{propo}\label{tt7}  Let $G=SO(V)$ and let $s\in G^*$ be a semisimple element such that $St^+_s\neq 0$ and $V^*_s\neq 0$.
 Let $U$ be a maximal totally singular (or totally isotropic) subspace of $V^*_s$, $P_U$ the stabilizer of $U$ in $G^*$
and  $L^*$ a Levi subgroup of $P_U$ such that $L^*V^*_s=V^*_s$. 
Then 

$$\sum_{({\mathbf T},\theta)\in {\mathbf R}^{L}_s}\frac{1}{|W_{L}(T)_\theta|}\ep_{{\mathbf
L}}\ep_{{\mathbf T}} R_{{\mathbf T},\theta}=\big( \sum_{({\mathbf T},\theta)\in {\mathbf R}^{L}_s}\frac{1}{|W_{L}(T)_\theta|}\ep_{{\mathbf
G}}\ep_{{\mathbf T}} R^L_{{\mathbf T},\theta}\big)^{\#G}
$$ and  $$\lam_s:=
\sum_{({\mathbf T},\theta)\in {\mathbf R}^{L}_s}\frac{1}{|W_{L}(T)_\theta|}\ep_{{\mathbf
L}}\ep_{{\mathbf T}} R^L_{{\mathbf T},\theta}$$
is a regular \ir character of $L$. In addition, $\lam_s$ is a unique regular character of ${\mathcal E}_s^L$.
\end{propo}
Proof.  
Let  $u=\dim U$. Observe that the group $C_{{\mathbf L}^*}(s)$ is connected. Indeed, 
if $\dim V$ is odd then ${\mathbf L}^*\cong GL_u( \overline{F}_q)\times   Sp_{2(n-u)}(  \overline{F}_q)$, and hence  
$C_{{\mathbf L}^*}(s)\cong GL_u (\overline{F}_q)\times C_{Sp_{2(n-u)}(  \overline{F}_q)}(s_2)$. So   the claim is true as the latter group is simply connected  (Lemma   \ref{os4}). If  $\dim V$ is even  then ${\mathbf L}^*={\mathbf G}_1^*\times   {\mathbf G}_2^*$, where ${\mathbf G}_1^* \cong GL_u( \overline{F}_q)$ and ${\mathbf G}_2^*\cong  SO^\pm _{2(n-u)}(  \overline{F}_q)$. Therefore, $C_{{\mathbf L}^*}(s)$ is connected \ii $C_{{\mathbf G}_2^*}(s_2)$ is connected. Again, this is the case by Lemma   \ref{os4}.

The expression for $\lam_s$ coincides with the class function on $L$ defined in \cite[14.40]{DM}. As $C_{{\mathbf L}^*}(s)$ is connected, it follows, by \cite[14.43]{DM}, that $\lam_s$ is the unique  regular  character of ${\mathcal E}_s^L$. 

 By \cite[7.4.4]{C}, if $ {\mathbf T}\subset  {\mathbf L}$ then $R_{{\mathbf T},\theta}=(R^L_{{\mathbf T},\theta})^{\# G}$, where
$R^L_{{\mathbf T},\theta}$ is a Deligne-Lusztig character of $L$. 
 Note that  $\ep_{{\mathbf
L}}=\ep_{{\mathbf G}}$ and $\ep_{{\mathbf T}}$ is unchange when we view ${\mathbf
T}$ as a torus of ${\mathbf L}$. So  the proposition follows. 

\med
The result of Proposition \ref{tt7} is not sufficient to prove that $\lam_s^{\#G}$ is \mult free. For this we shall additionally show that $\lam_s=St_{G_1}\otimes \rho_{s_2}$, where $\rho_{s_2}$ is a regular character of $G_2$.

\begin{propo}\label{tt8}  Let  $\lam_s$ be as in Proposition  $\ref{tt7}$. Then $\lam_s=
St_{G_1}\otimes \rho_{s_2}$, where $\rho_{s_2}$ is a regular character of $G_2$ from the Lusztig series ${\mathcal E}_{s_2}$.          
\end{propo}

Proof.  Recall that ${\mathbf L}={\mathbf G}_1\times {\mathbf G}_2$, where $ {\mathbf G}_1\cong GL_r( \overline{F}_q)$ and $ {\mathbf G}_2\cong SO_{\dim V-2r}( \overline{F}_q)$. Then  $T=T_1\times T_2$ and $\theta=\theta_1\otimes \theta_2$, where  $T_i$ is a maximal torus of $G_i$ and $ \theta_i$ is a linear character of $T_i$, $i=1,2$.

It is well known that $R_{{\mathbf T},\theta}^L=R^{G_1}_{{\mathbf T}_1,\theta_1}\otimes R^{G_2}_{{\mathbf T}_2,\theta_2}$, $\ep_{{\mathbf L}}=\ep_{{\mathbf G}_1}\cdot \ep_{{\mathbf G}_2}$ and $\ep_{{\mathbf T}}=\ep_{{\mathbf T}_1}\cdot \ep_{{\mathbf T}_2}$. Furthermore, $W_L(T)_\theta=W_{G_1}(T_1)_{\theta_1}\times W_{G_2}(T_2)_{\theta_2}=W_{G_1}(T_1)\times W_{G_2}(T_2)_{\theta_2} $ as $\theta_1=1_{T_1}$ and $s_2$
does not have \ei 1. \itf 

$$\sum_{({\mathbf T},\theta)\in {\mathbf R}^L_s}\frac{\ep_{{\mathbf G}}\ep_{{\mathbf T}}}{|W_{G}(T)|}R_{{\mathbf T},\theta}^L=
\big(\sum_{({\mathbf T}_1,\theta_1)\in {\mathbf R}^{{\mathbf G}_1}_{s_1}}\frac{\ep_{{\mathbf G}_1}\ep_{{\mathbf T}_1}}{|W_{G_1}(T_1)|} R^{G_1}_{{\mathbf T}_1,\theta_1}\big ) \otimes \big(\sum_{({\mathbf T}_2,\theta_2)\in {\mathbf R}^{{\mathbf G}_2}_{s_2}}\frac{\ep_{{\mathbf
G}_2}\ep_{{\mathbf T}_2}}{{|W_{G_2}(T_2)_{\theta_2}|}} R^{G_2}_{{\mathbf T}_2,\theta_2}\big). $$
The expression in the first parentheses yields the Steinberg character of $G_1$, see \cite[7.6]{C}. 
Therefore, $\lam'=St_{G_1}$ and $\lam''$ equals the expression in the second parentheses. The latter coincides with the class function on $G_2$ defined in \cite[14.10]{DM}. As $\lam''$ is irreducible, it is exactly the unique  regular character of $G_2$ that belongs to the Lusztig series ${\mathcal E}^{G_2}_{s_2}$, see \cite[14.49]{DM}.  (Alternatively, as $s=\diag(s_1,s_2)$, the regular character of the Lusztig series ${\mathcal E}_s$ is the product of  the regular characters in  ${\mathcal E}^{G_1}_{s_1}$ and ${\mathcal E}^{G_2}_{s_2}$. As  $s_1=\Id$, the regular characters in  ${\mathcal E}^{G_1}_{s_1}$ is $St_{G_1}$, whence the result.)

\med
Recall that if $-1$ is not an \ei of $ s$, $\dim V$ is even and $\dim V_s^*\equiv 0 \pmod 4$ (that is, $\dim U$ is even), then ${\mathcal L}$ consists of  two $C_{G^*}(s)$-orbits (Lemma \ref{no9}). We choose   $L_1^*,L_2^*$ from distinct orbits, and  write $L_1^*=G_1^{\prime *}\times G_2^*$ and $L_2^*=G_1^{\prime\prime  *}\times G_2^*$. (Note that $\rho  _{s_2} $, as defined in Proposition \ref{tt8}, is 
the same  in both the cases, as $\rho  _{s_2} $ depends only on  $s_2$.)  In the other cases ${\mathcal L}$ forms   a single  $C_{G^*}(s)$-orbit. 

\begin{propo}\label{mr9}  Let $G=SO(V)$, and let $s\in G^*$ be a semisimple element. 
Let ${\rm Spec}\, s$ denote the set of \eis of s on $V^*$.
Suppose that  $V^*_s\neq 0$ and  $St_s^+\neq 0$.

$(1)$ If  $\dim V$ is odd then $St_s^+=(St_{G_1}\otimes \rho_{s_2})^{\# G}$.

$(2)$ If   $\dim V$ is  even then 
$$St_s^+=\begin{cases}
(St_{G_1}\otimes \rho_{s_2})^{\# G}&if ~-1\in {\rm Spec}\, s, \\
2\cdot (St_{G_1}\otimes \rho_{s_2})^{\# G}&if 
~-1\notin {\rm Spec}\, s~ {\rm and}~\dim V_s^*\equiv 2 \pmod 4,
\cr
((St_{G'_1}+St_{G''_1})\otimes \rho_{s_2})^{\# G}&if~
~-1\notin {\rm Spec}\, s ~{\rm and}~\dim V_s^*\equiv 0 \pmod 4.
\end{cases}$$\end{propo}
Proof. This follows from Lemma \ref{pr8} and Proposition \ref{tt8}. (Note that $\rho_{s_2}$ has to be omitted if $s=1$.)
 
\med
If $q$ is odd, we denote by $1_G^-$ the only non-trivial one-dimensional character of $G$, and set $St_G^-=St_G\otimes 1_G^- $. It is well known that $St_G$ and $St_G^-$ are the only \ir characters of $G$ of defect 0.

\begin{corol}\label{co4} $(1)$ Let $G=SO_{2n+1}(q)$ , $n>0$. Then $(St^+_G,St_G)=1=(St^+_G,St^-_G)$, and $(St^+_G,1_G)\neq 0$ \ii $n=1$. In addition, $(St^+_G,1^-_G)=0$. 

$(2)$ Let  $G=SO^\al_{2n}(q)$, $n>1$.  Then  $(St^+_G,St_G)=1+\al$, and $(St^+_G,1_G)= 0=(St^+_G,1^-_G)$. In addition,
 $(St^+_G,St^-_G)=1$.
\end{corol}
Proof. It is well known that $1_G,St_G\in {\mathcal E}_1$, that is, $s=1$. If $q$ is odd then  $1^-_G,St^-_G\in {\mathcal E}_s$ for $s=-\Id$.  
So $1^-_G$ is a constituent of $St^+_G$ \ii $St^-_G$ is one-dimensional, and hence $G$ is abelian. As $n>1$, this is not the case. So the claims about $1^-_G$ follow. 

 By Proposition \ref{pp4}, $St^+_{-\Id}$ is the regular character 
in ${\mathcal E}_{-\Id}$, which is  exactly $St^-_G$.  It is also well known that the Harish-Chandra restriction of $St_G$ to any Levi subgroup $L$ of $G$ is $St_L$, see \cite[p. 72]{DM}.
If $s=1$ then $L=G_1$.

(1) In this case $St^+_1=St_{G_1}^{\#G}$, where $G_1\cong GL_n(q)$.   By Harish-Chandra reciprocity,  $(St_{G_1}^{\#G},St_G)=(St_{G_1},St_{G_1})=1$. In addition,    $(St_{G_1}^{\#G},1_G)=(St_{G_1},1_{G_1})  $, which is non-zero \ii $St_{G_1}=1_{G_1}$, that is, when  $G_1$ is abelian. This implies $n=3$. 

(2)  Let $\al=-1$. Then  $St^+_1=0$ by Lemma \ref{pp5}. So $1_G,St_G$ are not constituents of $St^+_G$. 

Let  $\al=1$. Then $St^+_1=2\cdot St_{G_1}^{\#G}$ if $\dim V^*/2$ is odd, and $St_{G'_1}^{\#G}+St_{G''_1}^{\#G}$ otherwise. Here $G_1\cong GL_n(q)$, as well as $G_1'$ and $G_1''$. As  above, $St_G$ is a constituent of each 
 $St_{G'_1}^{\#G}$ and $St_{G''_1}^{\#G}$ by Harish--Chandra reciprocity, whereas $1_G$ is a constituent of neither
 $St_{G'_1}^{\#G}$ nor $St_{G''_1}^{\#G}$, unless $G'_1\cong G_1''$  is a torus. The letter implies $n=1$. 

\section{The decomposition of $(St_{G_1}\otimes \rho_{s_2})^{\#G}$}
 
Let $U\subset V^*$ be a totally singular (or totally isotropic) subspace, $P_U$ its stabilizer in  $G^*$ and  $L^*$ a Levi subgroup of $P_U$. Then $L^*=G^*_1\times G^*_2$, where $G^*_1 \cong GL(U)$, $G^*_2\cong  SO(V^{\prime *})$ and $V^{\prime *}$ be a complement of $U$ in $U^\perp$. 
These correspond to respective objects in $G$: there is a totally singular subspace $R$  of $V$  of dimension $\dim U$   whose stabilizer in $G$ is a parabolic subgroup and its Levi subgroup $L$ is dual to 
$L^*$. If we set $V'=R^\perp/R$ then $SO(V')$ is dual to $G_2^*$. 

Let $\rho$ be a regular character of $SO(V')$. In this section  we show that $(St_{G_1}\otimes \rho )^{\#G}$ is a  \mult free character of $G$. We denote by $S_n$ the symmetric group of permutations  of $n$ objects. 

We start with a special case where 
 $G^*=SO^+_{2n}(q)$, $n$ even, and $s=1$.     Let  
  $U_1,U_2$ be two totally singular subspaces of $V^*$ of dimension $n$ such that $gU_1\neq U_2$ for any $g\in G^*$ (see Lemma \ref{no6}). Let $P_i$ be the stabilizer of 
$U_i$ for $i=1,2$. Then  $P_1,P_2$ are non-conjugate parabolic subgroups of $G $.  Let   $L^*_1,L^*_2$ be Levi subgroups of $P_1,P_2$, respectively.  By  Lemma \ref{no9}  for $s=1$,  the groups  $L^*_1,L^*_2$ are not conjugate in $G^*. $ Let $W$ be the Weyl group of type $D_n$. Then the  Weyl groups $W_1,W_2$ of  
$L^*_1,L^*_2$, respectively, can be viewed as subgroups of $W$, each isomorphic to $S_n$, 
the symmetric group. Groups $W_1,W_2$ are known to be conjugate in $W(B_n)$, but not in $W$ if $n$ is even. 


\begin{lemma}\label{aa6}  Let $W$ be the Weyl group of type $D_n$ for $n>2$ even.

$(1)$  $W_1,W_2$ are not conjugate in $W;$ 

$(2)$ 
$(1_{W_1}^W,1_{W_2}^W)=n/2$ and  $(1_{W_1}^W,1_{W_1}^W)=(n+2)/2$.  

\end{lemma}

Proof. View  $W$ as a group of monomial matrices with non-zero entries $\pm 1$, see comments prior Lemma \ref{ee1}.
Let  $D$ be the subgroup  of diagonal matrices in $W$. Then $D$  is of exponent $2$,  of  order $2^{n-1}$, and $\det d=1$ for every $d\in D$.
In addition, $W$ contains a subgroup $W_1$ consisting of monomial matrices with entries $0,1$ and $W=DW_1=W_1D$. 
Let $t=\diag(-1,1\ld 1)$. Then $t\notin W$, but $tWt\up=W$. We set $W_2=tW_1t\up$.
We first observe that $W_1,W_2$ are not conjugate in $W$. Indeed, if $W_2=xW_1x\up$ for $x\in W$ then 
$tx$ normalizes $W_1.$ Let $x=ds$, where $d\in D$, $s\in W_1$. So $td$  normalizes $W_1.$
One easily observes that a diagonal matrix normalizing $W_1$ must be scalar. Therefore,  $t=\pm d$, which is false as $\det d=\det (-d)=1$ whereas $\det t=-1.$

Thus, $W_1,W_2$ are not conjugate in $W$.  We  show that $(1_{W_2}^W,1_{W_1}^W)=n/2$.
The fact that $(1_{W_1}^W,1_{W_1}^W)=(n+2)/2$  is known, and also follows from our resoning below.
It is well known that  $(1_{W_2}^W,1_{W_1}^W)$ equals the number of the orbits of  $W_2$ on the
cosets $gW_1$. So we proceed with computing these orbits.

As $W=DW_1$, it follows that representatives of cosets $gW_1$ can be chosen in $D$. Obviously,
all cosets $dW_1$ are distinct when $d$ runs over $D$. This  implies that $W_1dW_1=W_1d'W_1$
\ii the ranks of the matrices $d-\Id$ and $d'-\Id$ coincide. The ranks may be equal $2i$ with $i=0\ld n/2$, so 
$(1_{W_1}^W,1_{W_1}^W)=(n+2)/2$.

Let $Y$ be the subgroup of $W_1$ consisting of matrices with $1$ at the $(1,1)$-position. 
Then $Y\cong S_{n-1}$ and $[t,Y]=1$. Note that $W_1=\cup_{i=1}^n a_iY$, where $a_i\in W_1$ is the matrix 
whose non-zero non-diagonal entries are at the $(1,i)$ and $(i,1)$-positions. Clearly, $Y$ is also contained in $W_2$.
It is easy to observe that representatives of the double cosets $YxW_1$ ($x\in W$) can be chosen 
to be $\diag(\pm 1, 1\ld 1, -1\ld -1)$, where the number of entries $-1$ may be 0. As $W_2$
contains the matrix whose entries  at the positions $(1,i)$ and $(i,1)$ are  $-1$ and $0$ or $1$ elsewhere, it follows that 
representatives of the double cosets $YxW_1$ ($x\in W$) can  be chosen 
to be of shape $d_i=\diag(1, 1\ld 1, -1\ld -1)$, where the number of 1-entries is even and ranges between $1$ and $n/2$.
Moreover, one observes that all double cosets $W_2d_iW_1$ are distinct, and hence $d_1\ld d_{n/2}$ are representatives of all double cosets $W_1xW_1$ ($x\in W)$. Therefore, $(1_{W_1}^W,1_{W_2}^W)=n/2$, as claimed.

As both $1_{W_1}^W$ and  $1_{W_2}^W$ are \mult free, it follows that  $1_{W_1}^W$ and  $1_{W_2}^W$ 
have exactly $n/2$ common \ir constituents. 

\medskip
Remark. One can identify the characters of $W(D_n)$ occurring in both $1_{W_1}^W$ and $1_{W_2}^W$. 
In the bijection between the conjugacy classes and irreducible characters of $W(D_n)$ described in  
\cite[Section 5]{GP} they correspond to the partitions  $[(k)],[(l)]$, where $k+l=n$ and $k\neq l$.
(This follows from \cite[Proposition 6.1.5]{GP}, which is more precise than our elementary lemma \ref{aa6}.)
So $k$ runs over the numbers $0,2 \ld n$ except $n/2$ (so $l=n-k$), one obtains exactly $n/2 $ characters.

\bl{in1a}   
With above notation, let $St_i$ be the Steinberg character of $L_i$.
Then $St_1^{\# G}$ and $St_2^{\# G}$ are \mult free characters, having $n/2$ common \ir constituents. 
In addition, $St_G$ is a constituent of both $St_1^{\# G}$ and $St_2^{\# G}$.  \el
It is known \cite[70.24]{CR2}  that the \ir constituents of $St_i^{\# G}$ for $i=1,2$ are in bijection with the \ir constituents of $\ep_i^W$, where 
$W=W(D_n)$,  $W_1,W_2$ are as in Lemma \ref{aa6} and $\ep_i$ is the sign character of $W_i$. 
Moreover, the bijection regards the multiplicities of the \ir constituents. Note that $\ep_i^W=1_{W_i}^W\otimes \nu_i$, where $\nu_i$ is a linear character of $W$ such that $\nu_i|_{W_i}=\ep_i$ (there are two such characters $\nu_i$ for every $i=1,2$). So the first statement of the lemma  follows from Lemma \ref{aa6}. The proof of  Corollary \ref{co4}  contains a proof of the additional  statement.

\medskip
Let $G$ be a group with BN-pair, which later will be specified to be $SO(V)$. Let $W(G)$ be the Weyl group of the BN-pair. In terms of algebraic groups $W$ coincides with $W(T)$ as defined above for a maximal torus $T$ of a Borel subgroup of $G$.  By Harish-Chandra theory, every \ir character $\lam$  of $L$ is a constituent of $\delta^{\#L}$, where $\delta$ is a cuspidal character of a  Levi  subgroup $M$, say, of some parabolic subgroup $Q$ of $L$. 
Furthermore, 
the  decomposition of $\delta ^{\#G}$ as a sum of irreducible constituents
is described in terms of the group $W(\delta):=\{w\in W(G): w(M)=M$ and $w(\delta)=\delta\}$, see \cite{HL}. 
This can be also applied to 
$\Irr L$, so we define $W_L(\delta)=\{w\in W(L): w(M)=M$ and $w(\delta)=\delta\}$. The standard embedding of $L$ into $G$ yields an embedding $W_L$ into $W_G$, and hence   
$W_L( \delta)=W(\delta)\cap W_L$.

Let $\zeta\in \Irr G$ be a constituent of $\lam^{\#G}$. Then $\zeta$ is a constituent of $\delta^{\#G}$ by transitivity of Harish-Chandra induction. By the above we can label $\zeta$ and $\rho$ by elements of $\Irr W(\delta)$ and $\Irr W_L(\delta)$, respectively. So let $\zeta=\zeta_\mu$ and $\lam=\lam_\nu$ for some $\mu\in \Irr W(\delta)$ and $\nu\in \Irr W_L(\delta)$. Then $(\zeta,\lam^{\#G})=(\mu,\nu^{W(\delta)})$ by   Howlett and Lehrer \cite[Theorem 5.9]{HL} and Geck \cite[Corollary 2]{Ge2}. 

If $\lam$ is a regular character of $L$ then $\delta$ is a regular character of $M$, see for instance \cite[Proposition 2.5]{Ko}.
Note that $ \delta$ is not unique in general, however, if $\lam\in  {\mathcal E}^L_{s}$ then $\delta$ can be chosen in $ {\mathcal E}^M_{s}.$  (Indeed, let $\Gamma_L$ be a Gelfand-Graev  character containing  $\rho$ as a constituent.
By \cite[Theorem 2.9]{DLM},  the Harish-Chandra restriction $\Gamma_M$ of $\Gamma_L$ to $M$ is a  Gelfand-Graev character of $M$.
Let  $\lam_M$ denote the  Harish-Chandra restriction of $\lam$ to $M$. (In general,  $\rho_M$  is reducible.) By  Harish-Chandra reciprocity, $\delta$ is a constituent of $\lam_M$. Let $ {\mathcal E}^M_{s'}$ be the rational  Lusztig series  of $\Irr M$ containing $\lam_M$. Then the \ir constituents of $\delta^{\#L}$ belongs to  $ {\mathcal E}^L_{s'}$ (see, for instance, 
\cite[Theorem 8.25]{CE}). In particular,  $\lam\in  {\mathcal E}^L_{s'}$, and hence $ {\mathcal E}^L_{s'}= {\mathcal E}^L_{s}$. \itf $s,s'$ are conjugate in $L$. Therefore,
we can take $s'$ for $s$.) 

Let now $G=SO(V)$ and let $L$ be a Levi subgroup of $G$  defined in the beginning of this section. We shall  prove that $\nu^{W(\delta)}$ is \mult   free for $\lam=\lam_s=St_{G_1}\otimes \rho_{s_2}$. 
As $L$ is the direct product of $G_1$ and $G_2$, it follows that $M=M_1\times M_2$ and $\delta=\delta_1\otimes \delta_2$, where $M_i$ (for $i=1,2$) is a Levi subgroup of some parabolic subgroup of $G_i$ and $\delta_i$ is a cuspidal character of $M_i$. We have to show that $(St_{G_1}\otimes \rho_{s_2})^{\#G}$ is a \mult free character of $G.$ We shall do this by analysis of $(\delta_1\otimes \delta_2)^{\#G}$
(this is meaningful as  $M$ is a quotient group of a parabolic subgroup of $G$ by a unipotent normal subgroup). 
We first obtain some information on the  decomposition of $(\delta_1\otimes \delta_2)^{\#G}$ (specified to our situation
where $\lam=\lam_s=St_{G_1}\otimes \rho_{s_2}$).

An advantage of our situation 
is that $W(\delta)\subset W_X=W_{X_1}\times W_{X_2}$.
As $St_{G_1}$ is well known to be a constituent of $1_{T_0}^{\#G_1}$, we have $\delta=1_{T_0}\otimes \delta_2$, where $\delta_2$ is a cuspidal character of $M_2$.

\bl{ch6} 
$W( \delta)\subset W(X_1)\times W(X_2)$.
\el  By the definition of $W(\delta)$ above, if $w\in W(\delta)$ then $w(M)=M$. 
By the comments prior the lemma, we can assume that $\delta$ is a regular character of $M$ and
  $\delta\in {\mathcal E}^M_{s}$. Recall that $C_{{\mathbf L}}(s)$ is connected (see the proof of Proposition \ref{tt7}), and hence 
 so is $C_{{\mathbf M}}(s)$ by Lemma \ref{os4}. \itf that $\delta$ is the only regular character in  $ {\mathcal E}^M_{s}$, 
see \cite[14.40, 14.43]{DM}, and

$$\delta=\sum _{({\mathbf T}, \theta)\in {\mathbf R}^{M}_{s}} a_{{\mathbf T}, \theta}\cdot R_{{\mathbf T}, \theta}^{M},$$
where $R_{{\mathbf T}, \theta}^{M}$ are Deligne-Lusztig characters of $M$ and $a_{{\mathbf T}, \theta}$ are rational numbers,
specified in \cite[14.40]{DM}.

 Let  $w\in W(\delta)$. Then $w$ acts on $M$ as an automorphism, and  this induces an action of $w$ on the characters of $M$. Therefore, 

$$\delta=w(\delta)=\sum _{({\mathbf T}, \theta)\in {\mathbf R}^{M}_{s}} a_{{\mathbf T}, \theta}\cdot  w\big(R_{{\mathbf T}, \theta}^{M}\big).$$
Observe that 
$w(R_{{\mathbf T}, \theta}^{M})=R_{w({\mathbf T}), w(\theta)}^{M}$.  Therefore,
$w(R_{{\mathbf T}, \theta}^{M})$ is a Deligne-Lusztig character of $M$. Furthermore, distinct Deligne-Lusztig characters of $M$ are orthogonal, and if 
$R_{{\mathbf T}, \theta}^{M} =R_{{\mathbf T}', \theta'}^{M}$ then $({\mathbf T}', \theta')\in {\mathbf R}^{M}_{s}$, see \cite[11.15]{DM}.
So $\delta=w(\delta)$   implies $(w({\mathbf T}), w(\theta))\in {\mathbf R}^{M}_{s}$.

Let $T^*$ be a torus in $M^*$ dual to $T$. (Note that  $T^*$ is also dual to $T$ under the duality  $G\ra G^*$ as $M$ is a Levi subgroup of $G$.) Recall that $({\mathbf T}, \theta)\in {\mathbf R}^{M}_{s}$ is equivalent to saying that $s\in T^*$ up to conjugacy in $M^*$. Therefore, we can choose $T^*\subset M^*$ so that  $s\in T^*$.  According with Section 3,
 let $V^*=V_1^*\oplus\cdots \oplus V^*_k \oplus V^*_{k+1}\oplus \cdots \oplus V^*_{k+l} $ be the  $T^*$-decomposition of $V^*$, and 
$T^*=T_1^*\times \cdots \times  T^*_k\times    T^*_{k+1}\oplus \cdots \times  T^*_{k+l} $ the respective decomposition of $T^*$. Correspondingly, let $s=(s_1\ld s_{k+l})$, where $s_i\in T^*_i$ for $i=1\ld k+l$. 
Recall that 
$V^*_s$ denotes the $1$-eigenspace of $s$ on $V^*$. So $V^*_s$ is the sum of the $V^*_i$ that 
are contained in $V^*_s$. Set $J=\{i: V_i\subset V^*_s\}$. Then  $s_i=1$ \ii 
$i\in J$. Taking the dual decomposition $T=T_1\times\cdots\times T_{k+l}$ and $\theta=(\theta_1\ld \theta_{k+l})$, where 
$\theta_i$ is the linear character of $T_i$ corresponding to $s_i$ $(i=1\ld k+l)$, we observe that $\theta_i(T_i)=1$ \ii $i\in J$. 

As $(w({\mathbf T}), w(\theta))\in {\mathbf R}^{M}_{s}$, we can apply the previous paragraph reasoning to $(w({\mathbf T}), w(\theta))$. We conclude that $w(\theta_i)(w(T_i))=1$ implies that $w(
T_i)^*$ acts trivially on $V^*_s$, and hence $w(T_i)^*\subset X_1^*$. (We may view $T_i$ as a subtorus of $T$, and $T^*_i$ as a subtorus of $T^*$.)  In turn, this implies $w(T_i)\subset X_1$ for $i\leq k$, and the lemma follows.

\begin{corol}\label{ch7}  $W(\delta)=W(X_1)\times W_2(\delta_2)$, where $W_2(\delta_2)=\{w\in W(X_2): w(M_2)=M_2, ~w(\delta_2)=\delta_2\}$. Similarly, $W_L(\delta)=W(G_1)\times W_2(\delta_2)$.

\end{corol}
Proof. By Lemma \ref{ch6}, $W(\delta)\subset W(X_1)\times W(X_2)$, and hence $W(\delta)=W_1(\delta_1)\times W_2(\delta_2)$.
It is well known that $St_{G_1} $ is a constituent of $(1_{T_0})^{\#G_1}$, where $T_0$ is a maximal torus of a Borel subgroup of $G_1$ and $1_{T_0}$ is the trivial character of $T_0$. Thus, here $M_1=T_0$ and $\delta_1=1_{T_0}=1_{M_1}$. 
\itf $W_1(\delta_1)$ coincides $W(X_1)$. 
This implies the lemma. 

\medskip
Now we come back to proving that $\lam_s^{\#G}$ is \mult free. As explained prior Lemma \ref{ch6}, it suffices
to show that $\nu^{W(\delta)}$ is \mult free. As $\nu\in\Irr (W(G_1)\times W_2(\delta_2))$,
we can write $\nu=\nu_1\otimes \nu_2$, where $\nu_1\in\Irr W(G_1)$, $\nu_2\in\Irr W_2(\delta_2)$. 
                                             
\bl{ju7} $\nu^{W(\delta)}=\nu_1^{W(X_1)}\otimes \nu_2$ is a \mult free character of $W(\delta)$. 
\el  We have $\nu^{W(\delta)}=\nu_1^{W(X_1)}\otimes \nu_2^{W_2(\delta_2)}$. 
As  $\nu_2\in \Irr W_2(\delta_2)$, 
$\nu_2^{W_2(\delta_2)}$ is simply $\nu_2$ itself. So it suffices  to prove  that $\nu_1^{W(X_1)} $ is \mult free. 
Recall that $\nu_1$ is the Harish-Chandra correspondent of $St_{G_1}$ and $G_1=GL_m(q)$ for some $m$.
It is well known that $W(GL_m(q))\cong S_m$ and $\nu_1$ is the sign character of $S_m$. 
Furthermore, $X_1\cong SO_{2m}^+(q)$, so $W(X_1)$ is of type $D_m$. The fact that $\nu_1^{W(X_1)}$ is \mult free 
is well known and has been already discussed in the proof of Lemma \ref{aa6}.

\begin{propo}\label{rp1}  Let $\lam_s$ be as in Lemma $ \ref{tt8}$. Then the character $\lam_s^{\# G}$ is \mult free. 
\end{propo}

{\bf Proofs of Theorems} \ref{thm1} and \ref{thm2}. Recall that $St_H|_G$ is denoted by $St^+_G$ throughout the paper. Furthermore, it follows from Deligne-Lusztig theory that $St^*_G=\sum_s St^+_s$, where $s$ runs over semisimple elements of $G^*$, and $St^+_s$ is either zero or the sum of all \ir constituents of $St^+_G$ that belong to the Lusztig (rational) series ${\mathcal E}_s$. Let $V^*$ be the natural module for $G^*$ and $V^*_s $ the 1-eigenspace of $s$ on $V^*$. Then $St^+_s=0$ \ii $V_s^*$ is non-zero and  of Witt defect    1 (Proposition \ref{pp5}). If $V^*_s=0$ then $St^+_s$ is the regular character of 
${\mathcal E}_s$ (Proposition \ref{pp4}). Suppose that $St^+_s\neq 0$ and $V^*_s\neq 0$. Then the character $St^+_s$ is described by Proposition \ref{mr9}. If $\dim V$ is odd then $St^+_s$ is Harish-Chandra induced from a regular character
of a suitable Levi subgroup of $G$, and it is \mult free by Proposition \ref{rp1}. This implies Theorem \ref{thm1} and Theorem \ref{thm2} in this case. Suppose that $\dim V$ is even. Then $St^+_s$ is the sum of at most two characters, each is Harish-Chandra induced from a regular character
of a suitable Levi subgroup of $G$. So again the results follow from Proposition \ref{rp1}. The nature of the regular characters of the Levi subgroups in question is explained in Proposition \ref{tt8}. 

\med Examples. (1) Let $G=SO_3(q)$. Then $St^+_G=1_G+\Gamma _G$, where $\Gamma _G$ is the Gelfand-Graev character of $G$. Moreover, the linear character of $G$ of order $2$ is the only \ir character of $G$ that is not a constituent of $St^+_G$. 

Indeed, in this case $G^*=Sp_2(q)=SL_2(q)$, and $\Id$ is the only semisimple element of $G^*$ which has \ei 1. 
 If $s=\Id$ then $L^*=T^*$, where $T^*$ is  a split maximal torus of $G^*$. Then $St_L=1_L=1_T$ so $St^+_1=1_T^{\#G}=1_B^G$, where $B$ is a Borel subgroup of $G$ containing $T$. It is well known that $1_B^G=1_G+St_G$. If $s\neq \Id$ then $s$ does not have \ei 1, and hence $St^+_s$ is a regular character in ${\mathcal E}_s$ by Proposition \ref{pp4}. Let $\Gamma_G$ be the sum of all regular characters of $G$ (as the center of 
 ${\mathbf G}$ is  trivial and hence connected, $\Gamma_G$  is the Gelfand-Graev character of $G$). Recall that $St_G$ is the  regular character of ${\mathcal E}_1$. This implies the first claim. We omit the proof of the second claim. 

Let $D$ denote the Curtis duality operation. Then $\om_G=D(1_G)+D(\Gamma _G)=St_G+D(\Gamma _G)$.
It is known that $D(\Gamma _G)$ vanishes at all semisimple elements of $G$
\cite[14.33]{DM}. Therefore, $w_G$ coincides with $St_G$
at the semisimple elements of $G$. \itf that $St^+_G=\om_G \cdot St_G =St^2_G$.

\med
(2) $G=SO^-_4(q)$, $q$ odd. Then $G\cong PSL_2(q^2)\times  \lan-\Id\ran$ (as $\Omega^-_4(q)\cong PSL_2(q^2)$). 

By Proposition  \ref{pp5}, $St^+_s\neq 0$ unless $s=1$ or $V_s^*$ is an anisotropic  subspace of $V^*$. Suppose that $St^+_s\neq 0$. If $V_s^*=0$ then $St^+_s$ is an \ir regular character of $G$ (Proposition \ref{pp4}). If $V_s^*\neq 0$ then $V_s^*$ is of Witt defect 0 and $\dim V_s^*=2$.   
We can assume that $ V_s^*$ is the same for all such $s$. 
By Lemma \ref{no9}, ${\mathcal L}$ consists of a single $C_{G^*}(s)$-orbit. If $L\in {\mathcal L}$,  then  $L\cong L^*\cong GL_1(q)\times SO^-_2(q)$ is a maximal torus of $G$ of order $q^2-1$. As the Steinberg character of $GL_1(q)$ is $1_{GL_1(q)}$, by Lemma \ref{pr8} we have $St_s^+=2\cdot (1_{GL_1(q)}\otimes \rho_{s_2})^{\#G}$ if $s_2\neq -\Id$, 
otherwise
$St_s^+=(1_{GL_1(q)}\otimes \rho_{s_2})^{\#G}$. (Here $s_2\neq 1_{SO^-_2(q)}$.)  
The \ir constituents of $(1_{GL_1(q)}\otimes \rho_{s_2})^{\#G}$ are in bijection with those of 
$1_{GL_1(q)} ^{\#SO^+_2(q)}$, and this character is trivial as $SO^+_2(q)=GL_1(q)$.  So each  character $\chi_{s_2}:=(1_{GL_1(q)}\otimes \rho_{s_2})^{\#G}$ is \ir of degree $q^2+1$. One observes that there are $(q-1)/2$  $G^*$-conjugacy classes of $s$ such that $V_s^*$ is of Witt defect 0,  $\dim V_s^*=2$ and   $s_2\neq -\Id$ and  a single class with $s_2=-\Id$.  
Let $\Gamma _G$ be the sum of all regular characters of $G$. (As ${\mathbf G}$ is not with connected center,
there are two Gelfand-Graev characters $\Gamma _1,\Gamma _2$ of $G$, and each regular character of $G$ is a constituent of $\Gamma _1$ or $\Gamma _2$.)
 Therefore, $St_G^+=\Gamma_G+\sum_{s_2\neq \pm \Id}\chi_{s_2}-\sum_i\chi_s$, where  $\chi_s$ runs over the regular characters of $G$ that does  not occur in $St_G^+$. This means that $s$ runs over representatives of the conjugacy classes
$s$ of $G^*$   such that either $s=1$ or $V_s^*$ is  anisotropic. The number of classes with $V_s^*$ anisotropic equals $(q-1)/2$.
So $St_G^+=\Gamma_G+\sum_{s_2\neq \pm \Id}\chi_{s_2}- St_G-\sum_s\chi_s$, where $s$ is such that $\dim V_s^*=2$ and $V_s^*$ is anisotropic (chosen single  in each semisimple conjugacy class).

\med
Acknowledgement. 
This work was partially supported by the Hausdorff Research Institute for Mathematics (Bonn University, Germany) during the Trimester program "Universality and Homogeneity" (Autumn 2013).

\med\noindent 
National Academy of Sciences, Minsk, Belarus and University of East Anglia, Norwich, UK 

\med
e-mail: alexandre.zalesski@gmail.com

\end{document}